\documentclass[12pt,a4paper]{article}
\usepackage{amssymb,amsthm,amsmath,supertabular,curves}

\theoremstyle{plain}
\newtheorem{Theorem}{Theorem}[section]
\newtheorem{Lemma}[Theorem]{Lemma}

\newtheorem{Claim}[Theorem]{Claim}
\newtheorem{Proposition}[Theorem]{Proposition}
\newtheorem{Corollary}[Theorem]{Corollary}

\newtheorem*{Theorem*}{Theorem}
\newtheorem*{Lemma*}{Lemma}

\theoremstyle{definition}
\newtheorem*{Definition*}{Definition}
\newtheorem*{Remark*}{Remark}
\newtheorem*{Example*}{Example}
\newtheorem{Definition}[Theorem]{Definition}
\newtheorem{Remark}[Theorem]{Remark}

\newcommand{\divides}{\mathbin{\mid}}

\newcommand{\gen}[1]{\mathopen{<}#1\mathclose{>}}
\newcommand{\abs}[1]{\lvert#1\rvert}
\newcommand{\mat}[4]{\bigl(\begin{smallmatrix}#1&#2\\#3&#4
                    \end{smallmatrix}\bigr)}

\DeclareMathOperator{\Aut}{Aut}
\DeclareMathOperator{\Out}{Out}

\DeclareMathOperator{\Gal}{Gal}
\DeclareMathOperator{\Isom}{Isom}

\DeclareMathOperator{\lcm}{lcm}

\DeclareMathOperator{\ind}{ind}

\newcommand{\CCC}{\mathbb C}
\newcommand{\FFF}{\mathbb F}            
\newcommand{\NNN}{\mathbb N}
\newcommand{\QQQ}{\mathbb Q}
\newcommand{\RRR}{\mathbb R}
\newcommand{\ZZZ}{\mathbb Z}
\newcommand{\PP}{\mathbb P^1}           

\newcommand{\hk}{\hat k}

\newcommand{\ok}{{\mathcal O}_k}        
\newcommand{\cC}{\mathcal C}

\newcommand{\cE}{\mathcal E}

\newcommand{\frakp}{\mathfrak p}
\newcommand{\frakP}{\mathfrak P}

\newcommand{\bn}[1]{\mathbf #1}         
\newcommand{\s}[1]{\sigma_#1}

\DeclareMathOperator{\End}{End}

\DeclareMathOperator{\gL}{\Gamma L}
\DeclareMathOperator{\GL}{GL}
\DeclareMathOperator{\sL}{\Sigma L}
\DeclareMathOperator{\SL}{SL}

\DeclareMathOperator{\GU}{GU}

\DeclareMathOperator{\PgL}{P\Gamma L}
\DeclareMathOperator{\PGL}{PGL}
\DeclareMathOperator{\PsL}{P\Sigma L}
\DeclareMathOperator{\PSL}{PSL}

\DeclareMathOperator{\PGU}{PGU}
\DeclareMathOperator{\PSU}{PSU}

\DeclareMathOperator{\PSp}{PSp}

\DeclareMathOperator{\AGL}{AGL}
\DeclareMathOperator{\ASL}{ASL}
\DeclareMathOperator{\AgL}{A\Gamma L}
\DeclareMathOperator{\AsL}{A\Sigma L}

\newcommand{\M}[2]{\mathrm{M}_{#1#2}}             

\DeclareMathOperator{\alt}{\mathcal{A}}
\DeclareMathOperator{\sym}{\mathcal{S}}
\DeclareMathOperator{\Spl}{Sp}
\DeclareMathOperator{\Red}{Red}

\newcommand{\ord}[1]{\text{ord}(#1)}
\newcommand{\cX}{\mathcal X}
\newcommand{\Si}{S}

\allowdisplaybreaks[1]

\begin{document}
\title{Permutation Groups with a Cyclic Two-Orbits Subgroup and
  Monodromy Groups of Siegel Functions} \author{Peter M\"uller}
\maketitle

\begin{abstract}
  We classify the finite primitive permutation groups which have a
  cyclic subgroup with two orbits. This extends classical topics in
  permutation group theory, and has arithmetic consequences. By a
  theorem of C.~L.~Siegel, affine algebraic curves with infinitely
  many integral points are parametrized by rational functions whose
  monodromy groups have this property. We classify the possibilities
  of these monodromy groups, and give an application to Hilbert's
  irreducibility theorem.
\end{abstract}

\section{Introduction}

The aim of this paper is twofold: First, it provides the group
theoretic and arithmetic classification results to obtain sharpenings
of Hilbert's irreducibility theorem, like the following:

\begin{Theorem}
  Let $f(t,X)\in\QQQ(t)[X]$ be an irreducible polynomial with Galois
  group $G$, where $G$ is a simple group not isomorphic to an
  alternating group or $C_2$. Then $\Gal(f(\bar t,X)/\QQQ)=G$ for all
  but finitely many specializations $\bar t\in\ZZZ$.
\end{Theorem}

More results of this kind are given in Section \ref{S:App}.

The second purpose is to obtain a group theoretic classification
result, namely to determine those primitive permutation groups which
contain a cyclic subgroup with only two orbits. This classification
completes and generalizes previous results about permutation groups.
The list of possibilities is quite long and involved, we give it in
Section \ref{S:Results}.

In the following we explain how these seemingly unrelated topics are
connected.

Let $k$ be a field of characteristic $0$, and $f(X)\in k[X]$ be a
functionally indecomposable polynomial. Let $t$ be a transcendental
over $k$, and $A$ the Galois group of $f(X)-t$ over $k(t)$, considered
as a permutation group in the action on the roots of $f(X)-t$. It is
easy to see that this group action is primitive. Furthermore, $A$
contains a transitive cyclic subgroup $I$, the inertia group of a
place of the splitting field of $f(X)-t$ over the place
$t\mapsto\infty$. By classical theorems of Schur and Burnside, it is
known that a primitive permutation group $A$ with a transitive cyclic
subgroup is either doubly transitive, or a subgroup of the affine
group $\AGL_1(p)$ for a prime $p$. Using the classification of the
doubly transitive groups, which rests on the knowledge of the finite
simple groups, these groups have been determined, see \cite[Theorem
4.1]{Feit:SC}. Furthermore, if $k$ is algebraically closed, then one
can determine the subset of these groups which indeed are Galois
groups as above, see \cite{PM:Mon}, which completes (and corrects)
\cite{Feit:BD}. Besides the alternating, symmetric, cyclic and
dihedral groups only finitely many cases arise. For arithmetic
applications of these results see \cite{Fried:SC}, \cite{PM:Kron},
\cite{MV:Dav}.

A more general situation arises if one considers Hilbert's
irreducibility theorem over a number field $k$. Let $\ok$ be the
integers of $k$, and $f(t,X)\in k(t)[X]$ an irreducible polynomial.
Then $f(\bar t,X)$ is irreducible over $k$ for infinitely many
integral specializations $\bar t\in\ok$. Using Siegel's deep theorem
about algebraic curves with infinitely many integral points, one can,
to some extent, describe the set $\Red_f(\ok)$ of those
specializations $\bar t\in\ok$ such that $f(\bar t,X)$ is defined and
reducible. There are finitely many rational functions $g_i(Z)\in
k(Z)$, such that $\abs{g_i(k)\cap\ok}=\infty$ and $\Red_f(\ok)$
differs by finitely many elements from the union of the sets
$g_i(k)\cap\ok$. Thus, in order to get refined versions of Hilbert's
irreducibility theorem, one has to get information about rational
functions $g(Z)$ which assume infinitely many integral values on $k$,
see \cite{PM:hitzt}. We call such a function a \emph{Siegel function}.
A theorem of Siegel shows that a Siegel function has at most two poles
on the Riemann sphere $\PP(\CCC)$. Suppose that $g(Z)$ is functionally
indecomposable, and let $A$ be the Galois group of the numerator of
$g(Z)-t$ over $k(t)$. It follows that $A$ is primitive on the roots of
$g(Z)-t$, and the information about the poles of $g(Z)$ yields a
cyclic subgroup $I$ of $A$ such that $I$ has at most two orbits.

This generalizes the situation coming from the polynomials, where $I$
has just one orbit. In the two-orbit situation not much was known.
There is a result by Wielandt \cite{Wielandt:2p} (see also
\cite[V.31]{Wielandt}), if the degree of $A$ is $2p$ for a prime $p$,
and both orbits of $I$ have length $p$. Then $A$ is either doubly
transitive, or $2p-1$ is a square, and a point stabilizer of $A$ has
three orbits of known lengths. Another case which has been dealt with
is that $I$ has two orbits of relatively prime lengths $i<j$. If
$i>1$, then a classical result by Marggraf \cite[Theorem
13.5]{Wielandt} immediately shows that $A$ contains the alternating
group in natural action. The more difficult case $i=1$ was treated in
\cite{PM:HS}. In Section \ref{S:2} we classify the primitive groups
with a cyclic two-orbit subgroup without any condition on the orbit
lengths.

As in the polynomial case, we again want to know which of these groups
indeed are Galois groups of $g(Z)-t$ over $k(t)$, with $g(Z)$ a Siegel
function as above. This amounts to finding genus $0$ systems (to be
defined later) in normal subgroups of $A$, see Section \ref{SS:g0}.

Finally, for applications to Hilbert's irreducibility theorem over the
rationals, we classify the Galois groups of $g(Z)-t$ over $\QQQ(t)$ of
those rational functions $g(Z)\in\QQQ(Z)$ with
$\abs{g(\QQQ)\cap\ZZZ}=\infty$. Section \ref{S:Rat} is devoted to
that. In contrast to the other mostly group theoretic sections, we
need several arithmetical and computational considerations related to
the regular inverse Galois problem over $\QQQ(t)$.

\tableofcontents

\section{Permutation Groups -- Notations, Definitions, and Elementary
Results}

Here we collect definitions and easy results about finite groups and
finite permutation groups, which are used throughout the work.

\begin{description}
\item[General notation:] For $a,b$ elements of a group $G$ set
  $a^b:=b^{-1}ab$. Furthermore, if $A$ and $B$ are subsets of $G$,
  then $A^b$, $a^B$ and $A^B$ have their obvious meaning. If $H$ is a
  subgroup of $G$, then for a subset $S$ of $G$ let $C_H(S)$ denote
  the centralizer of $S$ in $H$ and $N_H(S)$ denote the normalizer
  $\{h\in H|\;S^h=S\}$ of $S$ in $H$.
  
  If $A,B,\dots$ is a collection of subsets or elements of $G$, then
  we denote by $\gen{A,B,\dots}$ the group generated by these sets and
  elements.
  
  The order of an element $g\in G$ is denoted by $\ord{g}$.
  
\item[Permutation groups:] Let $G$ be a permutation group on a finite
  set $\Omega$. Then $\abs{\Omega}$ is the
  \emph{degree} of
  $G$. We use the exponential notation $\omega^g$ to denote the image
  of $\omega\in\Omega$ under $g\in G$. The stabilizer of $\omega$ in
  $G$ is denoted by $G_\omega$. If $G$ is transitive and $G_\omega$ is
  the identity subgroup, then $G$ is called
  \emph{regular}.
  
  The number of fixed points of $g$ on $\Omega$ will be denoted by
  $\chi(g)$.
  
  Let $G$ be transitive on $\Omega$ of degree $\ge2$, and let
  $G_\omega$ be the stabilizer of $\omega\in\Omega$. Then the number
  of orbits of $G_\omega$ on $\Omega$ is the
  \emph{rank} of $G$. In particular,
  the rank is always $\ge2$, and exactly $2$ if and only if the group
  is doubly transitive. The
  \emph{subdegrees} of $G$ are
  defined as the orbit lengths of $G_\omega$ on $\Omega$.
  
  Let $G$ be transitive on $\Omega$, and let $\Delta$ be a nontrivial
  subset of $\Omega$. Set $S:=\{\Delta^g|\;g\in G\}$. We say that
  $\Delta$ is a \emph{block} of $G$ if $S$ is a
  partition of $\Omega$.  If this is the case, then $S$ is called a
  \emph{block system} of $G$. A block (or block system) is called trivial if
  $\abs{\Delta}=1$ or $\Delta=\Omega$. If each block system of $G$ is
  trivial, then $G$ is called
  \emph{primitive}. Primitivity of $G$ is equivalent to
  maximality of $G_\omega$ in $G$. Note that the orbits of a normal
  subgroup $N$ of $G$ constitute a block system, thus a normal
  subgroup of a primitive permutation group is either trivial or
  transitive.
  
\item[Specific groups:] We denote by $C_n$ and
  $D_n$ the cyclic and dihedral group of order $n$
  and $2n$, respectively. If not otherwise said, then $C_n$ and $D_n$
  are regarded as permutation groups in their natural degree $n$
  action.  The alternating and symmetric group on $n$ letters is
  denoted by $\alt_n$ and
  $\sym_n$, respectively.
  
  We write $\sym(M)$ for the symmetric group
  on a set $M$.
  
  Let $m\ge1$ be an integer, and $q$ be a power of the prime $p$. Let
  $\FFF_q$ be the field with $q$ elements. We denote by $\GL_m(q)$ (or
  sometimes $\GL_m(\FFF_q)$) the general linear group of $\FFF_q^m$,
  and by $\SL_m(q)$ the special linear group. Regard these groups as
  acting on $\FFF_q^m$. The group $\Gamma:=\Gal(\FFF_q/\FFF_p)$ acts
  componentwise on $\FFF_q^m$. This action of $\Gamma$ normalizes the
  actions of $\GL_m(q)$ and $\SL_m(q)$. We use the following symbols
  for the corresponding semidirect products:
  $\gL_m(q):=\gen{\GL_m(q),\Gamma}=\GL_m(q)\rtimes\Gamma$,
  $\sL_m(q):=\gen{\SL_m(q),\Gamma}=\SL_m(q)\rtimes\Gamma$.
  
  Note that if $q=p^e$, then we have the natural inclusion
  $\gL_m(q)\le\GL_{me}(p)$.
  
  Let $G$ be a subgroup of $\gL_m(q)$, and denote by $N$ the action of
  $\FFF_q^m$ on itself by translation. Then $G$ normalizes the action
  of $N$. If $G=\GL_m(q)$, $\SL_m(q)$, $\gL_m(q)$, or $\sL_m(q)$, then
  denote the semidirect product of $G$ with $N$ by
  $\AGL_m(q)$,
  $\ASL_m(q)$,
  $\AgL_m(q)$, or
  $\AsL_m(q)$, respectively. A group $A$ with
  $N\le A\le\AgL_m(q)$ is called an \emph{affine permutation
    group}.
  
  Let $G\le\gL_m(q)$ act naturally on $V:=\FFF_q^m$. We denote by
  $\PP(V)$ the set of one--dimensional subspaces of $V$. As $G$
  permutes the elements in $\PP(V)$, we get an (in general not
  faithful) action of $G$ on $\PP(V)$. The induced faithful
  permutation group on $\PP(V)$ is named by prefixing a $P$ in front
  of the group name, so we get the groups $\PGL_m(q)$, $\PSL_m(q)$,
  $\PgL_m(q)$, or $\PsL_m(q)$, respectively.
  
  The group $\GL_m(q)$ contains, up to conjugacy, a unique cyclic
  subgroup which permutes regularly the non--zero vectors of
  $\FFF_q^m$.  This group, and also its homomorphic image in
  $\PGL_m(q)$, is usually called \emph{Singer group}. Existence of
  this group follows from the regular representation of the
  multiplicative group of $\FFF_{q^m}$ on $\FFF_{q^m}\cong\FFF_q^m$,
  uniqueness follows for example from Schur's Lemma and the
  Skolem--Noether Theorem (or by the Lang's Theorem).
  
  For $n\in\{11,12,22,23,24\}$ we denote by $\mathrm{M}_n$ the five
  \emph{Mathieu groups} of degree $n$, and let $\M10$ be a point
  stabilizer of $\M11$ in the transitive action on $10$ points.
\end{description}

\subsection{The Aschbacher--O'Nan--Scott Theorem}

The Aschbacher--O'Nan--Scott Theorem makes a rough distinction between
several possible types of actions of a primitive permutation group.
This theorem had first been announced by O'Nan and Scott on the Santa
Cruz Conference on Finite groups in 1979, see \cite{Scott:AONS}. In
their statement a case was missing, and the same omission appears in
\cite{Cameron}. To our knowledge, the first complete version is in
\cite{AschScott}. A very concise and readable proof is given in
\cite{LPS}, see also \cite{DixMort}.

Let $A$ be a primitive permutation group of degree $n$ on
$\Omega$. Then one of the following actions occurs:
\begin{description}
\item[Affine action.] We can identify $\Omega$ with a vector space
  $\FFF_p^m$, and $\FFF_p^m\le A\le\AGL_m(p)$ is an affine group as
  described above.
\item[Regular normal subgroup action.] $A$ has a non--abelian normal
  subgroup which acts regularly on $\Omega$. (There are finer
  distinctions in this case, see \cite{LPS}, but we don't need that
  extra information.)
\item[Diagonal action.] $A$ has a unique minimal normal subgroup of the 
  form $N=S_1\times S_2\times\dots\times S_t$, where the $S_i$ are
  pairwise isomorphic non--abelian simple groups, and the point
  stabilizer $N_\omega$ is a diagonal subgroup of $N$.
\item[Product action.] We can write
  $\Omega=\Delta\times\Delta\dots\times\Delta$ with $t\ge2$ factors,
  and $A$ is a subgroup of the wreath product
  $\sym(\Delta)\wr\sym_t=\sym(\Delta)^t\rtimes\sym_t$ in the natural
  product action on this cartesian product. In such a case, we will
  say that $A$ preserves a product structure.
\item[Almost simple action.] There is $S\le A\le\Aut(S)$ for a simple
  non--abelian group $S$. In this case, $S$ cannot act regularly.
\end{description}

\section{Elements with two cycles}\label{S:2}

\subsection{Previous Results}
Our goal is the classification of those primitive permutation groups
$A$ of degree $n$ which contain an element $\sigma$ with at most two
cycles. If $\sigma$ actually is an $n$--cycle, the result is a
well--known consequence of the classification of doubly transitive
permutation groups.

\begin{Theorem}[{Feit \cite[4.1]{Feit:SC}}]\label{T:Feit:nZykel}
  Let $A$ be a primitive permutation group of degree $n$ which
  contains an $n$--cycle. Then one of the following holds.
  \begin{itemize}
  \item[(a)] $A\le\AGL_1(p)$, $n=p$ a prime; or
  \item[(b)] $A=\alt_n$ or $\sym_n$; or
  \item[(c)] $\PSL_k(q)\le A\le\PgL_k(q)$, $k\ge2$, $q$ a prime power, 
    $A$ acting naturally on the projective space with
    $n=(q^k-1)/(q-1)$ points; or
  \item[(d)] $n=11$, $A=\PSL_2(11)$ or $\M11$; or
  \item[(e)] $n=23$, $A=\M23$.
  \end{itemize}
\end{Theorem}

If $\sigma$ has two cycles of coprime length, say $k$ and $l=n-k$ with
$k\le l$, then it follows immediately from Marggraf's theorem
\cite[Theorem 13.5]{Wielandt}, applied to the subgroup generated by
$\sigma^l$, that $\alt_n\le A$ unless $k=1$. The critical case thus is
$k=1$. We quote the classification result \cite[6.2]{PM:HS}.

\begin{Theorem}\label{T:PM:(n-1)}
  Let $A$ be a primitive permutation group of degree $n$ which
  contains an element with exactly two cycles, of coprime lengths
  $k\le l$. Assume that $\alt_n\not\le A$. Then $k=1$, and one of the
  following holds.
  \begin{itemize}
  \item[(a)] $n=q^m$ for a prime power $q$, $\AGL_m(q)\le
    A\le\AgL_m(q)$; or
  \item[(b)] $n=p+1$, $\PSL_2(p)\le A\le\PGL_2(p)$, $p\ge5$ a prime;
    or
  \item[(c)] $n=12$, $A=\M11$ or $\M12$; or
  \item[(d)] $n=24$, $A=\M24$.
  \end{itemize}  
\end{Theorem}

In the remainder of this chapter, we deal with the case where $k$ and
$l$ are not necessarily coprime. The assumptions in the Theorems
\ref{T:Feit:nZykel} and \ref{T:PM:(n-1)} quickly give that $A$ is
doubly transitive (or $A\le\AGL_1(p)$, a trivial case) --- this is
clear under existence of an $(n-1)$--cycle, and follows from Theorems
of Schur \cite[25.3]{Wielandt} and
Burnside \cite[XII.10.8]{Huppert3} under the
presence of an $n$--cycle. So one basically has to check the list of
doubly transitive groups.

In the general case however, $A$ no longer need to be doubly
transitive.  Excluding the case $A\le\AGL_1(p)$, we will obtain as a
corollary of our classification that $A$ has permutation rank $\le3$,
though I do not see how to obtain that directly. I know only two
results in the literature where this has been shown under certain
restrictions. The first one is by Wielandt
\cite[Theorem 31.2]{Wielandt}, \cite{Wielandt:2p}, under the
assumption that $k=l$, and $k$ is a prime, and the other one is by
Scott, see the announcement of the never
published proof in \cite{Scott:n2n}. In Scott's announcement, however,
there are several specific assumptions on $A$. First $k=l$, and $A$
has to have a doubly transitive action of degree $k$, such that the
point--stabilizer in this action is intransitive in the original
action, but that the element with the two cycles of length $k$ in the
original action is a $k$--cycle in the degree $k$ action.

\subsection{Classification Result}\label{S:Results}
Recall that if $t$ is a divisor of $m$, then we have $\gL_{m/t}(p^t)$
naturally embedded in $\GL_m(p)$. We use this remark in the main
result of this chapter.

\begin{Theorem}\label{T:kl}
  Let $A$ be a primitive permutation group of degree $n$ which
  contains an element with exactly two cycles of lengths $k$ and
  $n-k\ge k$. Then one of the following holds, where $A_1$ denotes the
  stabilizer of a point.
  \begin{enumerate}
  \item\label{T:kl:I} (Affine action) $A\le\AGL_m(p)$ is an affine
    permutation group, where $n=p^m$ and $p$ is a prime number.
    Furthermore, one of the following holds.
\begin{enumerate}
\item\label{T:kl:I:a} $n=p^m$, $k=1$, and $\GL_{m/t}(p^t)\le
  A_1\le\gL_{m/t}(p^t)$ for a divisor $t$ of $m$;
\item\label{T:kl:I:b} $n=p^m$, $k=p$, and $A_1=\GL_m(p)$;
\item\label{T:kl:I:c} $n=p^2$, $k=p$, and $A_1<\GL_2(p)$ is the group
  of monomial matrices (here $p>2$);
\item\label{T:kl:I:d} $n=2^m$, $k=4$, and $A_1=\GL_m(2)$;
\item (Sporadic affine cases)
\begin{enumerate}
\item $n=4$, $k=2$, and $A_1=\GL_1(4)$ (so $A=\alt_4$);
\item\label{T:kl:I:d:ii} $n=8$, $k=2$, and $A_1=\gL_1(8)$;
\item $n=9$, $k=3$, and $A_1=\gL_1(9)$;
\item $n=16$, $k=8$, and $[\gL_1(16):A_1]=3$ or $A_1\in\{\gL_1(16),
  (C_3\times C_3)\rtimes C_4, \sL_2(4), \gL_2(4), \alt_6, \GL_4(2)\}$;
\item $n=16$, $k=4$ or $8$, and $A_1\in\{(\sym_3\times\sym_3)\rtimes
  C_2, \sym_5, \sym_6\}$;
\item $n=16$, $k=2$ or $8$, and $A_1=\alt_7<\GL_4(2)$;
\item $n=25$, $k=5$, and $[\GL_2(5):A_1]=5$;
\end{enumerate}
\end{enumerate}
\item\label{T:kl:II} (Product action) One of the following holds.
  \begin{enumerate}
  \item $n=r^2$ with $1<r\in\NNN$, $k=ra$ with $\gcd(r,a)=1$,
    $A=(\sym_r\times\sym_r)\rtimes C_2$, and
    $A_1=(\sym_{r-1}\times\sym_{r-1})\rtimes C_2$;
  \item $n=(p+1)^2$ with $p\ge5$ prime, $k=p+1$,
    $A=(\PGL_2(p)\times\PGL_2(p))\rtimes C_2$, and
    $A_1=(\AGL_1(p)\times\AGL_1(p))\rtimes C_2$.
  \end{enumerate}
\item\label{T:kl:III} (Almost simple action) $\Si\le A\le\Aut(\Si)$
  for a simple, non--abelian group $\Si$, and one of the following
  holds.
\begin{enumerate}
\item\label{T:kl:III:An} $n\ge5$, $\alt_n\le A\le\sym_n$ in natural
  action;
\item\label{T:kl:III:A5(10)} $n=10$, $k=5$, and $\alt_5\le A\le\sym_5$
  in the action on the $2$--sets of $\{1,2,3,4,5\}$;
\item\label{T:kl:III:PSL2(p)} $n=p+1$, $k=1$, and $\PSL_2(p)\le
  A\le\PGL_2(p)$ for a prime $p$;
\item\label{T:kl:III:Singer/2} $n=(q^m-1)/(q-1)$, $k=n/2$, and
  $\PSL_m(q)\le A\le\PgL_m(q)$ for an odd prime power $q$ and $m\ge2$
  even;
\item\label{T:kl:III:M10} $n=19$, $k=2$, and $\M10\le A\le\PgL_2(9)$;
\item\label{T:kl:III:PsL3(4)} $n=21$, $k=7$, and $\PsL_3(4)\le
  A\le\PgL_3(4)$;
\item\label{T:kl:III:M11(12)} $n=12$, $k=1$ or $4$, and $A=\M11$ in its
  action on $12$ points;
\item\label{T:kl:III:M12} $n=12$, $k=1$, $2$, $4$, or $6$, and
  $A=\M12$;
\item\label{T:kl:III:M22} $n=22$, $k=11$, and $\M22\le
  A\le\Aut(\M22)=\M22\rtimes C_2$;
\item\label{T:kl:III:M24} $n=24$, $k=1$, $3$, or $12$, and $A=\M24$.
\end{enumerate}
\end{enumerate}
\end{Theorem}

The proof of this theorem is given in the following sections. We
distinguish the various cases of the Aschbacher-O'Nan-Scott theorem,
because the cases require quite different methods. The almost simple
groups comprise the most complex case. We split this case further up
into the subcases where the simple normal subgroup $\Si$ is
alternating, sporadic, a classical Lie type group, or an exceptional
Lie type group.

We note an interesting consequence of the previous theorem which
generalizes several classical results on permutations groups. It would
be very interesting to have a direct proof, without appealing to the
classification of the finite simple groups.

\begin{Corollary} Let $A$ be a primitive permutation group
  which contains an element with exactly two cycles. Then $A$ has rank
  at most $3$.\end{Corollary}

\subsection{Affine Action}

We need the following well--known

\begin{Lemma}\label{L:unip:ord}
  Let $K$ be a field of positive characteristic $p$, and
  $\sigma\in\GL_m(K)$ of order $p^b\ge p$. Then $p^{b-1}\le m-1$. In
  particular, $\ord{\sigma}\le p(m-1)$ if $m\ge2$.
\end{Lemma}

\begin{proof} $1$ is the only eigenvalue of $\sigma$, therefore
  $\sigma$ is conjugate to an upper triangular matrix with $1$'s on
  the diagonal.  So $\sigma-\bn1$ is nilpotent. Now
  $(\sigma-\bn1)^{p^{b-1}}=\sigma^{p^{b-1}}-\bn1\ne0$, thus
  $p^{b-1}<m$, and the claim follows.\end{proof}

We note the easy consequence

\begin{Lemma}\label{L:Aff:p}
  Let $A$ be an affine permutation group of degree $p^m$. Let $A$
  contain an element of order $p^r$ for $r\in\NNN$. Then $p^{r-1}\le
  m$. In particular, if $r=m$, then $m\le2$, and $m=1$ for $p>2$.
\end{Lemma}

\begin{proof} Without loss $A=\AGL_m(p)$. We use the well--known
embedding of $A$ in $\GL_{m+1}(p)$: Let $g\in\GL_m(p)$, $\mathbf{v}\in
N$. Then define the action of $g\mathbf{v}\in A$ on the vector space
$N\times\FFF_p$ via
$(\mathbf{w},w_{m+1})^{g\mathbf{v}}:=
(\mathbf{w}g+w_{m+1}\mathbf{v},w_{m+1})$ for $\mathbf{w}\in N$,
$w_{m+1}\in\FFF_p$.

This way we obtain an element $\sigma$ of order $p^r$ in
$\GL_{m+1}(p)$. The claim follows from Lemma
\ref{L:unip:ord}.\end{proof}

\begin{Lemma}\label{L:Impr:rk3}
  Let $H$ be a rank $3$ permutation group with subdegrees $1\le u\le
  v$. Suppose that $H$ is imprimitive. Then $1+u$ divides $v$.
\end{Lemma}

\begin{proof} Let $\Delta$ be a nontrivial block, and
  $\delta\in\Delta$. Let $\epsilon\ne\delta$ be in $\Delta$. Then
  $\Delta$ contains the orbit $\epsilon^{H_\delta}$. Also, $\Delta$
  does not meet a point from an $H_\delta$--orbit different from
  $\{\delta\}$ and $\epsilon^{H_\delta}$, for then $\Delta$ were the
  full set $H$ is acting on. Thus
  $\Delta=\{\delta\}\cup\epsilon^{H_\delta}$. But
\[
\abs{\epsilon^{H_\delta}}=\abs{\Delta}-1\le(1+u+v)/2-1<v,
\]
so the orbit $\epsilon^{H_\delta}$ has length $u$, and $1+u$ divides
$v$ because the orbit of size $v$ is a union of conjugates of
$\Delta$.\end{proof}

We need to know the doubly transitive permutation subgroups of the
collineation group of a projective linear
space.

\begin{Proposition}[{Cameron,
    Kantor \cite[Theorem I]{CamKan:2tranproj}}]\label{P:PGL:2tr} Let
  $m\ge3$, $p$ be a prime, and $H\le\GL_m(p)$ be acting doubly
  transitively on the lines of $\FFF_p^m$. Then $\SL_m(p)\le H$ or
  $H=\alt_7<\SL_4(2)$.
\end{Proposition}

Also, the primitive rank $3$ permutation subgroups of even order of
the collineation group of a projective linear space have been
classified by Perin in an unpublished thesis. We use a result of
Cameron and Kantor which extends this result. The odd order case,
which is not given there, is easily handled by a result of Huppert.

\begin{Proposition}\label{P:PGL:rk3}
  Let $m\ge3$, $p$ be a prime, and $H\le\GL_m(p)$ be acting
  primitively of rank $3$ on the lines of $\FFF_p^m$. Then
  $\Spl_m(p)\trianglelefteq H$, or $p=2$, $m=3$, $H=\gL_1(8)$, and the
  subdegrees are $1$, $3$, $3$.
\end{Proposition}

\begin{proof} If $H$ has even order see the remarks preceding
  \cite[Prop.\ 8.5]{CamKan:2tranproj}. So assume that $H$ has odd
  order. Then $H$ is solvable, so also $\tilde H=\FFF_p^\star H$ is
  solvable. Furthermore, $\tilde H$ is transitive on $\FFF_p^m$. These
  groups have been classified by Huppert
  \cite[XII.7.3]{Huppert3}. Either $\tilde H\le\gL_1(p^m)$, or
  $p^m=3^4$. (The other exceptional cases in Huppert's classification
  have $m=2$.) Let us look at the action of
  $\gL_1(p^m)=\FFF_{p^m}^\star\rtimes\Aut(\FFF_{p^m})$ on
  $\FFF_{p^m}^\star/\FFF_p^\star$. The stabilizer of the set
  $\FFF_p^\star$ is just $\FFF_p^\star\rtimes\Aut(\FFF_{p^m})$. So the
  orbit lengths of the point stabilizer on the projective space are at
  most $m$. Therefore we get that the rank is at least
\[
1+\frac{(p^m-1)/(p-1)-1}{m}>1+(p^{m-1}-1)/m.
\]
Thus $p^{m-1}-1<2m$, hence $p=2$ and $m=3$ or $4$. If $m=4$, then $H$
has no contribution coming from $\Aut(\FFF_{16})$ by the assumption of
odd order, hence the point stabilizer is even trivial. The case $m=3$
indeed does occur.

Now suppose $p=3$, $m=4$. The transitive soluble subgroups of
$\GL_4(3)$ are available for instance via GAP \cite{GAP}, and one
immediately checks that none of them induces a primitive group on the
lines of $\FFF_3^4$.\end{proof}

In order to handle the case $m=p^2$, we need the following

\begin{Lemma}\label{L:PGL2:irr}
  Let $p$ be a prime, and let $H\le\GL_2(p)$ act irreducibly on
  $\FFF_p^2$. Let $\omega$ be a generator of the multiplicative group
  of $\FFF_p$, and suppose that $\tau:=\mat{1}{0}{0}{\omega}\in H$.
  Then one of the following holds.
\begin{itemize}
\item[(a)] $H=\GL_2(p)$.
\item[(b)] $H$ is the group of monomial matrices.
\item[(c)] $p=5$, and $[\GL_2(5):H]=5$.
\item[(d)] $p=3$, and $H$ is a Sylow $2$--subgroup of $\GL_2(3)$.
\item[(e)] $p=2$, and $H\cong C_3$.
\end{itemize}
\end{Lemma}

\begin{proof} The cases $p=2$ and $3$ are straightforward. So assume
  $p\ge5$. If $\SL_2(p)\le H$, then $H=\GL_2(p)$ and (a) holds,
  because the determinant of $\tau\in H$ is a generator of
  $\FFF_p^\star$.
  
  So we assume in the following that $H$ does not contain $\SL_2(p)$.
  We first contend that $p$ does not divide the order of $H$. Suppose
  it does. Then $H$ contains a Sylow $p$--subgroup $P$ of $H$. If $P$
  is normal in $H$, then $H$ is conjugate to a group of upper
  triangular matrices, hence not irreducible. Therefore $P$ is not
  normal in $H$, thus $H$ contains at least $1+p$ Sylow $p$--subgroups
  of $\GL_2(p)$ (by Sylow's Theorem). But $\GL_2(p)$ has exactly $p+1$
  Sylow $p$--subgroups, so $H$ contains all the $p+1$ Sylow
  $p$--subgroups of $\GL_2(p)$. But these Sylow $p$--subgroups
  generate $\SL_2(p)$, contrary to our assumption.
  
  Set $C=\gen{\tau}$, and let $S\cong\FFF_p^\star$ be the group of
  scalar matrices. So $CS$ is the group of diagonal matrices. First
  assume that $H$ normalizes $CS$. Then, by irreducibility of $H$,
  some element in $H$ must switch the two eigenspaces of $C$. It
  follows quickly that $H$ is monomial.
  
  So finally suppose that $CS$ is not normalized by $H$. Then there is
  a conjugate $(CS)^h$ with $h\in H$, such that $(CS)\cap(CS)^h=S$. So
  we have $\abs{HS}\ge\abs{(CS)(CS)^h}=(p-1)^3$ and
  $(p-1)^2\divides\abs{HS}$. First note that we cannot have
  $\abs{HS}=(p-1)^3$ simply because $(p-1)^3$ does not divide
  $\abs{\GL_2(p)}=(p-1)^2p(p+1)$. So $\abs{HS}\ge p(p-1)^2$. But again
  equality cannot hold, for we noted already that $p$ does not divide
  $\abs{H}$. So $\abs{HS}\ge(p+1)(p-1)^2$, hence $[\GL_2(p):HS]\le p$.
  But $\PGL_2(p)=\GL_2(p)/S$ acts faithfully on the coset space
  $\PGL_2(p)/(HS/S)$ and has an element of order $p$, hence
  $[\GL_2(p):HS]=p$. A classical theorem of
  Galois \cite[II.8.28]{Huppert1} says that if
  $\PSL_2(p)$ has a transitive permutation representation of degree
  $p$, then $p\le11$.  But one checks that $\GL_2(p)$ does not have a
  subgroup of index $p$ for $p=7$ and $11$, thus $p=5$. So
  $\abs{HS}=96=16\cdot2\cdot3$.  Therefore $CS$ (of order $16$) has a
  proper normalizer in $HS$. By an argument as above, we thus obtain
  an element $h\in H$ which switches the eigenspaces of $C$. So
  $\gen{C,C^h}\le H$ is the group of diagonal matrices, in particular
  $S\le H$. The claim follows.\end{proof}

The following proposition is based on Hering's
\cite{Hering:2tran1}, \cite{Hering:2tran2} classification of
transitive linear groups. Note that his results are incomplete, and
the first complete treatment is given by
Liebeck in \cite[Appendix
1]{Liebeck:affinerank3}.

\begin{Proposition}\label{P:Hering}
  Let $m\ge5$, and $H\le\GL_m(2)$ be irreducible on $V:=\FFF_2^m$.
  Suppose there is an element $\tau\in H$ which is the the identity on
  a $2$--dimensional subspace $U$ of $V$, and cyclically permutes the
  nonzero elements of a complement $W$ of $U$ in $V$. Then either
  $H=\GL_m(2)$, or $m$ is even and $\GL_{m/2}(4)\le H\le\gL_{m/2}(4)$.
\end{Proposition}

\begin{proof} For a subspace $X$ of $V$ set $X^\sharp:=X\setminus\{0\}$. We
  first want to show that $H$ is transitive on $V^\sharp$. The cycles
  of $\tau$ on $V^\sharp$ are $\{u\}$, $u\in U^\sharp$, and
  $W^\sharp+u$, $u\in U$. If $C_1$ and $C_2$ are subsets of $V^\sharp$
  such that each $C_i$ lies completely in an $H$--orbit, then we say
  that $C_1$ and $C_2$ are \emph{connected} if they lie in the same
  $H$--orbit. The latter is equivalent to the existence of $h\in H$
  with $C_1\cap C_2^h\ne\emptyset$. Each of the cycles from above lies
  in an $H$--orbit. Consider the graph with vertice set these cycles,
  and let two vertices be connected if and only if the corresponding
  cycles are connected. The aim is to show that this graph is
  connected.

We first show that for each $u\in U$ there is $u\ne u'\in U$ such that
$W^\sharp+u$ and $W^\sharp+u'$ are connected. Suppose that were not
the case. Then, for each $h\in H$,
\[
(W^\sharp+u)^h \subseteq (W^\sharp+u)\cup U^\sharp,
\]
so
\[
W^h \subseteq (W+u-u^h)\cup (U^\sharp-u^h).
\]
First assume that $u\ne0$. Then not each element of $U^\sharp-u^h$ can
be contained in $W^h$, for this would imply the nonsense
$h^{-1}(U^\sharp)\subseteq W^\sharp+u$. Thus we get
\[
\abs{W^h\cap(W+u-u^h)}\ge 2^{m-2}-3.
\]
Let $r$ be the dimension of $W^h\cap(W+u-u^h)$ as an affine space. It
follows that $2^r\ge2^{m-2}-3$, so $r=m-2$ as $m\ge5$. Thus $W^h=W$
for all $h\in H$, contrary to irreducibility of $H$.

Now suppose that $u=0$. Then by the above
\[
(W\setminus h^{-1}(U))^h=(W^\sharp\setminus h^{-1}(U\sharp))^h
\subseteq W^\sharp+u\subset W',
\]
where $W'$ is the $(m-1)$--dimensional space $W\cup(W+u)$. But the
elements in $W\setminus h^{-1}(U))^h$ generate $W$, so $W^h\subset W'$
for all $h$, again contrary to irreducibility of $H$.

Let $u\in U^\sharp$ be such that $W^\sharp$ and $W^\sharp+u$ are
connected. We show that these two cycles must also be connected to
another $W^\sharp+u'$ for $u'\in U^\sharp$ different from $u$. Suppose
that this were not the case. Let $W'$ be the $(m-1)$--dimensional
space $W\cup(W+u)$. Then, similar as above, $(W')^h\subset W'\cup U$,
for all $h\in H$. But $W'\setminus(W'\cap h^{-1}(U))$ generates $W'$,
so $W'$ is $h$--invariant for all $h\in H$, again contrary to
irreducibility of $H$.

From these two steps we see that all the $W^\sharp+u$ for $u\in U$ are
connected. Finally, let $u'\in U^\sharp$. Then also $\{u'\}$ is
connected to some and hence all the $W^\sharp+u$, because $(u')^H$
generates $V$ by irreducibility, so $(u')^H\nsubseteq U^\sharp$.

Thus $H$ is transitive on $V^\sharp$. So we can use the
Hering--Liebeck list \cite[Appendix 1]{Liebeck:affinerank3} of such
groups. Let $L\subseteq\End(V)$ be a maximal field which is normalized
by $H$. ($L$ is unique, see \cite[Lemma 5.2]{Hering:2tran1}.) So
$\abs{L}=2^s$ where $s$ divides $m$, and $H\le\gL_{m/s}(2^s)$. We get
that $\SL_{m/s}(2^s)\le H\le\gL_{m/s}(2^s)$, or $H\le\Spl_m(2)$ with
$m$ even. (The Hering--Liebeck result is more precise, but this rough
version is good enough here.) We claim that $s=1$ or $2$ in the first
case. As $m\ge5$, we have $s\le m<2^{m-2}-1=\ord{\tau}$, hence
$\tau^s$ has exactly $4$ fixed points on $V$ and
$\tau^s\in\GL_{m/s}(2^s)$. Thus $2^s\le4$. If $s=2$, then note that
the determinant of $\tau$ as an $L$--endomorphism of $V$ is a
generator of $L^\star$, hence $\GL_{m/2}(4)\le H$ in this case.

Finally, we need to show that $H\le\Spl_m(2)$ cannot happen. Let
$(\cdot,\cdot)$ be the associated symplectic form on $V$. If $v\in V$
is non--zero, then the stabilizer of $v$ in $\Spl_m(2)$ has two orbits
on $V^\sharp\setminus\{v\}$ -- the orbit of length $2^{m-1}-2$ through
those $v'$ with $(v,v')=0$, and the orbit of length $2^{m-1}$ through
those $v'$ with $(v,v')=1$, see \cite[II.9.15]{Huppert1}. Thus for
$u\in U^\sharp$, either $(u,W^\sharp)=0$ or $(u,W^\sharp)=1$. We aim
to show that the restriction of the symplectic form to $W$ is not
degenerate. This is clear if $(U,W)=0$. So suppose there is $u\in U$
with $(u,W^\sharp)=1$. The orthogonal complement $W^\bot$ intersects
$U$ non-trivially (for if $u_1$ and $u_2$ are different elements in
$U$ with $(u_i,W)=1$, then $(u_1+u_2,W)=0$). So the radical of $W$ has
dimension $\le1$, hence in fact is trivial, because $W$ has even
dimension.

Therefore $W$ is a non-degenerate symplectic space, where $\tau$ acts
irreducibly on. So $\ord{\tau}$ divides $2^{(m-2)/2}+1=2^{\dim
  W/2}+1$, see Lemma \ref{L:Bil:irr}, contrary to
$\ord{\tau}=2^{m-2}-1$.\end{proof}

The proof of the following lemma is straightforward.

\begin{Lemma}\label{L:Aff:kl:Exi}
  Let $m\ge2$, $p$ a prime, and $\FFF_p^m=U\oplus W$ with $U$ and $W$
  invariant under $\tau\in\GL_m(p)$. Assume that $\tau$ acts as a
  Singer cycle on $W$.
\begin{itemize}
\item[(a)] Let $\dim U=1$, and suppose that $\tau$ act as the identity
  on $U$. Choose $u\in U^\sharp$. Then $(\tau,u)\in\AGL_m(p)$ acts as
  an element with cycle lengths $p$ and $p^m-p$ on $\FFF_p^m$.
\item[(b)] Let $\dim U=2$, $p=2$, and suppose that $\tau$ act as an
  involution on $U$. Choose $u\in U$ with $u\ne u^\tau$. Then
  $(\tau,u)\in\AGL_m(2)$ acts as an element with cycle lengths $4$ and
  $2^m-4$ on $\FFF_2^m$.
\end{itemize}
\end{Lemma}

\subsubsection{Proof of Part \protect\ref{T:kl:I} of Theorem
  \protect\ref{T:kl}}

We may assume that $A\le\AGL_m(p)$, acting on $N=\FFF_p^m$. Let
$\sigma$ be the element with the cycle lengths $k$ and $l$. First note
that $k$ divides $l$, for otherwise $\sigma^l$ would fix $l>p^m/2$
points, which of course is nonsense. So $k=p^r$, $l=p^r(p^{m-r}-1)$
for some $r\in\NNN_0$.

First suppose $k<l$. Then $\sigma^k$ fixes exactly $k=p^r$ points on
$N$. Without loss $\sigma^k\in\GL_m(p)$, so the fixed point set of
$\sigma^k$ is a subspace $N_1$ of $N$. But the elements of $N_1$
constitute the $k$--cycle of $\sigma$, so $\sigma$ acts as an affine
map of order $p^r$ on the $r$--dimensional space $N_1$. Apply Lemma
\ref{L:Aff:p} to see that $r\in\{0,1,2\}$, and $r\le1$ if $p>2$.

If $k=l$, then of course $p=2$ and $k=l=\ord{\sigma}=2^{m-1}$. Lemma
\ref{L:Aff:p} gives $2^{m-2}\le m$, hence $r=m-1\le3$.

We need to determine the possible groups $A$. If $k=1$ we use a result
of Kantor \cite{Kantor:Singer} which
classifies linear groups over a finite field containing an element
which cyclically permutes the non-zero elements. Note that $\sigma$ is
just such an element.

Now suppose $k=p$. The element $\tau:=\sigma^p\in\GL_m(p)$ fixes a
line $U\cong\FFF_p$ pointwise. As $\gcd(\ord{\tau},p)=1$, Maschke's
Theorem gives a complement $W$ of $U$ which is $\tau$--invariant. As
$\tau$ has cycles of length $\ord{\tau}=p^{m-1}-1=\abs{W^\sharp}$ on
$W^\sharp$, we see that $\tau$ permutes the elements of $W^\sharp$
cyclically. Look at the action which $A_1$ and $\tau$ induce on the
projective space $P(\FFF_p^m)$. The element $\tau$ has a fixed point
(corresponding to $U$), a cycle of length $(p^{m-1}-1)/(p-1)$
(corresponding to $W$) and a cycle of length $p^{m-1}-1$
(corresponding to $u+W$ for $0\ne u\in U$).

Lemma \ref{L:PGL2:irr} handles the case $m=2$. So for the rest of this
argument we assume $m\ge3$. We contend that $A_1$ is transitive on
$P(\FFF_p^m)$. By primitivity, $A_1$ is irreducible on $\FFF_p^m$, so
it moves the fixed point of $\tau$ as well as the cycle of length
$(p^{m-1}-1)/(p-1)$. So if $A_1$ were not transitive, then $A_1$ would
leave $U\cup W$ invariant. Let $a\in A_1$ with $W^a\ne W$, and choose
$w\in W$ with $w^a\in U$. Then $W\setminus\FFF_pw$ is invariant under
$a$. But $W\setminus\FFF_pw$ generates $W$ because of $m\ge3$, a
contradiction.

Looking at $\tau$ we see that $A_1$ is a transitive group on
$P(\FFF_p^m)$ of rank at most $3$. If $A_1$ is even doubly transitive,
then by Proposition \ref{P:PGL:2tr} either $p=2$, $m=4$, and
$A_1=\alt_7$, or $\SL_m(p)\le A_1$. It is easy to see that the
determinant of $\tau$ is a generator of the multiplicative group of
$\FFF_p$. Thus $A_1=\GL_m(p)$ in the latter case.

Next assume that $A_1$ has rank $3$ on $P(\FFF_p^m)$. We first use
Lemma \ref{L:Impr:rk3} to see that $A_1$ is also primitive. For if
not, then $1+\alpha$ divides $(p-1)\alpha$ with
$\alpha=(p^{m-1}-1)/(p-1)$, so $1+\alpha$ divides $p-1$, which of
course is nonsense.

So we can apply Proposition \ref{P:PGL:rk3}. Suppose that
$\Spl_m(p)\trianglelefteq A_1$. However, $\Spl_m(p)$ has rank $3$ on
$P(\FFF_p^m)$ and subdegrees $1$, $p^{m-1}$ and $p(p^{m-2}-1)/(p-1)$
(see \cite[II.9.15]{Huppert1}), which is not compatible with the cycle
lengths of $\tau$ we determined above. Thus the other possibility of
the proposition holds, that is $p=2$, $m=3$, $A_1=\gL_1(8)$, which
indeed gives case \ref{T:kl:I:d:ii} in the theorem.

To see that the groups listed in \ref{T:kl:I:b} and \ref{T:kl:I:c}
indeed have an element with cycle lengths $p$ and $p^m-p$ use Lemma
\ref{L:Aff:kl:Exi}, and similarly for the cycle lengths $4$ and
$2^m-4$ in \ref{T:kl:I:d}.

Next we look at the case $p=2$ and $k=4$. The case $m\le4$ is done by
inspection, so assume $m\ge5$. Set $\tau:=\sigma^4$. As above we see
that $\FFF_2^m=U\oplus W$ with $\dim U=2$, $\tau$ is trivial on $U$,
and acts as a Singer cycle on $W$. In view of Proposition
\ref{P:Hering} we need to show that $GL_{m/2}(4)\le
A_1\le\gL_{m/2}(4)$ is not possible. Suppose that were the case. Write
$\sigma=(\beta,v)$ with $\beta\in\gL_{m/2}(4)$, and
$v\in\FFF_4^{m/2}$. First note that $U^\beta+v=U$, so $v\in U$ and
$\beta$ leaves $U$ invariant. This easily implies $\beta^4=\tau$. The
intersection $W\cap W^{\beta^2}$ is not trivial (by dimension reasons)
and invariant under $\beta^2$, so in particular invariant under
$\tau$. Hence $W=W^{\beta^2}$ by irreducibility of $\tau$ on $W$. But
$W\cap W^\beta$ is also $\tau$--invariant, so $\beta$ leaves $W$
invariant by the same argument. As $\tau=\beta^4$ permutes the
elements of $W^\sharp$ cyclically, the same holds true for $\beta$.
Hence $\beta$ has odd order $2^{m-2}-1$ when restricted to $W$. In
particular, $\beta\in\GL_{m/2}(4)$. On the other hand, as $(\beta,v)$
has order $4$ on $U$, the order of $\beta$ on $U$ must be a divisor of
$4$. However, $\abs{\GL_1(4)}=3$, hence $\beta$ is trivial on $U$. But
then $(\beta,v)$ has order $2$ on $U$, a contradiction.

The case $k=l=8$, $p=2$, is most conveniently done by inspection using
GAP \cite{GAP}.

\subsection{Product Action}

Set $\Delta=\{1,2,\dots,r\}$ for $r\ge2$, and let $m\ge2$ be an
integer. Then the wreath product $\sym_r\wr\sym_m=
(\sym_r\times\sym_r\times\cdots\times\sym_r)\rtimes\sym_m$ acts in a
natural way on $\Omega:=\Delta\times\Delta\times\cdots\times\Delta$.
We say that a permutation group $A$ acts via the \emph{product
  action}, if it is
permutation equivalent to a transitive subgroup of $\sym_r\wr\sym_m$
in this action.

In order to avoid an overlap with the affine permutation groups, we
quickly note the easy

\begin{Lemma}\label{L:ProdAff}
  Let $A$ be a primitive subgroup of $\sym_r\wr\sym_m$ where $r\le4$.
  Then $A$ is affine.
\end{Lemma}

\begin{proof} Let $N$ be the minimal normal
  subgroup of $\sym_r\wr\sym_m$. Then $N$ is elementary abelian of
  order $r^m$. If $A$ intersects $N$ non-trivially, then $N\cap A$ is
  a minimal normal subgroup of $A$, and the claim follows. So suppose
  that $\abs{A\cap N}=1$. Then $A$ embeds into $(\sym_r\wr\sym_m)/N$.
  But $r^m$ divides $\abs{A}$ by transitivity, so $r^m$ divides
  $(r!)^mm!/r^m$. We get that $2^m$ divides $m!$ if $r=2$ or $4$, and
  $3^m$ divides $m!$ if $r=3$. But if $p$ is a prime, then the
  exponent of $p$ in $m!$ is
  $\sum_{\nu\ge0}\genfrac{[}{]}{}{}{m}{p^\nu} <
  \sum_{\nu\ge0}\frac{m}{p^\nu} = \frac{m}{p-1}\le m$, a
  contradiction.\end{proof}

\begin{Remark*} One might expect that any primitive subgroup of an
  affine group is affine. However, that is \emph{not} the case. There
  seem to be very few counter-examples. The smallest is as follows:
  Set $A=\AGL_3(2)=C_2^3\rtimes A_1$. Then it is known (see
  e.g.~\cite[page 161]{Huppert1}) that $H^1(\GL_3(2),C_2^3)=C_2$. So
  there is a complement $U$ of $C_2^3$ in $A$ which is not conjugate
  to $A_1$. One checks that $U$ acts primitively on the $8$ points via
  $U\cong\GL_3(2)\cong\PSL_2(7)$.\end{Remark*}

The following two lemmas are trivial but useful.

\begin{Lemma}\label{L:Prod:CyLe}
  Let $\Delta_1$, $\Delta_2$, \dots, $\Delta_m$ be finite sets, and
  $g_i$ be in the symmetric group of $\Delta_i$. Let $o_i$ be the
  cycle length of $g_i$ through $\delta_i\in\Delta_i$. Then the cycle
  length of $(g_1,g_2,\dots,g_m)$ through
  $(\delta_1,\delta_2,\dots,\delta_m)\in
  \Delta_1\times\Delta_2\times\dots\times\Delta_m$ is
  $\lcm(o_1,o_2,\dots,o_m)$.  In particular,
  $\delta_1^{\gen{g_1}}\times\delta_2^{\gen{g_2}}\times\cdots
  \times\delta_m^{\gen{g_m}}$ is the orbit of
  $\gen{(g_1,g_2,\dots,g_m)}$ through
  $(\delta_1,\delta_2,\dots,\delta_m)$ if and only if the $o_i$ are
  relatively prime.
\end{Lemma}

\begin{Lemma}\label{L:Prod:CyNu}
  Let $\Delta_1$, $\Delta_2$, \dots, $\Delta_m$ be finite sets, and
  $g_i$ be in the symmetric group of $\Delta_i$. Let $c_i$ be the
  number of cycles of $g_i$ on $\Delta_i$. Then $(g_1,g_2,\dots,g_m)$
  has at least $c_1c_2\cdots c_m$ cycles on
  $\Delta_1\times\Delta_2\times\cdots\times\Delta_m$.
\end{Lemma}

\begin{Lemma}\label{L:Prod:Cyc}
  Suppose $r\ge5$ and $m\ge2$. Let $g=(\s1,\s2,\dots,\s m)\tau$ be an
  element of $\sym_r\wr\sym_m$, with $\s i\in\sym_r$ and
  $\tau\in\sym_m$ an $m$--cycle. Then $g$ has at least $3$ cycles in
  the product action.
\end{Lemma}

\begin{proof}
  Set
  $\overline{g}:=g^m=(\overline{\s1},\overline{\s2},\ldots,\overline{\s
    m})\in\sym_r^m$.  Then
\[
\overline{\s i}=\s i\s{{i+1}}\cdots\s m\s1\cdots\s{{i-1}},
\]
so in particular the $\overline{\s i}$ are pairwise conjugate in
$\sym_r$.  Suppose that $g$ has at most $2$ cycles. Then
$\overline{g}$ has at most $2m$ cycles.
  
Let $\lambda$ be the number of cycles of $\overline{\s1}$. Then
$\overline{g}$ has at least $\lambda^m$ cycles by Lemma
\ref{L:Prod:CyNu}, hence $\lambda^m\le 2m$. This gives $\lambda=1$
unless $m=2$ and $\lambda=2$. If $\lambda=1$, then $\overline{g}$ has
$r^{m-1}$ cycles by Lemma \ref{L:Prod:CyLe}, so $r^{m-1}\le 2m$, hence
$r\le4$, a contradiction. So suppose that $m=2$ and $\overline{\s1}$
has two cycles.  Then $\overline{g}$ has obviously at least $6>2m$
cycles, a contradiction.
\end{proof}

\subsubsection{Proof of Part \ref{T:kl:II} of Theorem \ref{T:kl}}

We assume that $A\le\sym_r\wr\sym_m$ with $r\ge5$ (by Lemma
\ref{L:ProdAff}) and $m\ge2$. Let $g=(\s1,\s2,\ldots,\s m)\tau$ with
$\s i\in\sym_r$, $\tau\in\sym_m$.
  
Assume that $g$ has exactly $2$ cycles. By the previous lemmas, we get
that $m=2$ and $\tau=1$, one of the $\s i$ must be an $r$--cycle, and
the other $\s i$ has two cycles, with lengths relatively prime to $r$.
  
We need to determine the groups which arise this way. The description
of the product action as in \cite{LPS} shows that there is a primitive
group $U$ with socle $S$ acting on $\Delta=\{1,2,\dots,r\}$, such that
$S\times S\trianglelefteq A\le (U\times U)\rtimes C_2$. Let
$g=(\s1,\s2)$ be the element with the two cycles from above. Then
$(\s2,\s1)\in (U\times U)\rtimes C_2$.  Thus $U$ contains an
$r$--cycle, and an element with two cycles of coprime lengths. In
particular, $U$ is not contained in the alternating group $\alt_r$,
and so is not simple. Furthermore, $U$ is not affine. Taking
Theorems \ref{T:Feit:nZykel} and \ref{T:PM:(n-1)} together gives
that either $U=\PGL_2(p)$ for a prime $p\ge5$, or $U=\sym_r$ for
$r\ge5$. The element $g$ shows that $U\times U\le A$, but $U\times U$
is not primitive, so $A=(U\times U)\rtimes C_2$, and the claim
follows.

\begin{Remark} Case \ref{T:kl:I:c} of Theorem \ref{T:kl}, that
  is $A<\GL_2(p)$ for a prime $p>2$ and $A_1$ the group of monomial
  matrices, can also be seen as a product action, namely as
  $A=(\AGL_1(p)\times\AGL_1(p))\rtimes C_2$ on $p^2$ points.
\end{Remark}

\subsection{Regular Action} As an immediate consequence of the previous
section we obtain

\begin{Theorem}\label{T:Reg:kl}
  Let $A$ be a primitive non--affine permutation group with a regular
  normal subgroup. Then $A$ does not contain an element with two
  cycles.
\end{Theorem}

\begin{proof} Let $N$ be a regular normal subgroup of $A$. Then, by
  regularity, $N$ is a minimal normal subgroup of $A$, so $N\cong L^m$
  for some simple non--abelian group $L$ and $m\ge2$. Identify $N$
  with the set of points $A$ is acting on, and let $C$ be the
  centralizer of $N$ in the symmetric group $\sym(N)$ of $N$. If $N$
  acts from the right on $N$, then $C\cong N$ acts from the left on
  $N$. Set $H=L\times L$, and let the first and second component act
  from the left and from the right, respectively. Then $A$ is
  contained in the wreath product $H\wr\sym_m$ in product action, see
  \cite[page 392]{LPS}. Now apply Theorem \ref{T:kl} to see that
  this cannot occur, a distinguishing property of $H$ being that it is
  not doubly transitive (in contrast to $\PGL_2(p)$).\end{proof}

\subsection{Diagonal Action} Let $\Si$ be a non--abelian simple group, and
$m\ge2$ an integer. Set $N:=\Si^m$. Let $N$ act on itself by
multiplication from the right. Furthermore, let the symmetric group
$\sym_m$ act on $N$ by permuting the components, and $\Aut(\Si)$ act
on $N$ componentwise. Define an equivalence relation $\sim$ on $N$ by
$(l_1,l_2,\dots,l_m)\sim(cl_1,cl_2,\dots,cl_m)$ for $c\in\Si$. The
above actions respect the equivalence classes, so we get a permutation
group $D$ acting on the set $N/\negthickspace\sim$ of size
$\abs{\Si}^{m-1}$. Note that the diagonal elements of $N$ in right
multiplication induce inner automorphisms of $\Si$ on
$N/\negthickspace\sim$, for
$(i^{-1}l_1i,i^{-1}l_2i,\dots,i^{-1}l_mi)\sim
(l_1,l_2,\dots,l_m)(i,i,\dots,i)$.

We say that a permutation group $A$ acts in \emph{diagonal
action}, if it embeds as
a transitive group of $D$ with $N\le A$.

We begin with a technical

\begin{Proposition}\label{P:Diag:Orbs}
  Let $\Si$ be a non--abelian simple group, $m\ge2$ be an integer, and
  $D$ be the group in diagonal action as above. Let $o(\Out(\Si))$ and
  $o(\Si)$ be the largest order of an element in $\Out(\Si)$ and
  $\Si$, respectively. Then each element of $D$ has at least
  $\frac{1}{o(\Out(\Si))\abs{\Si}}(\abs{\Si}/o(\Si))^m$ cycles.
\end{Proposition}

\begin{proof} Choose an element in $D$. Raise it to the smallest power
  such that the contribution from $\Out(\Si)$ disappears. Let $\sigma\in
  N\rtimes\sym_m$ be this element. Set $o=o(\Si)$. We are done once we
  know that $\sigma$ has at least $\frac{1}{\abs{\Si}}(\abs{\Si}/o)^m$
  cycles. Write $\sigma=(\sigma_1,\sigma_2,\dots,\sigma_m)\tau$ with
  $\tau\in\sym_m$ and $\sigma_i\in\Si$. Let $\tau$ have $u$ cycles of
  lengths $\rho_1,\rho_2,\dots,\rho_u$.

Without loss assume that the first $\rho_1$ coordinates of $N=\Si^m$ are
permuted in an $\rho_1$--cycle $(1\;2\;\cdots\;\rho_1)$. Write $\rho$
for $\rho_1$. Then $\sigma^{\rho}$ acts by right multiplication with
\[
(\overline{\sigma_1},\overline{\sigma_2},\ldots,\overline{\sigma_\rho})=
(\sigma_1\sigma_2\cdots \sigma_\rho,\sigma_2\sigma_3\cdots
\sigma_\rho\sigma_1,\dots,\sigma_\rho\sigma_1\cdots
\sigma_{\rho-1})\in\Si^{\rho}
\]
on these first $\rho$ coordinates. Note that all the elements
$\overline{\sigma_i}$ have the same order $o'$ because they are
conjugate in $\Si$. So, by Lemma \ref{L:Prod:CyLe}, $\sigma^\rho$
induces $\abs{\Si}^\rho/o'\ge\abs{\Si}^\rho/o$ cycles on $\Si^\rho$, thus
$\sigma$ induces at least $\abs{\Si}^\rho/(\rho o)$ cycles on $\Si^\rho$.
Apply this consideration to the other $\tau$--cycles and use Lemma
\ref{L:Prod:CyNu} to see that the number of cycles of $\sigma$ on $N$
is at least
\begin{align*}
\prod_{i=1}^u\frac{\abs{\Si}^{\rho_i}}{\rho_io} &=
\frac{\abs{\Si}^m}{o^u}\prod_{i=1}^u\frac{1}{\rho_i}\\
&\ge \abs{\Si}^m\left(\frac{u}{mo}\right)^u,
\end{align*}
where we used the inequality between the arithmetic and geometric mean
in the last step. But the function $(x/(mo))^x$ is monotonously
decreasing for $0\le x\le mo/e$. Note that $o\ge5$ (because a group
with element orders $\le4$ is solvable), so $mo/e>m$, but $u\le m$. So
the above expression is $\ge(\abs{\Si}/o)^m$. Furthermore, the number of
cycles of $\sigma$ on $N$ is at most $\abs{\Si}$ times the number of
cycles on $N/\negthickspace\sim$. From that we get the
assertion.\end{proof}

\begin{Theorem}\label{T:Diag:kl}
  Let $A$ be a primitive permutation group in diagonal action. Then
  $A$ does not contain an element with at most two cycles.
\end{Theorem}

\begin{proof} Suppose there is a counterexample $A$, with associated
  simple group $\Si$. Proposition \ref{P:Diag:Orbs} gives, as $m\ge2$,
\[
\abs{\Si}\le2o(\Si)^2o(\Out(\Si)).
\]
If $\Si$ is sporadic, then use list \ref{Tb:Sp} on page \pageref{Tb:Sp}
along with the group orders given in the atlas \cite{ATLAS} to see
that this inequality has no solution. Next suppose that $\Si=\alt_n$ is
alternating. Then $\Out(\Si)=C_2$ if $n\ne6$, and
$\Out(\alt_6)=C_2\times C_2$, so $o(\Out(\Si))=2$ in any case (see
e.g.~\cite[II.5.5]{Huppert1}). Use the bound $o(\Si)\le e^{n/e}$ from
Proposition \ref{P:Sn:ord} to see that only $n=5$ is possible with
$m=2$. But it is easy to take into account the possible outer
automorphism and show along the lines of the previous proposition that
the minimal number of cycles of an element in $A$ is $4$ (all of
length $15$), or one checks that with a GAP computation.

So we are left with the case that $\Si$ is simple of Lie type. Using the
information about $\Out(\Si)$ and $o(\Si)$ in the Tables \ref{Tb:ExGr}
(page \pageref{Tb:ExGr}) and \ref{Tb:ClGr} (page \pageref{Tb:ClGr})
and in Section \ref{SS:ClGr} together with the order of $\Si$ given for
instance in the atlas \cite{ATLAS}, one sees that the only group which
does fulfill the above inequality is $\Si=\PSL_2(7)$. (One also has to
use the atlas \cite{ATLAS} in some small cases where the given bounds
for $o(\Si)$ are too coarse in order to exclude $\Si$.)

However, the proof of the proposition above shows that we have $m=2$,
$u=2$, and $\ord{\overline{\s1}}\ord{\overline{\s2}}\ge168/4=42$,
hence $\ord{\overline{\s1}}=\ord{\overline{\s2}}=7$, so $\sigma$ has
at least $168^2/(7\cdot168)=24$ cycles on $\Si^2/\negthickspace\sim$, a
clear contradiction.\end{proof}

\subsection{Almost Simple Action}

By what we have seen so far, the only remaining case is the almost
simple action. The aim of the following sections is to show that only
the cases listed in part \ref{T:kl:III} of Theorem \ref{T:kl} appear.
See Section \ref{SS:AS:Proof} where all the results achieved in the
following sections are bundled to give a proof of this assertion.

Many cases of almost simple permutation groups can be ruled out by
comparing element orders with indices of (maximal) subgroups of almost
simple groups, though some other require finer arguments.

We make the following

\begin{Definition} For a finite group $X$ let $\mu(X)$ be the
  smallest degree of a faithful, transitive permutation representation
  of $X$, and $o(X)$ the largest order of an element in $X$.
\end{Definition}

We use the trivial

\begin{Lemma}\label{L:kl:Deg-Ord}
  Let $A$ be a transitive permutation group of degree $n$, and let
  $\sigma\in A$ have two cycles in this action. Then
\begin{align}
  n &\le        2\ord{\sigma}\label{Eq:muo}\\
  n &\le 3\ord{\sigma}/2,\text{ if $n$ is odd.}\label{Eq:muo:odd}
\end{align}
\end{Lemma}

\subsection{Alternating Groups}\label{SS:Alt} Using methods and
results from analytic number theory, one can show that the logarithm
of the maximal order of an element in $\sym_n$ is asymptotically
$\sqrt{n\log n}$, see \cite[\S61]{Landau:Prim}. Here, the following
elementary but weaker result is good enough for us -- besides, we need
an exact bound rather than an asymptotic bound anyway.

\begin{Proposition}\label{P:Sn:ord}
  The order of an element in $\sym_n$ is at most $e^{n/e}$ for all
  $n\in\NNN$, and at most $(n/2)^{\sqrt{n/2}}$ for $n\ge6$. (Here
  $e=2.718\ldots$ denotes the Euler constant.)
\end{Proposition}

\begin{proof}
  Let $\nu_1,\nu_2,\dots,\nu_r$ be the different cycle lengths $>1$ of
  an element $g\in\sym_n$. Then
\begin{equation*}\label{Snord1}
\ord{g}=\lcm(\nu_1,\nu_2,\dots,\nu_r)\le \nu_1\nu_2\cdots\nu_r,
\end{equation*}
and
\begin{equation}\label{Snord2}
\nu_1+\nu_2+\dots+\nu_r\le n.
\end{equation}
The inequality between the arithmetic and geometric mean yields
\begin{align*}
\ord{g} &=   \lcm(\nu_1,\nu_2,\dots,\nu_r)\\
        &\le \nu_1\nu_2\cdots\nu_r\\
        &\le \left(\frac{\sum\nu_i}{r}\right)^r\\
        &\le \left(\frac{n}{r}\right)^r.
\end{align*}
The function $x\mapsto(n/x)^x$ is increasing for $0<x\le n/e$, and
decreasing for $x>n/e$. From that we obtain the first inequality.

Suppose that $\nu_1<\nu_2<\dots<\nu_r$. Then $\nu_i\ge i+1$, and we
obtain
\[
n\ge\sum\nu_i\ge 2+3+\dots+r+(r+1)=\frac{r^2+3r}{2}>\frac{r^2}{2}.
\]
If $n>2e^2=14.7\ldots$, then $r<\sqrt{2n}<n/e$, and the claim follows
from the monotonicity consideration above. Check the cases $6\le
n\le14$ directly.
\end{proof}

Now suppose that $\alt_n\le A\le\Aut(\alt_n)$ for $n\ge5$. Note that
except for $n=6$, $\Aut(\alt_n)=\sym_n$ by \cite[II.5.5]{Huppert1}. We
exclude $n=6$ in this section, and treat this case in Section
\ref{SS:ClGr} about classical groups, because $\alt_6\cong\PSL_2(9)$.

So $A_1$ is a maximal subgroup of $A$ not containing $\alt_n$. Let
$\sigma\in A$ have at most two cycles on $A/A_1$. We regard $A_1$ as a
subgroup of $\sym_n\ge A$ in the natural action on $\{1,2,\dots,n\}$
points. There are three possibilities for $A_1$ with respect to this
embedding: $A_1$ is intransitive, or transitive but imprimitive, or
primitive. We treat these three possibilities separately.

\paragraph{$A_1$ intransitive.} $A_1$ leaves a set of size $m$
invariant, with $1\le m<n$. Denote by $M_m$ the subsets of size $m$ of
$\{1,2,\dots,n\}$. By maximality of $A_1$ in $A$ and transitivity of
$A$ on $M_m$ we see that $A_1$ is the full stabilizer in $A$ of a set
of $m$ elements, thus the action of $A$ is given by the action on
$M_m$. If $m=1$, then we have the natural action of $A$, leading to
case \ref{T:kl:III:An} in Theorem \ref{T:kl:III}. So for the remainder
assume $m\ge2$.

First consider the case that $\sigma$ is an $n$--cycle in the natural
action. One of the two cycles of $\sigma$ has length at least
$\binom{n}{m}/2$, so $n\ge\binom{n}{m}/2\ge n(n-1)/4$, thus $n=5$.
This case really occurs, and gives case \ref{T:kl:III:A5(10)} in
Theorem \ref{T:kl:III}.

Next suppose that $\sigma$ is not an $n$--cycle. Then $\sigma$ leaves
(on $\{1,2,\dots,n\}$) a set $S$ of size $1\le\abs{S}\le n/2$
invariant. Without loss $m\le n/2$ (as the action on the $m$--sets is
the same as the action on the $(n-m)$--sets). Note that $\sigma$
cannot be an $(n-1)$--cycle by an order argument as above. So we can
assume $\abs{S}\ge2$. For $i=0,1,2$ choose sets $S_i$ of size $m$,
such that $i$ points of $S_i$ are in $S$, and the remaining $m-i$
points are in the complement of $S$. Then these three sets of course
are not conjugate under $\gen{\sigma}$.

\paragraph{$A_1$ transitive but imprimitive.} Let $1<u<n$ be the
size of the blocks of a non--trivial system of imprimitivity. Then
$v:=n/u$ is the number of blocks, and $A_1=(\sym_u\wr\sym_v)\cap
A=((\sym_u)^v\rtimes\sym_v)\cap A$ in the natural action (not to
mistake with the product action).

The index of $A_1$ in $A$ thus is $n!/((u!)^vv!)$. We will use the
bounds in Lemma \ref{L:kl:Deg-Ord} and Proposition \ref{P:Sn:ord} to
see that this case does not occur. The proof is based on the following

\begin{Lemma}\label{L:u!^vv!}
  Let $u,v\ge2$ be integers, then
\begin{equation}\label{Eq:u!^vv!}
u!^vv!<\frac{1}{2}\frac{(uv)!}{e^{uv/e}},
\end{equation}
except for $(u,v)=(2,2)$, $(3,2)$, $(4,2)$, and $(2,3)$.
\end{Lemma}

\begin{proof} We contend that if the inequality \eqref{Eq:u!^vv!} holds
for $(u,v)$, then it holds also for $(u,v+1)$. First
\[
3<4.31\ldots=\left(\frac{3}{e^{1/e}}\right)^2\le
\left(\frac{3}{e^{1/e}}\right)^u,
\]
hence
\[
e^{u/e}<3^{u-1}\le(v+1)^{u-1}.
\]
This implies
\begin{equation}\label{uv2}
(v+1)e^{u/e}<(v+1)^u.
\end{equation}
But
\[
v+1\le\frac{uv+i}{i}
\]
for $i=1,2,\dots,u$, so taking the product over these $i$ yields
\[
(v+1)^u\le\binom{uv+u}{u},
\]
so
\[
(v+1)e^{u/e}\le\binom{uv+u}{u}
\]
by \eqref{uv2}. Multiply the resulting inequality
\[
u!(v+1)<\frac{(uv+u)!}{(uv)!e^{u/e}}
\]
with \eqref{Eq:u!^vv!} to obtain the induction step for $v$.

Next we show that \eqref{Eq:u!^vv!} holds for $v=2$ and $u\ge7$. As
$\binom{2u}{u}$ appears as the biggest binomial coefficient in the
expansion of $(1+1)^{2u}$, we obtain
$\binom{2u}{u}\ge\frac{1}{2u+1}2^{2u}$. Inequality \eqref{Eq:u!^vv!} for
$v=2$ reduces to
\[
\binom{2u}{u}>4e^{2u/e}.
\]
So we are done once we know that
\[
\frac{1}{2u+1}2^{2u}>4e^{2u/e},
\]
which is equivalent to
\[
\left(\frac{2}{e^{1/e}}\right)^{2u}>4(2u+1).
\]
But it is routine to verify this for $u\ge7$.

In order to finish the argument, one verifies \eqref{Eq:u!^vv!}
directly for $u<7$ and the least value of $v$ where the inequality is
supposed to hold.\end{proof}

As $uv\ge5$ and $uv\ne6$ by our assumption, we have the only case
$u=4$, $v=2$. But $8!/(4!^22!)=35$, and the maximal order of an
element in $\sym_8$ is $15$, contrary to Lemma \ref{L:kl:Deg-Ord}.

\paragraph{$A_1$ primitive.} Now suppose that $A_1$ is primitive on
$\{1,2,\dots,n\}$, hence $\left[\frac{n+1}{2}\right]!\le [\sym_n:A_1]$
by a result of Bochert, see \cite{Bochert:prim}
or \cite[14.2]{Wielandt}. Here $[x]$ denotes the biggest integer less
than or equal $x$. As $A$ has index at most $2$ in $\sym_n$, we obtain
from Lemma \ref{L:kl:Deg-Ord} and Proposition \ref{P:Sn:ord}
\[
\left[\frac{n+1}{2}\right]!\le 2[A:A_1]\le 4e^{n/e}.
\]

However, one verifies that for $n=9$ and $12$ the following holds
\begin{equation}\label{boch}
\left[\frac{n+1}{2}\right]!> 4e^{n/e}.
\end{equation}
But if \eqref{boch} holds for some $n\ge9$, then it holds for $n+2$ as
well, as the left side grows by the factor $[(n+3)/2]$, whereas the
right side grows by the factor $e^{2/e}<[(n+3)/2]$.

So we are left to look at the cases $n\in\{5,7,8,10\}$.

Suppose $n=5$. The only maximal transitive subgroup of $\sym_5$ not
containing $\alt_5$ is $A:=C_5\rtimes C_4$, and the only maximal
transitive subgroup of $\alt_5$ is $A\cap\alt_5=C_5\rtimes C_2$. So
the index is $6$, and these cases indeed occur and give
\ref{T:kl:III:PSL2(p)} in Theorem \ref{T:kl:III} for $p=5$.

Now assume $n=7$. The only transitive subgroups of $\sym_7$ which are
maximal subject to not containing $\alt_7$ are $\AGL_1(7)$ and
$\PSL_3(2)$. Of course, the index of $\AGL_1(7)$ in $\sym_7$ is much
too big. The group $\PSL_3(2)$ is contained in $\alt_7$, and has index
$15$. But the maximal order of an element in $\alt_7$ is $7<15/2$, so
this case does not occur by Lemma \ref{L:kl:Deg-Ord}.

Now assume $n=8$. Similarly as above, we see that the only case which
does not directly contradict Lemma \ref{L:kl:Deg-Ord} is
$A_1=\AGL_3(2)$ inside $\PSL_4(2)\cong\alt_8$. But then $A=\PSL_4(2)$
in the natural degree $15$ action on the projective space. Lemma
\ref{L:PGL:kl} shows that this case actually does not occur.

Finally, if $n=10$, then we keep Bochert's bound, but use Proposition
\ref{P:Sn:ord} to see that the order of an element in $\sym_{10}$ is
at most $5^{\sqrt{5}}=36.55\dots$, hence at most $36$. (The exact
bound is $30$.) So $5!\le2\cdot36$ by Lemma \ref{L:kl:Deg-Ord}, a
contradiction.

\subsection{Sporadic Groups}\label{SS:Sp} Let $\Si$ be one of the
$26$ sporadic groups. Table \ref{Tb:Sp} on page \pageref{Tb:Sp}
contains information about small permutation degrees, big element
orders, and the outer automorphism group. The atlas \cite{ATLAS}
contains all this information except for the maximal subgroups of the
Janko group $J_4$, the Fischer groups $Fi_{22}$, $Fi_{23}$, and
$Fi_{24}'$, the Thompson group $Th$, the baby monster $B$, and the
monster group $M$. For the groups $J_4$, $Fi_{22}$, $Fi_{23}$, and
$Th$ we find the necessary information in \cite{KW:J4},
\cite{KW:Fi22}, \cite{KPW:Fi23}, and \cite{Linton:Th}, respectively.
The bounds for the groups $Fi_{24}'$, $B$, and $M$ are not sharp, and
have been obtained as follows from the character tables in
\cite{ATLAS}: If $M$ is a proper subgroup of $\Si$ with index $n$,
then the permutation character for the action of $\Si$ on $\Si/M$ is
the sum of the trivial character and a character of degree $n-1$ which
does not contain the trivial character. Thus $n-1$ is at least the
degree of the smallest non--trivial character of $\Si$. (In view of
the applications we have in mind we could have used this argument in
most other cases as well.)

Now $\Si\le A\le\Aut(\Si)$ for a sporadic group $\Si$. Let $\sigma\in
A$ be an element with only two cycles in the given permutation action.
By Lemma \ref{L:kl:Deg-Ord} we get
$\mu(\Si)\le2\abs{\Out(\Si)}o(\Si)$. We see that the only possible
candidates for $\Si$ are the five Mathieu groups.

The atlas \cite{ATLAS} provides the permutation characters of the
simple groups of not too big order on maximal subgroups of low index.
In the case of the Mathieu groups in the representations which are
possible, we thus can immediately read off the cycle lengths of an
element. Namely the atlas also tells in which conjugacy class a power
of an element lies, so we can compute the fixed point numbers of all
powers of a fixed element.

$\mathbf{\Si=\M11}$. Then $A=\M11$ either in the natural action of
degree $11$, or in the action of degree $12$. The degree $11$ case
cannot occur for the following reason. By Lemma \ref{L:kl:Deg-Ord}
$\ord{\sigma}\ge(2/3)11$, so $\ord{\sigma}=8$ or $11$. An element of
order $11$ is an $11$--cycle. An element of order $8$ has a fixed
point, so if it would have two cycles, the other cycle length had to
be $10$, which is nonsense. Now look at the degree $12$ action. Then
of course an element of order $11$ has cycle lengths $1$ and $11$, and
one readily checks that an element of order $8$ has cycle lengths $4$
and $8$, whereas an element of order $6$ has a fixed point, hence must
have more than $2$ cycles.

$\mathbf{\Si=\M12}$. The smallest degree of a faithful primitive
representation of $\Aut(\M12)$ is $144$ (see \cite{ATLAS}), which is
considerably too big. So we have $A=\M12$ in its natural action. As
$\M11<\M12$, the elements of order $11$ and $8$ in $\M11$ with only
two cycles appear also in $\M12$. Besides them, an element of order
$10$ has cycle length $2$ and $10$, and an element in one of the two
conjugacy classes of elements of order $6$ has cycle lengths $6$.

$\mathbf{\Si=\M22}$. We have the natural action of $\Si$ of degree $22$,
and $A\le\M22\rtimes C_2$. An element of order $11$ has two cycles of
length $11$. An element in $\Si$ of order $8$ has cycle lengths $2$,
$4$, $8$, $8$, so this element cannot be the square of an element with
only $2$ cycles. An element of order $7$ has one fixed point, so it
cannot arise either. And an element in $\Si$ of order $6$ has $6$
cycles, so is out too.

$\mathbf{\Si=\M23}$. Here $A=\M23$ in the natural action of degree $23$.
An element of order $23$ is a $23$--cycle. Looking at the fixed points
of elements of order $3$ and $5$ we see that an element of order $15$
has cycle lengths $3$, $5$, and $15$. Similarly, an element of order
$14$ has cycle lengths $2$, $7$, and $14$. So this group does not
occur at all.

$\mathbf{\Si=\M24}$. Here $A=\M24$ in the action on $24$ points. One
quickly checks that the elements of order $14$ and $15$ have a fixed
point, so they do not occur. The elements of order $23$, $21$, and
from one of the two conjugacy classes of elements of order $12$ have
indeed two cycles of the lengths as claimed.

\subsection{Classical Groups}\label{SS:ClGr} Suppose that $S$ is a
classical group. Our goal is to show that $\Si=\PSL_m(q)$, and that
except for a few small cases, the action is the natural one on the
projective space over $\FFF_q$. The main tool for doing that are good
upper bounds for element orders in automorphism groups of classical
groups.

\subsubsection{Element Orders in Classical Groups}
The following lemma controls the maximal possible orders of elements
in linear groups, if they are decorated with a field automorphism.

\begin{Lemma}\label{L:Lang:gamma}
  Let $q$ be a power of the prime $p$, $\overline{\FFF_p}$ be an algebraic
  closure of $\FFF_p$, and $G\le\GL_n(\overline{\FFF_p})$ be a connected
  linear algebraic group defined over $\FFF_p$. For $E$ a subfield of
  $\overline{\FFF_p}$, denote by $G(E)$ the group $G\cap\GL_n(E)$ of
  $E$--rational elements.
  
  Suppose that $E$ is finite, and let $\gamma\in\Aut(E)$. Then $G(E)$
  is normalized by $\gen{\gamma}$. Take $g=\gamma h$ in the semidirect
  product of $\gen{\gamma}$ with $G(E)$, where $h\in G(E)$. Let $f$ be
  the order of $\gamma$, and $F$ the fixed field in $E$ of $\gamma$.
  Then $g^f$ is conjugate in $G$ to an element in $G(F)$.
\end{Lemma}

\begin{proof} Clearly $\gen{\gamma}$ normalizes $G(E)$, as $G$ is
  defined over $\FFF_p$. We compute
\[
g^f=h^{\gamma^{f-1}}\cdots h^\gamma h,
\]
thus
\[
(g^f)^\gamma = hg^fh^{-1}.
\]
Extend $\gamma$ to $\overline{\FFF_p}$, and denote the induced action on
$G$ also by $\gamma$. By Lang's Theorem
(see~\cite[Theorem 10.1]{Steinberg:Mem}), the map $w\mapsto w^\gamma
w^{-1}$ from $G$ to $G$ is surjective. Thus there is $b\in G$ with
\[
h=b^\gamma b^{-1}.
\]
Therefore
\[
(b^{-1}g^fb)^\gamma = b^{-1}g^fb,
\]
so $b^{-1}g^fb$ is fixed under $\gamma$, hence contained in $G(F)$.
\end{proof}

In order to apply this lemma, we need the following easy estimate:

\begin{Lemma}\label{L:fq^r/f}
  Let $q,f,r$ be positive integers such that $2^f\le q$. Then $f\cdot
  q^{r/f}\le q^r$.
\end{Lemma}

\begin{proof}
  We have
\[
q^{r(1-1/f)}\ge 2^{r(f-1)}\ge 2^{f-1}\ge f,
\]
and the claim follows after multiplying with $q^{r/f}$.
\end{proof}

\begin{Lemma}\label{L:L:ind}
  Let $q$ be a power of the prime $p$. Let $\sigma\in\GL_n(q)$ act
  indecomposably on $V:=\FFF_q^n$. Then the order of $\sigma$ divides
  $p^b(q^u-1)$, where $u$ divides $n$, and $p^{b-1}\le n/u-1$ if
  $b>0$. Furthermore, $\sigma^{p^b(q^u-1)/(q-1)}$ is a scalar, and
  $p^b(q^u-1)\le q^n-1$. So in particular $\ord{\sigma}\le q^n-1$, and
  the order of the image of $\sigma$ in $\PGL_n(q)$ is at most
  $(q^n-1)/(q-1)$.
\end{Lemma}

\begin{proof}
  Write $\sigma=\sigma_{p'}\sigma_p$, where $\sigma_{p'}$ and
  $\sigma_p$ are the $p'$--prime part and $p$--part of $\sigma$,
  respectively. Let
\[
V=U_1\oplus U_2\oplus\dots\oplus U_m,
\]
be a decomposition into irreducible $\sigma_{p'}$--modules. Such a
decomposition exists by Maschke's Theorem.

Let $U$ be the sum of those $U_i$ which are $\sigma_{p'}$--isomorphic
to $U_1$. As $\sigma_p$ commutes with $\sigma_{p'}$, we get that
$U_i^{\sigma_p}$ is $\sigma_{p'}$--isomorphic to $U_i$ for each $i$.
By Jordan--H\"older, $U$ is a $\sigma$--invariant direct summand of
$V$. The indecomposability of $V$ with respect to $\sigma$ gives
$U=V$, so all $U_i$ are $\sigma_{p'}$--isomorphic.

Let $u$ be the common dimension of $U_i$, so $n=um$. By Schur's Lemma,
the restriction of $\sigma_{p'}$ to each $U_i$ can be identified with
an element of the multiplicative group of $\FFF_{q^u}$. As $\sigma$
commutes with $\sigma_{p'}$, we can consider $\sigma$ and $\sigma_{p}$
as elements in $\GL_m(q^u)$. So either $\sigma_{p}=1$, or
$p^b:=\ord{\sigma_p}\le p(m-1)$ by Lemma \ref{L:unip:ord}. Also, with
respect to this identification, $\sigma_{p'}$ is a diagonal matrix. So
$\sigma_{p'}^{(q^u-1)/(q-1)}$ acts as a scalar
$\lambda_i\in\FFF_q^\star$ on $U_i$. However, the $\lambda_i$ are
independent of $i$, because the $U_i$ are $\sigma_{p'}$--isomorphic.

To finish the claim, we need to show that $p^b(q^u-1)\le
q^{um}-1$. This is clear for $b=0$. For $b\ge1$, this follows from
$p^b\le p(m-1)$ and
\[
\frac{q^{um}-1}{q^u-1}=1+q^u+\dots+q^{u(m-1)}\ge 1+q^u(m-1)>p^b.
\]
(Note that $b\ge1$ implies $m>1$.)
\end{proof}

We obtain the following consequence

\begin{Proposition}\label{P:L:ord}
  Let $q$ be a prime power, and $n\ge2$.
  \begin{enumerate}
  \item\label{P:L:ord:gL} If $\sigma\in\gL_n(q)$, then
    $\ord{\sigma}\le q^n-1$.
  \item\label{P:L:ord:PgL} If $\overline{\sigma}\in\PgL_n(q)$, then
    $\ord{\overline{\sigma}}\le(q^n-1)/(q-1)$, except for
    $(n,q)=(2,4)$.
  \end{enumerate}
\end{Proposition}

\begin{proof}
  First assume that $\sigma\in\GL_n(q)$, and denote by
  $\overline{\sigma}$ the image of $\sigma$ in $\PGL_n(q)$. Let
  $\FFF_q^n=:V=V_1\oplus\dots\oplus V_r$ be a decomposition of $V$
  into $\sigma$--invariant and $\sigma$--indecomposable modules $V_i$.
  Let $n_i$ be the dimension of $V_i$. By Lemma \ref{L:L:ind}, the
  order of the restriction of $\sigma$ to $V_i$ divides
  $a_i:=p^{b_i}(q^{u_i}-1)$, where $u_i$ divides $n_i$, and $a_i\le
  q^{n_i}-1$. The order of $\sigma$ divides the least common multiple
  of the $a_i$. First suppose that $r>1$. Then $q-1$ divides each
  $a_i$, so
\begin{align*}
\ord{\sigma} &\le \lcm(a_1,\dots,a_r)\\
             &\le (a_1\cdots a_r)/(q-1)\\
             &\le (q^{n_1}-1)\cdots(q^{n_r}-1)/(q-1)\\
             &\le (q^n-1)/(q-1).
\end{align*}
If however $r=1$, then Lemma \ref{L:L:ind} applies directly. So in
either case, (a) and (b) hold for $\GL_n(q)$ and $\PGL_n(q)$,
respectively.

Now assume that $\sigma\in\gL_n(q)\setminus\GL_n(q)$, and let $f$ be
the smallest positive integer with $\sigma^f\in\GL_n(q)$. Note that
$f\ge2$. By Lemma \ref{L:Lang:gamma}, $\tau:=\sigma^f$ is conjugate to
an element $\tau'\in\GL_n(r)$, where $r:=q^{1/f}$. (We take the
natural inclusion $\GL_n(r)<\GL_n(q)$.) Part (a) is clear, as, by what
we saw already, $\ord{\sigma}\le f\ord{\sigma^f}< fr^n\le q^n$, where
we used Lemma \ref{L:fq^r/f} in the last step.

Part (b) requires a little more work. We have, similarly as above,
\[
\ord{\overline{\sigma}}\le f\frac{r^n-1}{r-1},
\]
and are done once we know that
\begin{equation*}
f\frac{r^n-1}{r-1}\le\frac{r^{nf}-1}{r^f-1}=\frac{q^n-1}{q-1}
\end{equation*}
which is equivalent to
\begin{equation}\label{qr}
f\frac{r^f-1}{r-1}\le\frac{r^{nf}-1}{r^n-1}.
\end{equation}
Note that $(x^f-1)/(x-1)=1+x+\dots+x^{f-1}$ is strongly monotonously
increasing for $x>1$, so inequality \eqref{qr} holds once it holds for
$n=2$. In this case, we have to show that $f\le(r^f+1)/(r+1)$. It is
easy to see that this last inequality holds except for $f=2$, $r=2$.
But then \eqref{qr} is equivalent to $6\le 2^n+1$, which is clearly
the case for $n\ge3$.
\end{proof}

\begin{Remark*}
  $\PgL_2(4)$ is indeed an exception for part (b) of the previous
  theorem. Note that $\PgL_2(4)\cong\sym_5$, so this group contains an
  element of order $6>5=(4^2-1)/(4-1)$.
\end{Remark*}

\begin{Lemma}\label{L:Bil:irr}
  Let $V$ be a vector space of dimension $n\ge2$ over $\FFF_q$ with a
  non--degenerate bilinear form $\kappa=(\cdot,\cdot)$. Let
  $\tau\in\Isom(V,\kappa)$ be an isometry with
  respect to this form, and assume that $\tau$ is irreducible on $V$.
  Then $n$ is even and the order of $\tau$ divides $q^{n/2}+1$.
\end{Lemma}

\begin{proof} By Schur's Lemma we have $V\cong\FFF_{q^n}$, and the
  action of $\tau$ induced on $\FFF_{q^n}$ is by multiplication with
  $\lambda\in\FFF_{q^n}^*$, where $\FFF_q[\lambda]=\FFF_{q^n}$. The
  eigenvalues of $\tau$ then are the powers $\lambda^{q^i}$ for
  $i=0,1,\dots,n-1$. Let $v_i\in V\otimes\FFF_{q^n}$ be an eigenvector
  to the eigenvalue $\lambda^{q^i}$. The form $(\cdot,\cdot)$ extends
  naturally to a non--degenerate form on $V\otimes\FFF_{q^n}$. Thus
  there exists $i$ with $(v_0,v_i)=c\ne0$. This gives
  $c=(v_0^\tau,v_i^\tau)=(\lambda
  v_0,\lambda^{q^i}v_i)=\lambda^{1+q^i}(v_0,v_i)=\lambda^{1+q^i}c$, so
  $\lambda^{1+q^i}=1$. Thus $\lambda\in\FFF_{q^{2i}}$, so $n\divides
  2i$. But $i<n$, hence $2i=n$, and the claim follows.\end{proof}

\begin{Lemma}\label{L:kappa:ind:p'}
  Let $V$ be a vector space over the finite field $F$ with a
  non--degenerate symmetric, skew--symmetric, or hermitian form
  $\kappa=(\cdot,\cdot)$. Write $F=\FFF_q$ if $\kappa$ is bilinear,
  and $F=\FFF_{q^2}$ if $\kappa$ is hermitian. Let
  $\sigma\in\Isom(V,\kappa)$ be an isometry with respect to $\kappa$.
  Suppose that $\sigma$ is semisimple and orthogonally indecomposable,
  but reducible on $V$. Then the following holds:
  
  $V=Z\oplus Z'$, where $Z$ and $Z'$ are $\sigma$--irreducible and
  totally isotropic spaces of the same dimension. Let $\Lambda$ and
  $\Lambda'$ be the set of eigenvalues of $\sigma$ on $Z$ and $Z'$,
  respectively. Then
\[
\Lambda'=
\begin{cases}
  \{\lambda^{-1}|\;\lambda\in\Lambda\} &\text{if $\kappa$ is
    bilinear,}\\
  \{\lambda^{-q}|\;\lambda\in\Lambda\} &\text{if $\kappa$ is
    hermitian.}
\end{cases}
\]

Furthermore, if $\kappa$ is not skew--symmetric, then $Z$ is not
$\sigma$--isomorphic to $Z'$.
\end{Lemma}

\begin{proof}
  Let $Z$ be a $\sigma$--invariant subspace of minimal positive
  dimension, in particular $Z$ is $\sigma$--irreducible. Also $Z^\bot$
  is $\sigma$--invariant. Furthermore, $Z$ is totally isotropic, for
  otherwise $V=Z\bot Z^\bot$ by irreducibility of $Z$.  As $\sigma$ is
  semisimple, there is a $\sigma$--invariant complement $Z'$ of
  $Z^\bot$ in $V$. From $\dim(Z')=\dim(V)-\dim(Z^\bot)=\dim(Z)$ and
  the minimality of $\dim(Z)$ we get that $Z'$ is
  $\sigma$--irreducible as well. We get $V=Z\oplus Z'$ once we know
  that $Z\oplus Z'$ is not degenerate. But this follows from
  \begin{align*}
    (Z\oplus Z')\cap(Z\oplus Z')^\bot &= (Z\oplus Z')\cap
    Z^\bot\cap(Z')^\bot \\
    &= Z\cap (Z')^\bot\\
    &= \{0\},
  \end{align*}
where the latter equality holds because $Z'$ is a complement to
$Z^\bot$, therefore $Z$ is not contained in $(Z')^\bot$.

Next we show the assertion about the eigenvalues if $\kappa$ is
bilinear. Let $\lambda$ be an eigenvalue of $\sigma$ with eigenvector
$v\in Z\otimes\overline{\FFF_q}$. Let $w\in
Z'\otimes\overline{\FFF_q}$ be such that $V\otimes\overline{\FFF_q}$
is the span of $w$ and $v^\bot$, and that $w$ is an eigenvector of
$\sigma$. Let $\mu$ be the corresponding eigenvalue. By construction,
$\rho:=(v,w)\ne0$, hence
\[
\rho=(v,w)=(v^\sigma,w^\sigma)=(\lambda v,\mu w)=\lambda\mu\rho,
\]
and the claim follows, as we can also switch the role of $Z$ and $Z'$
in this argument.

The case that $\kappa$ is hermitian is completely analogous.

Finally, suppose that $\kappa$ is not skew--symmetric, and assume in
contrary that there is a $\sigma$--isomorphism $\phi:\;Z\mapsto Z'$.
Let $R=\FFF_q[\sigma]\le\End(V)$ be the algebra generated by $\sigma$.
As $\kappa$ is not skew--symmetric, there is an element $v\in V$ with
$(v,v)\ne0$.  Write $v=z+z'$ with $z$ and $z'$ in $Z$ and $Z'$,
respectively.  Clearly $z$ and $z'$ are non--zero. By Schur's Lemma,
$R$ acts sharply transitively on the non--zero elements of $Z'$, in
particular, there is $\rho\in R$ such that $(z^\phi)^\rho=z'$. Let
$\psi:\;Z\mapsto V$ be the homomorphism defined by
$w^\psi:=w+(w^\phi)^\rho$. This map is clearly injective, $\psi$
commutes with $\sigma$, so the image $Z^\psi$ has the same dimension
as $Z$, and of course is $\sigma$--irreducible as well. By
construction, the element $v=z^\psi$ is not isotropic, so $Z^\psi$ is
not totally isotropic, thus $\kappa$ restricted to $Z^\psi$ is not
degenerate. We get $V=Z^\psi\bot(Z^\psi)^\bot$, contrary to
indecomposability.
\end{proof}

\begin{Remark*} Let $V$ be $2$--dimensional with a non--degenerate
  skew--symmetric form, and $\sigma$ the identity map. As $V$ is
  clearly not the orthogonal sum of two $1$--dimensional spaces, we
  cannot dispense of the assumption that $\kappa$ is not skew--symmetric in
  the last part of the lemma.
\end{Remark*}

We now extend the previous lemma to those $\sigma$ which are not
necessarily semisimple.

\begin{Lemma}\label{L:kappa:ind} Let $V$ be a vector space over
  $\FFF_q$ with a non--degenerate symmetric, skew--symmetric, or
  hermitian form $\kappa=(\cdot,\cdot)$. Let
  $\sigma\in\Isom(V,\kappa)$ be an isometry with respect to this form.
  Assume that $\sigma$ is orthogonally indecomposable, but reducible
  on $V$. Denote by $\sigma_{p'}$ the $p'$--part of $\sigma$. Then the
  following holds:
\[
V=(U_1\bot U_2\bot\dots\bot U_r)\bot((Z_1\oplus
Z_1')\bot\dots\bot(Z_s\oplus Z_s')),
\]
where the $U_i$, $Z_i$ and $Z_i'$ are $\sigma_{p'}$--irreducible, the
$U_i$ and $(Z_i\oplus Z_i')$ are not degenerate, the $Z_i$ and $Z_i'$
are totally isotropic and the $U_i$, $Z_i$ and $Z_i'$ have all the
same dimension. Also, $r+2s\ge2$.
\end{Lemma}

\begin{proof}
  Choose an orthogonal decomposition of $V$ into non--trivial
  $\sigma_{p'}$--invariant subspaces of maximal length, so these
  subspaces do not decompose orthogonally into smaller
  $\sigma_{p'}$--invariant spaces. Let the $U_i$ be those subspaces
  which are $\sigma_{p'}$--irreducible, and let the $(Z_i\oplus Z_i')$
  be the remaining ones according to the previous lemma.
  
  The $\sigma_{p'}$--homogeneous components $H_1, H_2,\dots$ are
  $\sigma$--invariant as a consequence of Jordan--H\"older. Let $H$ be
  the sum of those $H_k$ where the irreducible summands of $H_k$ have
  the same dimension as those of $H_1$. Then $Z_i$ appears in $H$ if
  and only if $Z_i'$ appears in $H$. The orthogonal indecomposability
  of $\sigma$ forces $H=V$.
  
  Suppose that $r+2s<2$. Then $s=0$ and $r=1$, that is $\sigma$ is
  irreducible on $V=U_1$, a contradiction.
\end{proof}

\begin{Lemma}\label{L:q^mi+1}
  Let $q\ge2$ and $m_1,m_2,\dots,m_\rho$ be distinct positive integers
  with sum $m$. Then
\[
\prod_{i=1}^\rho(q^{m_i}+1)\le e^{1/(q-1)}q^m.
\]
\end{Lemma}

\begin{proof} For $x$ real we have $1+x\le e^x$. Substitute
  $x=1/q^{m_i}$ and multiply by $q^{m_i}$ to obtain
\[
q^{m_i}+1\le q^{m_i}e^{1/q^{m_i}}.
\]
Multiply these inequalities for $i=1,2,\dots,\rho$ to obtain
\[
\prod(q^{m_i}+1)\le q^me^\Sigma,
\]
with
\[
\Sigma=\sum_{i=1}^\rho\frac{1}{q^{m_i}}
\le\sum_{k=1}^\infty\frac{1}{q^k}=\frac{1}{q-1},
\]
as the $m_i$ are distinct. The claim follows.\end{proof}

\begin{Lemma}\label{L:Bil:ind:ord}
  Use the notation from Lemma \ref{L:kappa:ind} with $\kappa$
  bilinear, and let $z$ be the common dimension of the spaces $Z_i$,
  $Z_i'$, $U_i$. Set $w:=r+2s$, thus $v:=\dim(V)=wz$. Then there is a
  non--negative integer $b$, such that $\ord{\sigma}$ divides
  $p^b(q^z-1)$. Furthermore,
\[
\ord{\sigma} \le
\begin{cases}
2q^{[v/2]}                &   \text{in any case,}\\
q^{[v/2]}                 &   \text{if $\ord{\sigma}$ is odd,}\\
q^{[v/2]}                 &   \text{if $q$ is even, and
  $(q,w,z)\ne(2,2,2)$ or $(2,3,2)$,}
\end{cases}
\]
If $q=2$ and $v=4$ or $6$ and $\ord{\sigma}>2^{v/2}$, then
$\ord{\sigma}=6$ if $v=4$, and $\ord{\sigma}=12$ if $v=6$.
\end{Lemma}

\begin{proof} As the spaces $Z_i$, $Z_i'$, and $U_i$ are all
  $\sigma_{p'}$--irreducible of dimension $z$, it follows that the
  order of $\sigma_{p'}$ divides $q^z-1$. Let $p^b$ be the order of
  the $p$--part of $\sigma$. As $w\ge2$, hence $z\le[v/2]$, the stated
  inequalities clearly hold for $b=0$. Thus assume $b\ge1$ from now
  on.
  
  First assume $p>2$. We are clearly done except if
\begin{equation}\label{qwz}
p^b(q^z-1) > 2q^{[wz/2]}.
\end{equation}
From \eqref{qwz} we obtain
\[
p^bq^z > 2q^{[wz/2]}.
\]
As each factor except $2$ is divisible by $p$, we obtain from that
even sharper
\[
p^bq^z \ge pq^{[wz/2]},
\]
hence
\begin{equation}\label{qwz2}
p^{b-1}q^z\ge q^{[wz/2]}.
\end{equation}
Let $w'$ be the number of elements in a maximal subset of the summands
$Z_i$, $Z_i'$, and $U_i$ which are pairwise $\sigma_{p'}$--isomorphic.
Then the restriction of $\sigma_p$ to the sum of these spaces can be
seen as an element in $\GL_{w'}(q^z)$, so the order of this
restriction is bounded by $p(w'-1)$, see Lemma \ref{L:unip:ord}.
Clearly $w'\le w$, hence $p^{b-1}\le w-1$. So with \eqref{qwz2} we
obtain further
\[
w-1 \ge q^{[wz/2]-z}.
\]
We first contend that $w\le5$, and that $z=1$ if $w>2$. For suppose
$z\ge2$. Then $[wz/2]-z\ge w-2$, as $w\ge2$. So $w-1\ge q^{w-2}$,
which gives $w=2$. Is is easy to see that $w-1\ge q^{[w/2]-1}$ gives
$w\le5$. Suppose $w=4$ or $5$. We obtain $q=3$. Furthermore, $b\le2$,
so $b=2$ for otherwise we are done (check \eqref{qwz2}). As $V$
decomposes into $1$--dimensional eigenspaces for $\sigma_{3'}$, the
eigenvalues are in $\FFF_3\setminus\{0\}$, so we have that
$\ord{\sigma_{3'}}$ is at most $2$, hence the order of $\sigma$ is at
most $2\cdot3^2=18$, the exact bound we wanted to prove (and which is
sharp indeed).

Now suppose $w=3$. Clearly $b=1$. We have either $r=3$, $s=0$, or
$r=1$, $s=1$. In the first case $\sigma_{p'}$ restricts to an element
of order at most $2$ on each $U_i$, so the order of $\sigma$ divides
$2p$, and the claim follows. Thus assume $r=1$, $s=1$. Let $\lambda$
be the eigenvalue of $\sigma_{p'}$ on $U_1$. Clearly $\lambda=\pm1$.
Also, $\lambda$ is an eigenvalue on $Z_1$ or $Z_1'$, for otherwise
$U_1$ were $\sigma$--invariant, contrary to orthogonal irreducibility.
By Lemma \ref{L:kappa:ind:p'} the eigenvalues on $Z$ and $Z'$ then are
$\pm1$, so the order of $\sigma_{p'}$ is at most $2$, and we are done
again.

Finally, we have to look at $w=2$. Here we have not necessarily $z=1$.
First suppose that $s=0$, that is $V=U_1\oplus U_2$. The order of
$\sigma_{p'}$ on $V$ divides $q^{[z/2]}+1$. The claim follows as
$p(q^{[z/2]}+1)\le 2q^z=2q^{[v/2]}$. Thus suppose that $r=0$, $s=1$,
so $V=Z_1\oplus Z_1'$. Let $\lambda\in\FFF_{q^z}$ be an eigenvalue of
$\sigma_{p'}$ on $Z_1$. By irreducibility, the eigenvalues of
$\sigma_{p'}$ on $Z_1$ are $\lambda^{q^i}$ for $i=0,1,\dots,z-1$. By
Lemma \ref{L:kappa:ind:p'}, the inverses of these eigenvalues are the
eigenvalues of $\sigma_{p'}$ on $Z_1'$. We contend that these two sets
are the same.  Namely as $\sigma$ is not semisimple, it cannot leave
invariant both $Z_1$ and $Z_1'$. So without loss $Z_1^{\sigma_p}\ne
Z_1$, and we obtain that $Z_1$ and $Z_1'$ are
$\sigma_{p'}$--isomorphic by Jordan--H\"older. So the set of
eigenvalues on $Z_1$ is closed under inversion, in particular there is
an $i$ such that $\lambda^{-1}=\lambda^{q^i}$. This gives
$\lambda^{q^{2i}-1}=1$, so $\lambda\in\FFF_{q^{2i}}$. We obtain that
$z$ divides $2i<2z$, as $\FFF_{q^z}=\FFF_q[\lambda]$. If $i=0$, then
$\lambda=\pm1$, so $\sigma_{p'}$ has order at most $2$, and the claim
clearly follows, as $b=1$. If $i>0$, then $z=2i$, so the order of
$\sigma_{p'}$ divides $q^{z/2}+1$, and the claim follows again from
$(q^{z/2}+1)p < 2q^{z/2}q\le 2q^z=2q^{[v/2]}$.

We are left to look at the case $p=2$. As the form is not degenerate,
we have necessarily $v=wz$ even. We proceed similarly as above. Recall 
that $b\ge1$. We are done unless
\begin{equation}
  \label{wz2}
  2^b(q^z-1) > q^{wz/2}.
\end{equation}
From that we obtain
\[
2^b > q^{wz/2-z},
\]
hence
\[
2^{b-1} \ge q^{wz/2-z}
\]
and
\begin{equation}
  \label{wz22}
  w-1 \ge q^{wz/2-z},
\end{equation}
as $2^{b-1}\le w-1$. If $z\ge2$, then $w-1\ge q^{w-2}$, hence either
$w=3$, $q=2$, $z=2$; or $w=2$. The first case gives $2^{b-1}\le
w-1=2$, so $b\le2$, hence $\ord{\sigma}=12$ or $\le6<2^3=2^{v/2}$.

Thus we have $z=1$ except possibly for $w=2$. First assume $w>2$, so
$w\ge4$ is even. We obtain $w\le6$ from \eqref{wz22}. Suppose $w=6$.
Then $q=2$ and $b\le3$, and we obtain a contradiction to \eqref{wz2}.
Next suppose $w=4$. Again $q=2$. From \eqref{wz2} we obtain $2^b>2^2$,
hence $b\ge3$, a contradiction to $4\le 2^{b-1}\le w-1=3$.

Finally, suppose $w=2$. Clearly $b=1$. The argument from the last
paragraph in the case $p>2$ shows that the critical case is when $z$
is even and $\ord{\sigma}$ divides $2(q^{z/2}+1)$. Now
\[
2(q^{z/2}+1) = q^z - ((q^{z/2}-1)^2-3) \le q^z = q^{[v/2]}
\]
except for $q=2$, $z=2$.
\end{proof}

\begin{Proposition}\label{P:Bil:ord}
  Let $\sigma\in\GL_n(q)$ be an isometry with respect to a
  non--degenerate skew--symmetric or symmetric bilinear form on
  $\FFF_q^n$. Then
\[
\ord{\sigma}\le
\begin{cases}
            2q^{[n/2]}  &  \text{if $q$ is odd,}\\
             q^{[n/2]}  &  \text{if $q$ and $\ord{\sigma}$ are odd,}\\
  e^{1/(q-1)}q^{[n/2]}<2q^{[n/2]}  &  \text{if $q\ne2$ is even,}\\
  (3e/2)2^{[n/2]}       &  \text{if $q=2$.}
\end{cases}
\]
\end{Proposition}

\begin{proof}
Choose a decomposition of $V$ into orthogonally indecomposable
$\sigma$--invariant subspaces. The order of $\sigma$ is the least
common multiple of the orders of the restriction of $\sigma$ to these
subspaces. Lemmas \ref{L:Bil:irr} and \ref{L:Bil:ind:ord} give upper
bounds for these orders.

In the following we use several times the trivial inequality
\[
[u_1/2]+[u_2/2]+\dots+[u_k/2]\le[(u_1+u_2+\dots+u_k)/2]
\]
for integers $u_i$.

First suppose that $q$ is odd. Let $U$ be such a subspace of dimension
$u$. If $U$ is $\sigma$--irreducible, then $\ord{\sigma|_U}$ is at
most $q^{[u/2]}+1$, so the order is at most $(q^{[u/2]}+1)/2\le
q^{[u/2]}$ if $\ord{\sigma|_U}$ is odd, and at most $q^{[u/2]}+1\le
2q^{[u/2]}$ otherwise. The assertion follows if $U=V$. So suppose
$U<V$. By induction, the stated bound holds for the restriction of
$\sigma$ to $U^\bot$. Let $\overline{u}=\dim(U^\bot)$. If the orders
of the restriction of $\sigma$ to $U$ and $U^\bot$ are relatively
prime, then at least one of the orders is odd, and we obtain the claim
by multiplying the corresponding upper bounds. If these orders are not
relatively prime, then the product of these orders divided by $2$ is
an upper bound for the order of $\sigma$, so the claim holds as well.

Now suppose that $q$ is even. Let $W_i$ be those subspaces from above
on which $\sigma$ acts irreducibly, and let $W$ be the sum of these
spaces. Set $w:=\dim(W)$, and let $1<w_1<w_2<\dots$ be the distinct
dimensions of the spaces $W_i$. Note that if $\dim(W_i)=1$, then the
restriction of $\sigma$ to $W_i$ is trivial. By Lemma \ref{L:Bil:irr},
the $w_i$ are even, and the order of the restriction of $\sigma$ to
the associated space divides $q^{w_i/2}+1$. Thus the order of
$\sigma|_W$ divides the product of the $q^{w_i/2}+1$. This product is
less than $e^{1/(q-1)}q^{[w/2]}$ by Lemma \ref{L:q^mi+1}. If $q\ne2$,
then apply the bounds in Lemma \ref{L:Bil:ind:ord} to the summands of
$W^\bot$ to get the claim. Finally suppose $q=2$. We are done except
if one of the summands $Q$ of $W^\bot$ has dimension $4$ or $6$, and
$\sigma|_Q$ has order $6$ or $12$, respectively. The stated inequality
then holds for $W\bot Q$. If there are more such summands $Q'$ in
$W^\bot$, then they do contribute at most by a factor
$2<2^{[\dim(Q')/2]}$ to the order of $\sigma$. All other summands of
dimension $r$ contribute by a factor of at most $2^{[r/2]}$, so the
claim follows.
\end{proof}

At a few places we need the following trivial

\begin{Lemma}\label{L:q^m+-1}
  Let $1\le i<m$ and $q\ge2$ be integers. Let $\varepsilon$ be $-1$ or
  $1$. Then
\[
\frac{(q^m+\varepsilon)(q^{m-1}-\varepsilon)}{q^i-1}>q^{2m-1-i}.
\]
\end{Lemma}

\begin{proof}
  Clearly $q^{m-i}-1\ge\varepsilon(q-1)$. Multiply by $q^{m-1}$ to get
  $q^{2m-1-i}-q^{m-1}\ge\varepsilon(q^m-q^{m-1})$, hence
  $q^{2m-1-i}-1>\varepsilon(q^m-q^{m-1})$. But this inequality is
  equivalent to the stated one.
\end{proof}

As before denote by $\mu(\Si)$ and $o(\Si)$ a lower bound for the
degree of a faithful permutation representation and an upper bound for
the order of an element, respectively. The minimal permutation degrees
$\mu(\Si)$ have been determined by Cooperstein and Patton -- we use
the ``corrected'' list \cite[Theorem 5.2.2]{KL:Sub} which still
contains a mistake (giving the wrong $\mu$ for $P\Omega_{2m}^+(3)$).
We exclude the group $\PSL_2(5)$, as $\PSL_2(5)\cong\alt_5$, a case we
already dealt with.  Besides that, the list \cite[Theorem
5.2.2]{KL:Sub} contains a few duplications.  Accordingly, we drop
$\PSp_4(3)$ in view of $\PSp_4(3)\cong\PSU_4(2)$ and $\Spl_4(2)'$ in
view of $\Spl_4(2)'\cong\PSL_2(9)$.

Assume that the almost simple group $A$ acts primitively and contains
an element with at most two cycles. We consider the case that the
minimal normal subgroup $\Si$ of $A$ is a classical group. The aim of
this section is to show that $\Si$ is isomorphic $\PSL_m(q)$, a case
to be handled afterwards.

\subsubsection{$\Si$ Symplectic}

\begin{Lemma}\label{L:Sp:ord}
  Let $\Si=\PSp_{2m}(q)$ be the simple symplectic group, and
  $\sigma\in\Aut(\Si)$. Then
\[
\ord{\sigma} \le
\begin{cases}
4q^m                  & \text{if $q$ is odd,} \\
e^{1/(q-1)}q^m        & \text{if $q\ne2$ is even, $m\ge3$,}\\
2e^{1/(q-1)}q^2       & \text{if $q\ne2$ is even, $m=2$,}\\
(3e/2)2^m             & \text{if $q=2$, $m\ge3$.}
\end{cases}
\]
In particular, $\ord{\sigma}\le4q^m$ if $q\ne2$.
\end{Lemma}

\begin{proof}
  Let $q=p^r$ with $p$ a prime. If $q$ is odd, then
  $\Out(\Si)=C_2\times C_r$, see \cite[Theorem~2.1.4,
  Prop.~2.4.4]{KL:Sub}, where $C_r$ comes from a field automorphism.
  Thus $\sigma^2$ has a preimage $\tau$ in
  $\Spl_{2m}(q)\rtimes\Aut(\FFF_q)$. Let $f$ be the order of the
  associated field automorphism. By Lemma \ref{L:Lang:gamma}, $\tau^f$
  is conjugate to an element in the group $\Spl_{2m}(q^{1/f})$, whose
  element orders are bound by $2q^{m/f}$ by Proposition
  \ref{P:Bil:ord}. Thus $\tau$ has order at most $2fq^{m/f}\le2q^m$,
  where we used Lemma \ref{L:fq^r/f}. The claim follows in the odd
  case.
  
  If $q$ is even, then $\Out(\Si)=C_r$ if $m\ge3$. Argue as above. If
  $m=2$, then $\Out(\Si)$ is cyclic of order $2r$, and the square of a
  generator is a field automorphism, see \cite[Chapter
  12]{Carter:Simple}. The claim follows as above.
\end{proof}

Now we rule out the symplectic groups in the order as they appear in
Table \ref{Tb:ClGr} on page \pageref{Tb:ClGr}.

$\mathbf{m\ge2,\;q\ge3,\;(m,q)\ne(2,3).}$ Let $\sigma\in\Aut(\Si)$. The
minimal faithful permutation degree of $\Si$ is $(q^{2m}-1)/(q-1)$. As
$q\ge3$, we get $\ord{\sigma}\le4q^m$ by the previous Lemma. So Lemma
\ref{L:kl:Deg-Ord} gives
\[
\frac{q^{2m}-1}{q-1}\le2\ord{\sigma}\le 2\cdot4\cdot q^m.
\]
Note that the left hand side is bigger than $q^{2m-1}$, so it follows
that $q^{m-1}<8$. Thus $m=2$ and $q\le7$. But
$(7^4-1)/(7-1)=400>392=8\cdot 7^2$, so $q=7$ is out. Thus $q=4$ or
$5$. But $\ord{\sigma}\le20$ for $q=4$, and $\ord{\sigma}\le30$ for
$q=5$, see the atlas \cite{ATLAS}. These improved bounds contradict
the above inequality.

$\mathbf{m\ge3,\;q=2.}$ We get $\mu(\Si)=2^{m-1}(2^m-1)\le 2(3e/2)2^m$,
hence $2^m-1\le6e$, so $m=3$ or $4$. If $m=4$, then the atlas gives
$\ord{\sigma}\le30$, contrary to $\mu(\Si)\le2\ord{\sigma}$. Thus $m=3$.
The atlas gives $\ord{\sigma}\le15$, and the next biggest element
order is $12$. Also, there is a maximal subgroup of index $28$, and
the next smallest has index $36$. So $\ord{\sigma}=15$ and $n=28$. But
$15=\lcm(k,28-k)$ has no solution, therefore $\sigma$ must have more
than $2$ cycles in this representation.

\subsubsection{$\Si$ Orthogonal in Odd Dimension}

\begin{Lemma}\label{L:O2m+1:ord}
  Let $\Si=\Omega_{2m+1}(q)$ be the simple orthogonal group with $q$
  odd, $m\ge3$, and $\sigma\in\Aut(\Si)$. Then
\[
\ord{\sigma} \le 2q^m.
\]
\end{Lemma}

\begin{proof}
  Set $V=\FFF_q^{2m+1}$, $\overline{V}=V\otimes\overline{\FFF_q}$, and
  let $\kappa$ be the standard bilinear form on $\overline{V}$. The
  algebraic group $G:=\SL(\overline{V})\cap\Isom(\overline{V},\kappa)$
  is connected. Let $\sigma$ be in $\Aut(\Si)$. By the structure of
  the automorphism group of $\Si$ (see \cite[Prop.~2.6.3]{KL:Sub}), we
  find a preimage $\tau$ of $\sigma$ in
  $\Isom(V,\kappa)\rtimes\Aut(\FFF_q)$.  As $\Isom(V,\kappa)$ is an
  extension of $G(\FFF_q)$ by the scalar $-1$, we may assume that
  $\tau\in G(\FFF_q)\rtimes\Aut(\FFF_q)$. Now use Lemma
  \ref{L:Lang:gamma} together with Proposition \ref{P:Bil:ord} and
  Lemma \ref{L:fq^r/f} to get the conclusion.
\end{proof}

$\mathbf{m\ge3,\;q\ge5\text{ odd}.}$ We get a stronger inequality as
in the previous case $\Si=\PSp_{2m}(q)$, where we saw that there is no
solution for $m\ge3$.

$\mathbf{m\ge3,\;q=3.}$ We get $3^m(3^m-1)/2\le 2\cdot2\cdot3^m$,
hence $3^m\le9$, so $m\le2$, a contradiction.

\subsubsection{$\Si$ Orthogonal of Plus Type}

\begin{Lemma}\label{L:O+2m:ord}
  Let $\Si=P\Omega^+_{2m}(q)$ be the simple orthogonal group with Witt
  defect $0$, and $\sigma\in\Aut(\Si)$. Write $q=p^f$ for $p$ a prime.
  Then
\[
\ord{\sigma} \le
\begin{cases}
4fq^m\le2q^{m+1}          & \text{if $q$ is odd, $m\ge5$,}\\
8fq^4\le4q^5              & \text{if $q$ is odd, $m=4$,}\\
2fq^m\le q^{m+1}          & \text{if $q\ne2$ is even, $m\ge5$,}\\
(9/2)fq^4\le(9/4)q^5      & \text{if $q\ne2$ is even, $m=4$,}\\
(3e/2)2^m                 & \text{if $q=2$, $m\ge5$,}\\
30                        & \text{if $q=2$, $m=4$}.
\end{cases}
\]
\end{Lemma}

\begin{proof}
  Let $\kappa$ be the bilinear form associated to $\Si$. First suppose
  that $m\ge5$. Assume $q$ odd first. Then $\sigma^{2f}$ has a
  preimage in $\Isom(\FFF_q^{2m},\kappa)$, this follows from the
  structure of the automorphism group of $\Si$, see \cite[Theorem 2.1.4,
  Table~2.1.D]{KL:Sub}. Now apply Proposition \ref{P:Bil:ord}, and
  note that $2f\le q$. If $q$ is even, then $\sigma^f$ already has a
  preimage in $\Isom(\FFF_q^{2m},\kappa)$, hence if $q\ne2$, then
  $\ord{\sigma}\le fe^{1/(q-1)}q^m<2fq^m\le q^{m+1}$ by Proposition
  \ref{P:Bil:ord}, or $\ord{\sigma}\le(3e/2)2^m$ if $q=2$.
  
  Now suppose that $m=4$. We have
  $\Out(P\Omega^+_8(q))\cong\sym_3\times C_f$ if $q$ is even, and
  $\cong\sym_4\times C_f$ if $q$ is odd, see \cite[p.38]{KL:Sub}. Thus
  if $q$ is odd, then either $\sigma^{3f}$ or $\sigma^{4f}$ has a
  preimage in $\Isom(\FFF_q^{2m},\kappa)$, so $\ord{\sigma}$ is at
  most $4f$ times the maximal order of an element in
  $\Isom(\FFF_q^{2m},\kappa)$, and we use Proposition \ref{P:Bil:ord}
  again. If $q\ne2$ is even, then analogously
  $\ord{\sigma}\le3fe^{1/(q-1)}q^4\le3e^{1/3}fq^4<(9/2)fq^4$. If
  $q=2$, then use the atlas information \cite{ATLAS}.
\end{proof}

$\mathbf{m\ge4,\;q\ge4.}$ First suppose that $m\ge5$. We get
\[
\frac{(q^m-1)(q^{m-1}+1)}{q-1}\le2\ord{\sigma}\le 4q^{m+1}.
\]
The left hand side is bigger than $q^{2m-2}$ by Lemma \ref{L:q^m+-1},
so we obtain further $q^{2m-2}<4q^{m+1}$, hence $q^2\le q^{m-3}<4$, a
contradiction.

Next assume $m=4$. First assume $q$ odd. Similarly as above we obtain
$q^6<16fq^4\le 8q^5$. Note that if $f=1$, then $q<4$, a contradiction.
Thus assume $f\ge2$. We obtain $q<8$, so $f=2$, hence $q^2<32$, thus
$q\le5$, giving the contradiction $f=1$.

Now assume that $q\ne2$ is even. We obtain $q^6<2\cdot(9/4)q^5$, hence 
$q=4$. But $\ord{\sigma}\le(9/4)4^5=2304$, whereas
$\mu(\Si)=5525>2\cdot2304$, a contradiction.

$\mathbf{m\ge4,\;q=3.}$ First consider $m=4$. One verifies that
$o(P\Omega^+_{8}(3))=40$, so $\ord{\sigma}\le4\cdot40=160$, because
$\Out(P\Omega^+_{8}(3))=\sym_4$. In view of $\mu(S)=1080$ this case is
out. Suppose $m\ge5$. We obtain $3^{m-1}(3^m-1)/2\le
2\cdot2\cdot3^{m+1}$, hence $m<5$, a contradiction.

$\mathbf{m\ge4,\;q=2.}$ If $m=4$, then $\ord{\sigma}\le30$, whereas
$\mu(\Si)=120$, so this case is out. Suppose $m\ge5$. We obtain
$2^{m-1}(2^m-1)\le 2\cdot(3e/2)2^m$, hence $2^m\le 6e+1=17.3\dots$,
thus $m\le4$, a contradiction.

\subsubsection{$\Si$ Orthogonal of Minus Type}

\begin{Lemma}\label{L:O-2m:ord}
  Let $\Si=P\Omega^-_{2m}(q)$ be the simple orthogonal group with Witt
  defect $1$, and $\sigma\in\Aut(\Si)$. Write $q=p^f$ for $p$ a prime.
  Then
\[
\ord{\sigma} \le
\begin{cases}
4fq^m\le2q^{m+1}          & \text{if $q$ is odd, $m\ge4$,}\\
2fq^m\le q^{m+1}          & \text{if $q\ne2$ is even, $m\ge4$,}\\
(3e/2)2^m                 & \text{if $q=2$, $m\ge4$,}\\
30                        & \text{if $q=2$, $m=4$,}\\
60                        & \text{if $q=2$, $m=5$.}
\end{cases}
\]
\end{Lemma}

\begin{proof}
  The proof follows exactly as in Lemma \ref{L:O+2m:ord}, except that
  for $m=4$, there is no exceptional (graph) automorphism of order
  $3$. For $q=2$ and $m=4$ or $5$ use the atlas \cite{ATLAS}.
\end{proof}

Now $\Si=P\Omega^-_{2m}(q)$ for $m\ge4$. From Lemma \ref{L:q^m+-1} we
get $\mu(\Si)>q^{2m-2}$. First suppose $q\ne2$. We obtain
$q^{2m-2}<2\cdot2q^{m+1}$, hence $q^{m-3}<4$. Thus $m=4$ and $q=3$.
But this contradicts the sharper bound $\ord{\sigma}\le4\cdot3^4=324$.
If $q=2$, then $2^{2m-2}<2\cdot(3e/2)2^m$, hence
$2^{m-2}\le3e=8.1\dots$, so $m\le5$. Arrive at a contradiction using
the upper bounds for $\ord{\sigma}$ from Lemma \ref{L:O-2m:ord}.

\subsubsection{$\Si$ Unitary} 

\begin{Lemma}\label{L:GU:irr}
  Suppose that $\sigma\in\GU_n(q)$ acts irreducibly on $\FFF_q^n$.
  Then $n$ is odd, and $\ord{\sigma}$ divides $q^n+1$. The order of
  the image of $\sigma$ in $\PGU_n(q)$ divides $(q^n+1)/(q+1)$.
\end{Lemma}

\begin{proof}
  Let $\lambda$ be an eigenvalue of $\sigma$. Then
  $\FFF_{q^2}[\lambda]=\FFF_{q^{2n}}$. All the eigenvalues of $\sigma$
  are $\lambda^{q^{2i}}$ with $i=1,\cdots,n$. Similarly as in the
  proof of Lemma \ref{L:Bil:irr}, there exists an index $i$ in the
  given range such that $\lambda^{-q}=\lambda^{q^{2i}}$, so
\begin{equation}\label{2i+1}
\lambda^{q^{2i-1}+1}=1.
\end{equation}
It follows that $\lambda^{q^{4i-2}-1}=1$, so
$\lambda\in\FFF_q^{4i-2}$. Therefore $n\divides 2i-1<2n$, so $n=2i-1$.
The assertion about the order of $\sigma$ follows from \eqref{2i+1}.
By the irreducibility, the element $\sigma$ is a subgroup of a Singer
group of order $q^{2n}-1$ on $\FFF_{q^2}^n$. The (unique) subgroup of
order $q+1$ of this Singer group consists of scalars, because $q+1$
divides $q^2-1$. Also, $q+1$ divides $q^n+1$, so modulo scalars
$\sigma$ has order at most $(q^n+1)/(q+1)$.
\end{proof}

\begin{Lemma}\label{L:GU:n<5:ord}
  Let $\sigma\in\GU_n(q)$, and denote by $\overline{\sigma}$ the image of
  $\sigma$ in $\PGU_n(q)$. Let $q=p^f$ with $p$ prime. The following
  holds.
\begin{enumerate}
\item If $n=1$, then $\ord{\sigma}$ divides $q+1$.
\item If $n=2$, then $\ord{\sigma}$ divides $q^2-1$ or $p(q+1)$.
\item If $n=3$, then $\ord{\sigma}$ divides $q^3+1$, $q^2-1$, or
  $p^r(q+1)$ with $r\le2$ and $r=1$ if $p>2$. Furthermore,
  $\ord{\overline{\sigma}}$ divides $q^2-q+1$, $q^2-1$ or $p(q+1)$. For
  $p=2$, there is the additional possibility $\ord{\overline{\sigma}}=4$.
\item If $n=4$, then $\ord{\overline{\sigma}}$ divides $q^3+1$,
  $q^3-q^2+q-1$, or $p^r(q^2-1)$ where $r\le2$ and $r=1$ if $p>2$. For
  $p=3$, there is the additional possibility $\ord{\overline{\sigma}}=9$.
\end{enumerate}
\end{Lemma}

\begin{proof} Denote by $\sigma_{p'}$ the $p'$--part of $\sigma$. Set
  $F=\FFF_{q^2}$, so $\GU_n(q)$ is the isometry group of the unique
  hermitian from on $F^n$.
  
  The case $n=1$ is trivial.
  
  Suppose that $n=2$. By Lemma \ref{L:GU:irr}, $\sigma$ is reducible
  on $V=F^n$. If $\sigma$ is semisimple, then the eigenvalues of
  $\sigma$ are in $F$, so $\ord{\sigma}\divides q^2-1$. If $\sigma$ is
  not semisimple, then $\sigma_{p'}$ is the centralizer of an element
  of order $p$, hence $\sigma_{p'}$ is a scalar, and the claim follows
  again.
  
  Now assume $n=3$. If $\sigma$ is irreducible, then apply Lemma
  \ref{L:GU:irr}. If $\sigma$ is orthogonally decomposable, then apply
  (a) and (b) to get that $\ord{\sigma}$ divides $q^2-1$ or $p(q+1)$.
  Next assume $\sigma$ reducible, but orthogonally indecomposable.
  Choose a maximal orthogonal decomposition of $V$ in
  $\sigma_{p'}$--invariant subspaces. By Lemma \ref{L:kappa:ind:p'}
  and the notation from there, either $V=U_1\bot U_2\bot U_3$, or
  $V=U_1\bot(Z_1\oplus Z_1')$. Assume the first possibility. By
  orthogonal irreducibility of $\sigma$, the $U_i$ are pairwise
  $\sigma_{p'}$--isomorphic, thus $\sigma_{p'}$ is a scalar on $V$,
  with order dividing $q+1$. Let $p^r$ be the order of the $p$--part
  of $\sigma$. Then $p^{r-1}\le2$ by Lemma \ref{L:unip:ord}, and we
  get the divisibilities as stated. If we have the latter orthogonal
  decomposition, then $U_1$ must be $\sigma_{p'}$--isomorphic to $Z_1$
  or $Z_1'$, say to $Z_1$.  On the other hand, $Z_1$ and $Z_1'$ are
  not $\sigma_{p'}$--isomorphic by Lemma \ref{L:kappa:ind:p'}. We get
  that $\sigma_p$ leaves invariant $U_1\bot Z_1$ and $Z_2$, thus the
  order of $\sigma_p$ divides $p$. The order of $\sigma_{p'}$ divides
  $q+1$, because the restriction to $U_1$ satisfies this, so this
  holds also for the restriction to $Z_1$, and then also for the
  restriction to $Z'_1$ by Lemma \ref{L:kappa:ind:p'}.
  
  Now assume $n=4$. Let $p^b$ be the order of $\sigma_p$. First assume
  that $\sigma_{p'}$ is orthogonally decomposable. From (a), (b), and
  (c), we get that $\ord{\sigma_{p'}}$ divides $q^2-1$ or $q^3+1$. If
  the latter occurs, then $b=0$. If $b\ge2$, then $b=2$, and either
  $p=3$, and $\sigma_3$ acts indecomposably on $V$, or $p=2$. In the
  former case $\sigma_{3'}$ must be a scalar, so
  $\ord{\overline{\sigma}}$ divides $9$. Next assume that
  $\sigma_{p'}$ acts orthogonally indecomposably on $V$. Then
  $V=Z_1\oplus Z_1'$ with $\dim(Z_1)=2$.  Let $\lambda\in\FFF_{q^4}$
  be an eigenvalue on $Z_1$. Then the other eigenvalue is
  $\lambda^{q^2}$, and Lemma \ref{L:kappa:ind:p'} tells us that the
  eigenvalues on $Z_1'$ are $\lambda^{-q}$ and $\lambda^{-q^3}$. Set
  $m=q^3-q^2+q-1$. Raising these $4$ eigenvalues to the $m$--th power
  gives equal values (use $\lambda^{q^4}=\lambda$), hence
  $\sigma_{p'}^m$ is a scalar. Also, $\sigma_p=1$, because $Z_1$ and
  $Z_1'$ are not $\sigma_{p'}$--isomorphic by Lemma
  \ref{L:kappa:ind:p'}. We get the stated divisibilities.
\end{proof}

\begin{Lemma}\label{L:PSU4:ord}
  Let $q=p^f\ge3$ for a prime $p$. Then each element in
  $\Aut(\PSU_4(q))$ has order at most $\text{max}(2,f)\cdot(q^3+1)$.
\end{Lemma}

\begin{proof}
  Let $\sigma\in\GU_4(q)\rtimes\Gal(\FFF_{q^2}/\FFF_p)$ be a preimage
  of a given element $\overline{\sigma}\in\Aut(\PSU_4(q))$. Let $r$ be
  smallest positive integer with $\sigma^r\in\GU_4(q)$, so $r$ divides
  $2f$. If $r<2f$, then $r\le f$, and
  $\overline{\sigma}^r\in\PGU(4,q)$, so the claim follows from
  $\ord{\overline{\sigma}}\le f\ord{\overline{\sigma}^r}$ and Lemma
  \ref{L:GU:n<5:ord}. Also, if $f=1$, we are obviously done. Therefore
  we are concerned with $r=2f$ with $f\ge2$.
  
  By Lemma \ref{L:Lang:gamma}, we get that $\sigma^{2f}$ is conjugate
  to an element in $\GL_4(p)$, and an upper bound for the element
  orders in the latter group is $p^4$, see Proposition
  \ref{P:L:ord}. Thus $\ord{\sigma}\le 2fp^4$. From $f\ge2$ we obtain
  $2fp^4<f(p^6+1)\le f(q^3+1)$, and we are done.
\end{proof}

\begin{Lemma}\label{L:PSU:ord}
  Let $\Si=\PSU_n(q)$ be the simple unitary group with $n\ge3$, and
  $\sigma\in\Aut(\Si)$. Then
\[
\ord{\sigma} \le
\begin{cases}
  2q^n    & \text{if $q$ is odd,}\\
(3e/2)q^n & \text{in any case.}
\end{cases}
\]
\end{Lemma}

\begin{proof}
  Write $q=p^f$ with $p$ a prime. Then $\sigma$ has a preimage $\tau$
  in $\GU_n(q)\rtimes\Gal(\FFF_{q^2}/\FFF_p)$. Under restricting the
  scalars to $\FFF_p$, we obtain an embedding of the latter group into
  $\Isom(\FFF_{p}^{2fn},\kappa)$, where $\kappa$ is a symmetric
  non--degenerate $\FFF_p$--bilinear form. Now apply the bounds in
  Proposition \ref{P:Bil:ord} to obtain the claim.
\end{proof}

We rule out the unitary groups in the order as they appear in the list
\ref{Tb:ClGr} on page \pageref{Tb:ClGr}. So suppose that
$\Si=\PSU_m(q)$.

$\mathbf{m=3,\;q\ne2,5.}$ First suppose that $f\ge2$, so $q>p$. By
Lemma \ref{L:GU:n<5:ord} and the structure of the automorphism group
of $\PSU_m(q)$ given in \cite[Prop.~2.3.5]{KL:Sub} we get
$\ord{\sigma}\le 2f(q^2-1)$. But $\mu(\Si)=q^3+1$, so
$q^3+1\le2\cdot2f(q^2-1)$, hence $q^2-q+1\le4f(q-1)$. This shows
$q^2-q<4f(q-1)$, so $3^f\le q<4f$, contrary to $f\ge2$.

Next suppose $f=1$, so $q=p$. We obtain $\ord{\sigma}\le 2p(p+1)$.
Thus $p^3+1\le 4p(p+1)$, so $p^2-p+1\le4p$, therefore $p-1<4$, so
$p=3$. Check the atlas \cite{ATLAS} to see that $\ord{\sigma}\le12$,
so this case is out by $3^3+1>2\cdot12$.

$\mathbf{m=3,\;q=5.}$ Then $\Out(\Si)=\sym_3$ and $o(\Aut(\Si))=30$. Thus
the degree is at most $60$. But the only representation of $\Si$ with
degree $\le60$ has degree $50$, see \cite{ATLAS}. Now $o(\Si)=10$, so
$A>\Si$. As $\Si.3$ does not have a permutation representation of degree
$50$, we have $A=\Si.2$. However, $o(\Si.2)=20$, and this case is out too.

$\mathbf{m=4.}$ Suppose $q\ne2$ for the moment. First suppose $f\ge2$.
Then $\ord{\sigma}\le f(q^3+1)$ by Lemma \ref{L:PSU4:ord}. We obtain
$(q+1)(q^3+1)=\mu(\Si)\le 2f(q^3+1)$, hence $q+1\le 2f$. But $q\ge
2^f\ge 2f$, so there is no solution. Next suppose $f=1$, so $q=p$. We
obtain $p+1\le4$, so $p=3$.  However, the maximal element order in
$\Aut(\PSU_4(3))$ is $28$, see the atlas \cite{ATLAS}, a
contradiction. Similarly, if $q=2$, then $o(\Aut(\PSU_4(2)))=12$,
which is too small.

$\mathbf{6\divides m,\;q=2.}$ Use Lemma \ref{L:PSU:ord} to get
$2^{m-1}(2^m-1)/3\le2(3e/2)2^m=6e2^{m-1}$, hence $2^m-1\le
18e=48.9\dots$, so $m\le5$, a contradiction.

$\mathbf{m\ge5,\;(m,q)\ne(6m',2).}$ From Lemma \ref{L:q^m+-1} we
obtain $\mu(\Si)>q^{2m-3}$. On the other hand,
$\ord{\sigma}\le(3e/2)q^m$ by Lemma \ref{L:PSU:ord}, so $q^2\le
q^{m-3}\le 3e=8.1\dots$, thus $q=2$ and $m=5$. (Also $m=6$ would
fulfill the inequality, but this is excluded here.) However, in this
case $\mu(\Si)=165$, whereas $o(\Aut(\Si))=24$, see the atlas
\cite{ATLAS}, a contradiction.

\subsubsection{Projective Special Linear Groups}\label{SSS:PSL}

Now we assume that $\Si=\PSL_n(q)$, and show that except for some
small cases, only the expected elements can act with at most $2$
cycles in the natural representation.

In this section, we use results by Tiep and
Zalesskii \cite[Section 9]{TiepZ:MinChar}
on the three smallest faithful permutation degrees for the simple
groups $\PSL_n(q)$. Unfortunately, their result is mis--stated.
Apparently they mean to give the degrees of the three smallest
faithful \emph{primitive} permutation representations. In order to
make use of their result, we need a little preparation.

\begin{Lemma}\label{L:T:impr:A:pri}
  Let $\Si$ be a simple non--abelian group, and $n=\mu(\Si)$ be the degree
  of the smallest faithful permutation representation. Let $A$ be a
  group between $\Si$ and $\Aut(\Si)$. If $A$ has a primitive permutation
  representation on $\Omega$ such that $\Si$ is imprimitive on $\Omega$,
  then $\abs{\Omega}\ge3n$.
\end{Lemma}

\begin{proof}
  Suppose that $\Si$ acts imprimitively on $\Omega$, and assume that
  $\abs{\Omega}<3n$. Let $\Delta$ be a non--trivial block for $\Si$,
  and $M$ be a setwise stabilizer in $\Si$ of this block. Primitivity of
  $A$ forces transitivity of $\Si$ on $\Omega$, in particular $\Si$ is
  transitive on the block system containing $\Delta$. As there must be
  at least $n$ blocks by assumption,
\[
n\abs{\Delta}\le\abs{\Omega}<3n,
\]
hence $\abs{\Delta}<3$, so $\abs{\Delta}=2$. Let $A_1$ be the
stabilizer of a point in $A$. Set $\Si_1=\Si\cap A_1$, a point--stabilizer
in $\Si$. Clearly $[M:\Si_1]=\abs{\Delta}=2$, so $\Si_1$ is normal in $M$.
Also, $\Si_1$ is normal in $A_1$, and maximality of $A_1$ in $A$ forces
$A_1=N_A(\Si_1)$. So $M\le A_1$, a contradiction.
\end{proof}

\begin{Lemma}\label{L:PSL:impr:A:pri}
  Let $\Si=\PSL_n(q)$ with $(n,q)\ne(4,2)$, $(2,2)$, $(2,3)$, $(2,4)$,
  $(2,5)$, $(2,7)$, $(2,9)$, or $(2,11)$. Let $A$ be a group with
  $\Si\le A\le\Aut(\Si)$. Suppose that $A$ acts primitively, and there is
  $\sigma\in A$ with at most two cycles in this action. Then $\Si$ is
  primitive as well.
\end{Lemma}

\begin{proof}
  In these cases the natural action of $\Si$ on the
  $\mu=(q^n-1)/(q-1)$ lines of $\FFF_q^n$ is the one of smallest
  possible degree. Let $N$ be the degree of the action of $A$. Suppose
  that $\Si$ is imprimitive. From Lemma \ref{L:T:impr:A:pri} we obtain
  $N\ge3\mu$. If $\sigma\in\PgL_n(q)$, then $\ord{\sigma}\le\mu$ by
  Proposition \ref{P:L:ord}, contrary to Lemma \ref{L:kl:Deg-Ord}.
  Thus $\sigma$ involves a graph automorphism of $\PSL_n(q)$, hence
  also $n\ge3$.
  
  As $\sigma^2\in\PgL_n(q)$, we have $\ord{\sigma}\le2\mu$, hence
  $N\le4\mu$. Let $A_1$ be a point--stabilizer in $A$, and set
  $\Si_1=A_1\cap\Si$. Let $M$ be a maximal subgroup of $\Si$ containing
  $\Si_1$. Then $[\Si:M]\le[\Si:\Si_1]/2\le2\mu$, so it follows easily from
  \cite[Section 9]{TiepZ:MinChar} that $M$ fixes a line (or
  hyperplane) with respect to the natural action, except possibly for
  $(n,q)=(3,2)$. Exclude this single exception for a moment. As
  $A=A_1\Si$ by transitivity of $\Si$, also $A_1$ involves a graph
  automorphism $\tau$. As $A_1$ normalizes $\Si_1$, and the action of
  $\tau$ on $\Si$ interchanges point--stabilizers with
  hyperplane--stabilizers, we get that there is a hyperplane
  $H<\FFF_q^n$ and a line $L<\FFF_q^n$, such that $\Si_1$ fixes $H$ and
  $L$. Clearly, $\Si$ acts transitively on the $q^{n-1}(q^n-1)/(q-1)$
  non--incident line--hyperplane pairs, and also transitively on the
  $(q^n-1)(q^{n-1}-1)/(q-1)^2$ incident line--hyperplane pairs. The
  latter size is smaller than the former, so $N\ge
  (q^n-1)(q^{n-1}-1)/(q-1)^2=(q^{n-1}-1)/(q-1)\mu$. From
  $N\le2\ord{\sigma}\le4\mu$ we obtain $1+q+\dots+q^{n-2}\le4$. Hence
  $n=3$ and $q=3$ or $2$. However, for $q=3$ we have
  $\ord{\sigma}\le13$ by \cite{ATLAS}, contrary to $N\ge52$. If $q=2$,
  then $\Aut(\Si)\cong\PGL_2(7)$, so $\ord{\sigma}\le8$, but $N\ge21$, a
  contradiction.
  
  It remains to check the case $(n,q)=(3,2)$. Then
  $\Aut(\Si)\cong\PGL_2(7)$, so $\ord{\sigma}\le8$, hence $N\le16$. But
  this contradicts the above estimation $N\ge3\mu=21$.
\end{proof}

\begin{Lemma}\label{L:PSL:kl}
  Let $\Si=\PSL_n(q)$ with $(n,q)\ne(4,2)$, $(2,2)$, $(2,3)$, $(2,4)$,
  $(2,5)$, $(2,7)$, $(2,9)$, $(2,11)$ and $\Si\le A\le\Aut(\Si)$. Assume
  that $A$ acts primitively on $\Omega$. Suppose that $\sigma\in A$
  has at most $2$ cycles on $\Omega$. Then either $A\le\PgL_n(q)$ and
  $A$ acts naturally on the lines of $\FFF_q^n$, or $(n,q)=(3,2)$, and
  $A\le\Aut(\PSL_3(2))\cong\PGL_2(7)$ acts naturally of degree $8$.
\end{Lemma}

\begin{proof}
  Let $N=\abs{\Omega}$ be the permutation degree of $A$, and suppose
  that we do not have the natural action of $\Si=\PSL_n(q)$ on the
  points of the projective space.
  
  As $\sigma^2\in\PgL_n(q)$, we get $\ord{\sigma}\le2(q^n-1)/(q-1)$ by
  Proposition \ref{P:L:ord}.
  
  $\Si$ is primitive by the previous lemma, so we can use the results by
  Tiep and Zalesskii \cite[Section 9]{TiepZ:MinChar} on the three
  smallest primitive permutation degrees for the simple groups
  $\PSL_n(q)$, see the comment before Lemma \ref{L:T:impr:A:pri}.
  
  First suppose that $n\ge4$, and if $n=4$, then $q\ne2$. Then
\[
N\ge\frac{(q^n-1)(q^{n-1}-1)}{(q^2-1)(q-1)}.
\]
(This second smallest primitive representation is given by the action
on the $2$--spaces in $\FFF_q^n$.)  Now use $N\le2\ord{\sigma}$ to
obtain $q^{n-1}-1\le4(q^2-1)$. So $n=4$ and $q=3$. (Note that
$(n,q)=(4,2)$ is already excluded from the statement of the lemma.)
But $o(\Aut(\PSL_4(3)))=40$ by the atlas \cite{ATLAS}, whereas
$N=130>2\cdot40$, so this case is out.

Next assume $n=3$. Using \cite[Section 9]{TiepZ:MinChar}, one easily
verifies that $N\ge q^3-1$ except for $q=4$ and $2$. Exclude $q=2$ and
$4$ for a moment. So $q^3-1\le 4(q^3-1)/(q-1)$, hence $q=3$ or $5$.
But for $q=5$, we actually have $N\ge5^2(5^3-1)$, but
$\ord{\sigma}\le2(5^3-1)/(5-1)$, clearly a contradiction. If $q=3$,
then $N\ge144$, contrary to $\ord{\sigma}\le 2(3^3-1)/(3-1)=26$. Now
suppose $q=4$. The atlas \cite{ATLAS} gives $\ord{\sigma}\le21$,
whereas $N\ge56>2\cdot21$ by \cite[Section 9]{TiepZ:MinChar}, a
contradiction. If $q=2$, and we do not have the natural action, then
necessarily $N=8$, which corresponds to the natural action of
$\PGL_2(7)\cong\Aut(\PSL_3(2))$.

Finally we have to look at $n=2$. As $A\le\PgL_2(q)$ now, we have
$\ord{\sigma}\le(q^2-1)/(q-1)=q+1$.

We go through the cases in \cite[Section 9]{TiepZ:MinChar}. If $q>4$
is an even square, then $2(q+1)\ge N\ge\sqrt{q}(q+1)$, hence $q\le4$,
a contradiction. If $q$ is an odd square $\ne9$, $49$, then $2(q+1)\ge
N\ge\sqrt{q}(q+1)/2$, hence $q\le16$, a contradiction. If
$q\in\{19,29,31,41\}$, then $2(q+1)\ge N\ge q(q^2-1)/120$, so
$q\le16$, a contradiction. If $q=17$ or $q=49$, then $N=102$ or $175$,
respectively, so these cases do not occur. If $q$ is not among the
cases treated already (and $\ne7$, $9$, and $11$,) then $N\ge
q(q-1)/2$, so $q(q-1)\le4(q+1)$, hence $q\le5$, a contradiction.
\end{proof}

\begin{Lemma}\label{L:PSL2:kl}
  Let $\Si=\PSL_2(q)$ with $q=9$ or $11$ and $\Si\le A\le\Aut(\Si)$. Assume
  that $A$ acts primitively on $\Omega$, and that there exists
  $\sigma\in A$ with at most $2$ cycles on $\Omega$. Then either
  $A\le\PgL_2(q)$ and $A$ acts naturally on the lines of $\FFF_q^n$,
  or $q=9$, $A\le\sym_6<\Aut(\PSL_2(9))$ acting naturally on $6$
  points, or $q=11$, $\abs{\Omega}=11$, $A=\PSL_2(11)$, and $\sigma$
  is an $11$--cycle.
\end{Lemma}

\begin{proof}
  Suppose $q=9$. We have $\Si\cong\alt_6$, and the maximal subgroups of
  $\Si$ have index $6$, $10$, and $15$, respectively. Of course, the
  degree $6$ occurs. Degree $10$ corresponds to the natural action of
  $\Si$. The degree $15$ corresponds to $\alt_6$ acting on $2$--sets.
  Then $A\le\sym_6$, and one verifies easily that each element has
  $\ge3$ cycles. This settles the case that $\Si$ is primitive. If $\Si$
  is imprimitive, then $N\ge3\cdot6=18$ by Lemma \ref{L:T:impr:A:pri},
  but also $N\le2\ord{\sigma}\le20$. As $A$ contains no element of
  order $9$, we actually have $N=20$. Hence $\ord{\sigma}=10$, so
  $\PGL_2(9)\le A$. But neither $\PGL_2(9)$ nor $\PgL_2(9)$ act
  primitively on $20$ points, e.g.\ by the argument in the proof of
  Lemma \ref{L:T:impr:A:pri}.
  
  Next suppose $q=11$. As $\ord{\sigma}\le12$, we have $N\le24$, but
  $24<33=3\cdot\mu(\Si)$, so $\Si$ is primitive. The maximal subgroups of
  $\Si$ of index $\le24$ have index $11$ and $12$, and correspond to the
  actions covered by the claim.
\end{proof}

\begin{Lemma}\label{L:PgL:kl}
  Let $\PgL_n(q)$ act naturally on the lines of $\FFF_q^n$ for
  $n\ge2$. Suppose that an element
  $\sigma\in\PgL_n(q)\setminus\PGL_n(q)$ has at most $2$ cycles. Then
  $(n,q)=(3,4)$, $(2,4)$, $(2,8)$, or $(2,9)$.
\end{Lemma}

\begin{proof}
  Let $\gamma g\in\GL_n(q)\rtimes\Aut(\FFF_q)$ be a preimage of such a
  $\sigma$, with $\gamma\in\Aut(\FFF_q)$ and $g\in\GL_n(q)$. Then
  \[
\ord{\gamma g}\ge\frac{1}{2}\frac{q^n-1}{q-1}.
\]
Let $f\ge2$ be the order $\gamma$. By Lemma \ref{L:Lang:gamma},
$(\gamma g)^f$ is conjugate to an element in $\GL_n(q^{1/f})$, and the
orders of elements in this latter group are at most $q^{n/f}-1$ by
Proposition \ref{P:L:ord}. Thus
\begin{equation}\label{fqn}
f(q^{n/f}-1)\ge \ord{\gamma g}\ge \frac{1}{2}\frac{q^n-1}{q-1}.
\end{equation}
This gives 
\[
2fq>2f(q-1)\ge \frac{q^n-1}{q^{n/f}-1}>q^{n-n/f},
\]
hence 
\[
2f>q^{n-n/f-1}.
\]
Now use $2f\le2^f$ and $q\ge2^f$ to obtain
\[
2^f>2^{nf-n-f},
\]
hence
\begin{equation}
  \label{fn}
  n<\frac{2f}{f-1}=4-2\frac{f-2}{f-1}\le4,
\end{equation}
so $n\le3$.

First suppose $n=3$. Then \eqref{fn} shows $f<3$, hence $f=2$. Set
$r=q^{1/2}$. Then \eqref{fqn} gives
$2(r^3-1)\ge\frac{1}{2}\frac{r^6-1}{r^2-1}$, so $4(r^2-1)\ge r^3+1$,
hence $r<4$. One verifies easily that $r=3$ is not possible, because
the maximal order of an element in $\PgL_3(9)\setminus\PGL_3(9)$ is
$26$, see e.g.~\cite{ATLAS}.

Next assume $n=2$. Again set $r=q^{1/f}\ge2$. Let $h$ be an element in
$\GL_2(r)<\GL_2(q)$ which is conjugate (in $\GL_2(\overline{\FFF_q})$)
to $(\gamma g)^f$. Denote by $\overline{h}$ the image of $h$ in
$\PGL_2(q)$.  There are three possibilities for $h$: If $h$ is
irreducible on $\FFF_q^2$, then $\ord{h}$ divides $r^2-1$, so
$\ord{\overline{h}}$ divides $(r^2-1)/\gcd(r^2-1,q-1)$. But $r-1$
divides the denominator, so $\ord{\overline{h}}$ divides $r+1$. Next
assume that $h$ is reducible.  If $h$ is semisimple, then clearly
$\ord{\overline{h}}\divides\ord{h}\divides r-1$. If however $h$ has a
unipotent part, then this $p$--part has order $p$, and its centralizer
is the group of scalar matrices. Hence in this case,
$\ord{\overline{h}}=p\le r$.

We have seen that $\ord{\overline{h}}\le r+1$ in any case, hence
$\ord{\sigma}\le f(r+1)$. We obtain
\[
f(r+1)\ge\frac{q+1}{2}=\frac{r^f+1}{2},
\]
hence
\[
\frac{r^f+1}{r+1}\le2f.
\]
The left hand side is monotonously increasing in $r$. For $r=2$ we
obtain $2^f+1\le6f$, hence $f\le4$. For $f=3$ and $4$ there are only
the solutions $r=2$. If $r>2$, then $f=2$ and $r=3$ or $4$. In order
to obtain the claim, we have to exclude the possibility $q=r^f=16$.
The previous consideration shows that each element in
$\PgL_2(16)\setminus\PGL_2(16)$ has order at most $12$. But then we
clearly cannot have at most $2$ cycles in a representation of odd
degree $17$.
\end{proof}

\begin{Lemma}\label{L:PGL:kl}
  Let $2\le n\in\NNN$. Suppose that $\sigma\in\PGL_n(q)$ has at most
  $2$ cycles in the action on the lines of $\FFF_q^n$. Then one of the
  following holds:
  \begin{enumerate}
  \item $q$ is a prime, $n=2$, and $\sigma$ has order $q$.
  \item $\sigma$ is a Singer cycle or the square of a Singer cycle.
  \end{enumerate}
\end{Lemma}

\begin{proof}
  For a subset $S$ of $\FFF_q^n$ denote by $P(S)$ the
  ``projectivization'' of $S$, namely the set of $1$--dimensional
  spaces through the non--zero elements of $S$. Denote by
  $\hat{\sigma}\in\GL_n(q)$ a preimage of $\sigma$. If $\hat{\sigma}$
  is irreducible on $\FFF_q^n$, then Schur's Lemma shows that (b)
  holds. Thus assume that $\hat{\sigma}$ is reducible, and let
  $0<U<\FFF_q^n$ be a $\hat{\sigma}$--irreducible subspace.  The
  assumption shows that $\gen{\sigma}$ permutes transitively the
  elements in $P(U)$, as well as those of $P(\FFF_q^n\setminus U)$.
  The transitivity of this latter action shows
  \begin{equation}
    \label{qu}
    q^u\text{ divides $\ord{\hat{\sigma}}$, where $u=\dim(U)$}.
  \end{equation}
  Denote by $\hat{\sigma}_p$ and $\hat{\sigma}_{p'}$ the $p$--part and
  $p'$--part of $\hat{\sigma}$, respectively. Let $W$ be a
  $\hat{\sigma}_{p'}$--invariant complement to $U$ in $\FFF_q^n$. As
  $\hat{\sigma}$ is transitive on $P(U)$ and $P(\FFF_q^n/U)$, we have
  in particular that $\hat{\sigma}$ is irreducible on the quotient
  space $\FFF_q^n/U$, so $\hat{\sigma}_p$ is trivial on this quotient,
  hence $\hat{\sigma}_{p'}$ is irreducible on $W$. From (\ref{qu}) we
  get that $\hat{\sigma}_p$ is not trivial, in particular $W$ is not
  $\hat{\sigma}_p$--invariant. Then we see from Jordan--H\"older that
  $U$ and $W$ are $\hat{\sigma}_{p'}$--isomorphic, so
  $\hat{\sigma}_p\in\GL_2(q^u)$. Thus $\ord{\hat{\sigma}_p}=p$.
  Combine this with \eqref{qu} to get $n=2u=2$, and $q=p$. Finally,
  $\hat{\sigma}_{p'}$ centralizes $\hat{\sigma}_p$, so must be a
  scalar, that is $\sigma$ has order $p$.
\end{proof}

\subsection{Exceptional Groups of Lie Type}\label{SS:ExGr} Here we
rule out the case that $\Si$ is an exceptional group of Lie type.
Table \ref{Tb:ExGr} on page \pageref{Tb:ExGr} contains the exceptional
group of Lie type $\Si$, a lower bound $\mu(\Si)$ for the degree of a
non-trivial transitive faithful permutation representation, an upper
bound $o(\Si)$ for the orders of elements, the order of the outer
automorphism group, and finally restricting condition on $q$. In the
list $q=p^f$ for a prime $p$.

The lower bound for $\mu(\Si)$ has been computed as follows. If $\Si$
has a permutation representation of degree $m$, and $F$ is any field,
then the permutation module of $\Si$ over $F$ has a submodule of
dimension $m-1$. So $m-1$ is at least the dimension of the
lowest--dimensional projective representation of $\Si$ in
characteristic different from the defining characteristic. But these
minimal dimensions have been determined by
Landazuri and Seitz in
\cite{LandSeitz}. We use the corrected list \cite[Theorem
5.3.9]{KL:Sub}. Note that if $\Si$ does not have a doubly transitive
representation, then the $(m-1)$--dimensional module is reducible, so
one summand has dimension at most $(m-1)/2$, see
\cite[4.3.4]{Gorenstein}. This is the case for all $\Si$ except for
$^2B_2(q)$ and $^2G_2(q)$. So $\mu(\Si)$ is then at least $1$ plus $2$
times the minimal dimension of a representation of $\Si$.

The upper bound for $o(\Si)$ has been obtained as follows. Each
element of $\Si$ is the product of a $p$--element with a commuting
$p'$--element, so we multiply upper bounds for each. If $\ell$ is the
Lie rank of $\Si$, then the order of $p'$--elements is at most
$(q+1)^\ell$, see \cite[1.3A]{LS:excmax}. The order of a $p$--element
$g$ is bounded as follows. Suppose $\Si\le\PGL_w(F)$ for a field $F$
of characteristic $p$. Then the order of $g$ is a $p$--power at most
$p(w-1)$, see Lemma \ref{L:unip:ord}. Small values $w$ with an
embedding as above are classically known, see \cite[Prop.\ 
5.4.13]{KL:Sub}. However, for the Suzuki groups
$^2B_2(q)$ we used \cite[XI, \S3]{Huppert3} to determine $\mu$ and
$o$. To determine $\mu$ for $G_2(q)$ and $^3D_4(q)$ we use the papers
by Kleidman \cite{Kleidman:G2} and
\cite{Kleidman:3D4} respectively.

Now assume that $\Si\le A\le\Aut(\Si)$ and $\sigma\in A$ has at most
two cycles in a transitive action of $A$. Then
$\mu(\Si)\le2o(A)\le2\abs{\Out(\Si)}o(\Si)$. Comparing with the
information in the Table \ref{Tb:ExGr} on page \pageref{Tb:ExGr} rules
out all but a few little cases, which require extra data obtained from
the atlas \cite{ATLAS}.

$\mathbf{\Si=^2\!B_2(q)}$. We get $1+q^2\le 2f(q+\sqrt{2q}+1)$. As
$q\ge8$, we have $\sqrt{2q}+1\le\frac{5}{8}q$. So we get $q^2<1+q^2\le
2f(q+\frac{5}{8}q)$, hence $2^f<\frac{13}{4}f$. This implies $f=3$.
But $o(\Aut(^2\!B_2(8)))=15$ (see the atlas \cite{ATLAS}), contrary to
$\mu(^2\!B_2(8))=65>2\cdot15$.

$\mathbf{\Si=^2\!G_2(q)}$. We get $1+q(q-1)\le 2f\cdot9(q+1)$. Now
$q+1\le\frac{28}{27}q$, which gives $3^f=q<\frac{56}{3}f+1$, hence
$f=3$. But $\mu(^2\!G_2(27))=19684$, see \cite{ATLAS}, whereas
$o(^2\!G_2(27))=37$, so this case is clearly out.

$\mathbf{\Si=G_2(q)}$. Obviously $q\ge5$. First assume that $q$ is
odd. Bound $(q^6-1)/(q-1)$ from below by $q^5$, and $q+1$ from above
by $6q/5$. We then obtain $q^5\le2\cdot2f\cdot6p(q+1)^2\le
24q(6q/5)^2$, hence $p^{2f}\le864f/25$, which gives $q=5$. But then
$\Out(\Si)$ has order $1$, and when we use the estimations in the table,
we get a contradiction. The case $p=2$ and $f\ge3$ also does not occur
by a similar calculation.

$\mathbf{\Si=^3\!D_4(q)}$. We get $(q+1)(q^8+q^4+1)/2\le 2\cdot3f\cdot
8p(q+1)^2$. One quickly checks that this holds only for $q=2$. But
$\mu(^3\!D_4(2))=819$, whereas $o(^3\!D_4(2))=28$ (see \cite{ATLAS}),
so this case does not occur.

$\mathbf{\Si=^2\!F_4(2)'}$. This clearly does not occur.

$\mathbf{\Si=^2\!F_4(q)}$. One gets $1+q^4\sqrt{2q}(q-1)\le 2f\cdot
32(q+1)^2$, and one easily checks that this inequality has no
solutions.

$\mathbf{\Si=F_4(q)}$. The case $q=2$ does not occur. We have
$1+2q^6(q^2-1)\le 2(2,p)f\cdot25p(q+1)^4$, which implies that $q=3$ or
$4$. However, Theorem \cite[5.3.9]{KL:Sub} for even $q$ shows that the
minimal degree of a $2'$--representation of $F_4(4)$ is $1548288$, so
$\mu(\Si)\ge3096577$. But this violates the estimation
$o(F_4(4))\le31250$. So $q=3$. The maximal order of a $3'$--element is
$\le73$, see \cite[page 316]{Carter:LNM}. Furthermore, the $3$--order
is at most $27$. Thus $o(\Si)\le1971$. But $\mu(\Si)\ge11665$, a
contradiction.

$\mathbf{\Si=^2\!E_6(q)}$. We get quickly $q=2$. But $o(^2E_6(2))=35$,
which is much too small compared to $\mu(^2E_6(2))=3073$.

$\mathbf{\Si=E_6(q)}$. We quickly get that $q=2$, $3$, or $4$. The
$p'$--part is bounded by $91$, $949$, and $5061$, respectively (again
by \cite[page 316]{Carter:LNM}), and the $p$--part is bounded by $32$,
$27$, and $32$, respectively. So $o(\Si)$ is at most $2912$, $25623$,
and $161952$, respectively. If we compare this with the estimation for
$\mu(\Si)$, then only $q=2$ survives. We get $\mu(\Si)\le
2\cdot2\cdot2912=11648$. However, $E_6(2)$ contains $F_4(2)$, and
$\mu(E_6(2))\ge\mu(F_4(2))=69615$ (\cite{ATLAS}), a contradiction.

$\mathbf{\Si=E_7(q)}$. We get $q=2$. Use \cite[page 316]{Carter:LNM}
to obtain $o(\Si)\le 171\cdot64=10944$. But from the table we have
$\mu(\Si)\ge 196609$, which is clearly too big.

$\mathbf{\Si=E_8(q)}$ gives also no examples.

\subsection{Proof of Part \ref{T:kl:III} of Theorem
  \ref{T:kl}}\label{SS:AS:Proof}

Now we are ready to prove part \ref{T:kl:III} of Theorem \ref{T:kl},
by collecting the information achieved in the last sections. Thus
suppose that $A$ acts primitively, $\Si\le A\le\Aut(\Si)$ for a
non--abelian simple group $\Si$, and that $A$ contains an element
$\sigma$ which has exactly $2$ cycles.

If $\Si$ is sporadic, then Section \ref{SS:Sp} gives the
possibilities. This is the easiest case, as the result can be directly
read off from the atlas information \cite{ATLAS}. Only the Mathieu
groups $\M11$, $\M12$, $\M22$, and $\M24$ give rise to examples.

Section \ref{SS:Alt} treats the case that $\Si=\alt_n$, the
alternating group with $n\ge5$. The case $n=6$ has been excluded
there, and postponed to the analysis of the linear groups, in view of
$\alt_6\cong\PSL_2(9)$. The only examples coming not from the natural
action of $\Si$ are as follows: $\Si=\alt_5$ acting on the $2$--sets
of $\{1,2,3,4,5\}$, hence of degree $10$ (case \ref{T:kl:III:A5(10)}),
or $\Si=\alt_5$ acting on $6$ points (case \ref{T:kl:III:PSL2(p)} for
$p=5$, note that $\alt_5\cong\PSL_2(5)$).

By Section \ref{SS:ExGr}, $\Si$ cannot be of exceptional Lie type.

In Section \ref{SS:ClGr} it is shown that if $\Si$ is a classical group,
then $\Si$ is isomorphic to some $\PSL_n(q)$.

This is dealt with in Section \ref{SSS:PSL}. We can exclude a couple of
small pairs $(n,q)$ in view of exceptional isomorphisms, see
\cite[Prop.~2.9.1]{KL:Sub}. As $\Si$ is simple,
$(n,q)\ne (2,2)$, $(2,3)$. Also, $(n,q)\ne(2,4)$, $(2,5)$, as
$\Si=\alt_5$ has been dealt with already. Also $(n,q)\ne(4,2)$, as
$\alt_8$ had been ruled out in Section \ref{SS:Alt}. Furthermore, we
assume $(n,q)\ne(2,7)$ in view of $\PSL_2(7)\cong\PSL_3(2)$.

Suppose that $q\ne9$, or $11$, if $n=2$. Then $A\le\PgL_n(q)$ acting
naturally on the projective space, or $(n,q)=(3,2)$, and we have the
natural action of $\PSL_2(7)\cong\PSL_3(2)$ of degree $8$, see Lemma
\ref{L:PSL:kl}. Lemma \ref{L:PSL2:kl} shows that for $(n,q)=(2,9)$ the
action is either the natural one, or the natural one of
$\alt_6\cong\PSL_2(9)$, and for $(n,q)=(2,11)$, only the natural
action is possible.

In conclusion, we are left to look at the natural action of
$\PSL_n(q)\le A\le\PgL_n(q)$, and to determine the possibilities for
$\sigma$. By Lemma \ref{L:PgL:kl}, we have actually
$\sigma\in\PGL_n(q)$, except possibly for $(n,q)=(3,4)$, $(2,8)$, or
$(2,9)$. The case $(n,q)=(3,4)$ accounts for \ref{T:kl:III:PsL3(4)} in
Theorem \ref{T:kl:III}. One easily verifies that $\PgL_2(8)$ does not
contain an element with just $2$ cycles (but it does contain
$9$--cycles not contained in $\PGL_2(8)$!). Similarly, if an element
in $\PgL_2(9)\setminus\PGL_2(9)$ has only $2$ cycles, then
$\sigma\in\M10$, and the cycle lengths are $2$ and $10$. This gives
case \ref{T:kl:III:M10} of Theorem \ref{T:kl:III}.

So in addition to the assumption that $A\le\PgL_n(q)$ acts naturally,
we may finally assume $\sigma\in\PGL_n(q)$. Lemma \ref{L:PGL:kl}
finishes this case: Either $q$ is a prime, $n=2$, $\ord{\sigma}=q$ (so
$\sigma$ has cycle lengths $1$ and $q$, case \ref{T:kl:III:PSL2(p)} of
Theorem \ref{T:kl:III}), or $\sigma$ is the square of a Singer cycle
(case \ref{T:kl:III:Singer/2} of Theorem \ref{T:kl}).

By the classification theorem of the finite simple groups, we have
covered all possibilities of $\Si$.\newpage

\subsection{Tables on Minimal Permutation Degrees, Maximal Element
  Orders, etc.}

\begin{table}[h!tp]
\caption{Classical Groups}
\label{Tb:ClGr}
\setlength{\doublerulesep}{1pt}
\newcommand{\epq}{\epsilon(q)}
{\scriptsize
\[
\begin{array}{|c|c|c|c|}
\hline
\Si & \mu(\Si) &  \abs{\Out(\Si)} & m,q\\
\hline
\hline
\PSL_m(q) & (q^m-1)/(q-1) & 2(m,q-1)f,\;m\ge3 &
(m,q)\ne(2,5),\\
       &               & (m,q-1)f,\;m=2    & (2,7),(2,9),\\
       &               &                   & (2,11),(4,2)\\\hline
\PSL_2(7) & 7 & 2 & \\\hline
\PSL_2(9) & 6 & 4 & \\\hline
\PSL_2(11)& 11& 2 & \\\hline
\PSL_4(2)\cong\alt_8 & 8 & 2 & \\\hline
\text{PSp}_{2m}(q) & (q^{2m}-1)/(q-1) & (2,q-1)f,\;m\ge3 &   m\ge2,q\ge3,\\
            &                  & 2f,\;m=2         &  (m,q)\ne(2,3)\\
\hline
\text{Sp}_{2m}(2)  & 2^{m-1}(2^m-1)   & 1        & m\ge3\\
\hline
\Omega_{2m+1}(q) & (q^{2m}-1)/(q-1) & 2f & m\ge3,\\
      &                    &                & q\ge5\text{ odd}\\
\hline
\Omega_{2m+1}(3) & 3^m(3^m-1)/2 & 2 & m\ge3\\
\hline
\text{P}\Omega^+_{2m}(q) & (q^m-1)(q^{m-1}+1)/(q-1) &
2(4,q^m\!-\!1)f,m\ne4 & m\ge4,\;q\ge4\\
 & & 6(4,q^m\!-\!1)f,m=4 &    \\
\hline
\text{P}\Omega^+_{2m}(2) & 2^{m-1}(2^m-1) & 2,m\ne4 & m\ge4\\
                  &                & 6,m=4   &       \\
\hline
\text{P}\Omega^+_{2m}(3) & 3^{m-1}(3^m-1)/2 & 4, m>4\text{ odd} & m\ge4\\
                  &                         & 8, m>4\text{ even} &     \\
                  &                         & 24, m=4          &       \\
\hline
\text{P}\Omega^-_{2m}(q) & (q^m+1)(q^{m-1}-1)/(q-1) & 2(4,q^m+1)f & m\ge4 \\
\hline
\PSU_3(q) & q^3+1 & 2(3,q+1)f       & q\ne2,5 \\
\hline
\PSU_3(5) & 50    & 6               &       \\
\hline
\PSU_4(q) & (q+1)(q^3+1) & 2(4,q+1)f &      \\
\hline
\PSU_m(2) & 2^{m-1}(2^m-1)/3 & 6 &
6\divides m \\
\hline
\PSU_m(q) & \frac{(q^m-(-1)^m)(q^{m-1}-(-1)^{m-1})}{q^2-1} &
2(m,q+1)f & m\ge5, \\
 & & & (m,q)\ne(6m',2)           \\
\hline
\end{array}
\]
}
\end{table}

\newpage

\begin{table}[h!tp]
\caption{Exceptional Groups}
\label{Tb:ExGr}
\[
\begin{array}{|c|c|c|c|c|}
\hline
\text\Si   & \mu(\Si)\ge             & o(\Si)\le    & \abs{\Out(\Si)} & q\\
\hline
^2B_2(q)  & 1+q^2                 & q+\sqrt{2q}+1 & f      & q=2^{2u+1}>2\\
^2G_2(q)  & 1+q(q-1)              & 9(q+1)     & f         & q=3^{2u+1}>3\\
G_2(3)    & 351                   & 13         & 2         &             \\
G_2(4)    & 416                   & 21         & 2         &             \\
G_2(q)    & (q^6-1)/(q-1)         & 8(q+1)^2   & f         &
q\ge8\text{ even} \\
G_2(q)    & (q^6-1)/(q-1)         & 6p(q+1)^2  & \le2f     &
q\ge5\text{ odd} \\
^3D_4(q)  & (q+1)(q^8+q^4+1)/(2,q-1)& 7p(q+1)^2& 3f        &
\\
^2F_4(2)' & 1600                  & 16         & 2         &             \\
^2F_4(q)  & 1+q^4\sqrt{2q}(q-1)   & 32(q+1)^2  & f         & q=2^{2u+1}>2\\
F_4(2)    & 69615                 & 30         & 2         &             \\
F_4(q)    & 1+2q^6(q^2-1)         & 25p(q+1)^4 & (2,p)f    & q\ge3       \\
^2E_6(q)  & 1+2q^9(q^2-1)         & 26p(q+1)^4 & 2(3,q+1)f &             \\
E_6(q)    & 1+2q^9(q^2-1)         & 26p(q+1)^6 & 2(3,q-1)f &             \\
E_7(q)    & 1+2q^{15}(q^2-1)      & 55p(q+1)^7 & (2,q-1)f  &             \\
E_8(q)    & 1+2q^{27}(q^2-1)      & 247p(q+1)^8& f         &             \\
\hline
\end{array}
\]
\end{table}

\newpage

\begin{table}[h!tp]
\caption{Sporadic Groups}
\label{Tb:Sp}
\[
\begin{array}{|l|l|l|c|}
\hline
\text{Group }\Si & \text{Orders of} & \text{Indices of ma-} &
\abs{\Out(\Si)}\\
{ }  & \text{elements} & \text{ximal subgroups} & {}\\
\hline
\M11  &  11,\;8,\;6,\;\le5                  & 11,\;12,\;\ge55    & 1\\
\M12  &  11,\;10,\;8,\;6,\;\le5             & 12,\;\ge66         & 2\\
\M22  &  11,\;8,\;7,\;6,\;\le5              & 22,\;\ge77         & 2\\
\M23  &  23,\;15,\;14,\;\le11               & 23,\;\ge253        & 1\\
\M24  &  23,\;21,\;15,\;14,\;12,\;\le11     & 24,\;\ge276        & 1\\
\text{J}_1 & \le19                          & \ge266             & 1\\
\text{J}_2 & \le15                          & \ge100             & 2\\
\text{J}_3 & \le19                          & \ge6156            & 2\\
\text{J}_4 & \le66                          & \ge173067389       & 1\\
\text{HS}  & \le20                          & \ge100             & 2\\
\text{Suz} & \le24                          & \ge1782            & 2\\
\text{McL} & \le30                          & \ge275             & 2\\
\text{Ru}  & \le29                          & \ge4060            & 1\\
\text{He}  & \le28                          & \ge2058            & 2\\
\text{Ly}  & \le67                          & \ge8835156         & 1\\
\text{O'N} & \le31                          & \ge122760          & 2\\
\text{Co}_1& \le60                          & \ge98280           & 1\\
\text{Co}_2& \le30                          & \ge2300            & 1\\
\text{Co}_3& \le60                          & \ge276             & 1\\
\text{Fi}_{22}& \le30                       & \ge3510            & 2\\
\text{Fi}_{23}& \le60                       & \ge31671           & 1\\
\text{Fi}'_{24}& \le60                      & \ge8672            & 2\\
\text{HN}  & \le40                          & \ge1140000         & 2\\
\text{Th}  & \le39                          & \ge143127000       & 1\\
\text{B}   & \le70                          & \ge4372            & 1\\
\text{M}   & \le119                         & \ge196883          & 1\\
\hline
\end{array}
\]
\end{table}

\newpage

\section{Genus $0$ Systems}
\subsection{Branch Cycle Descriptions}\label{SS:g0}
\subsubsection{Algebraic Setting}
Let $k$ be a subfield of the complex numbers $\CCC$, $t$ be a
transcendental over $\CCC$, and $L/k(t)$ be a finite Galois extension
with groups $G$. We assume that $L/k(t)$ is regular, that means $k$ is
algebraically closed in $L$. Let $\frakp_1,\frakp_2,\dots,\frakp_r$ be
the places of $k(t)$ which are ramified in $L$. Then, by a consequence
of Riemann's Existence Theorem (see \cite{MM}, \cite{H:Buch}), we can
choose places $\frakP_i$ of $L$ lying above $\frakp_i$,
$i=1,2,\dots,r$, and elements $\s i\in G$ such that $\s i$ is a
generator of the inertia group of $\frakP_i$, so that the following
holds:
\begin{quote}
  The $\s i$, $i=1,2,\dots,r$ generate $G$, and $\s1\s2\dots\s r=1$.
\end{quote}
We call the tuple $(\s1,\s2,\dots,\s r)$ a \emph{branch cycle
  description} of the extension
$L/k(t)$.

Now let $E$ be a field between $L$ and $k(t)$, and consider $G$ as a
permutation group on the conjugates of a primitive element of
$E/k(t)$. Set $n:=[E:k(t)]$. For $\sigma\in G$ let
$\ind(\sigma)$ be $n$ minus the number of
cycles of $\sigma$. We call $\ind(\sigma)$ the
\emph{index} of
$\sigma$. This notion obviously applies to any permutation group of
finite degree.

Let $g_E$ be the genus of the field $E$. The Riemann--Hurwitz genus
formula gives
\begin{equation}
  \label{RH}
  2(n-1+g_E)=\sum_{i=1}^r\ind(\sigma_i).
\end{equation}

We will frequently use this relation for the case that $E$ is a
rational field, so that in particular $g_E=0$, and will call the
corresponding equation \emph{genus $0$ relation}, and the tuple $(\s1,\s2,\dots,\s r)$ a \emph{genus $0$
  system}.

The process of constructing a branch cycle description from the
extension $L/k(t)$ can be reverted to some extent. Namely let $G$ be
any finite group, generated by $\s1,\s2,\dots,\s r$, such that
$\s1\s2\dots\s r=1$. Then there exists a finite extension $k/\QQQ$,
and a regular Galois extension $L/k(t)$, such that the $\s i$ arise
exactly as described above. This again follows from (the difficult
direction of) Riemann's Existence Theorem. Modern references are
\cite{MM} and \cite{H:Buch}, where the latter one contains a
self--contained treatment.

\subsubsection{Topological Setting}
For explicit computations and a conceptual understanding of branch
cycle descriptions, the topological interpretation of the $\s i$ is
indispensable. Also $\CCC L/\CCC(t)$ has Galois group $G$. Again let
$E$ be a field between $k(T)$ and $L$. There is a composition of
ramified coverings of Riemann surfaces $\hat{\mathcal
  X}\rightarrow\cX\stackrel{\pi}\rightarrow\PP(\CCC)$, such that the
natural inclusion of the fields of meromorphic functions
$\CCC(t)=M(\PP(\CCC))\subseteq M(\cX)\subseteq M(\hat{\mathcal X})$ is
just the extension $\CCC(t)\subseteq \CCC E\subseteq \CCC L$. If we
identify the places of $\CCC(t)$ with the elements in $\PP(\CCC)$ in
the natural way, then the branch points of $\hat{\mathcal
  X}\to\PP(\CCC)$ are exactly the places of $\CCC(t)$ ramified in
$\CCC L$. Choose a point $\frakp_0\in\PP(\CCC)$ away from the branch
points $\frakp_i$, and choose a standard set of generators
$\gamma_1,\gamma_2,\dots,\gamma_r$ of the fundamental group $\Gamma$
of $\PP(\CCC)\setminus\{\frakp_1,\dots,\frakp_r\}$ with base point
$\frakp_0$, where $\gamma_i$ comes from a path starting and ending in
$\frakp_0$, winding clockwise around $\frakp_i$ just once and not
around any other branch point, see the diagram.

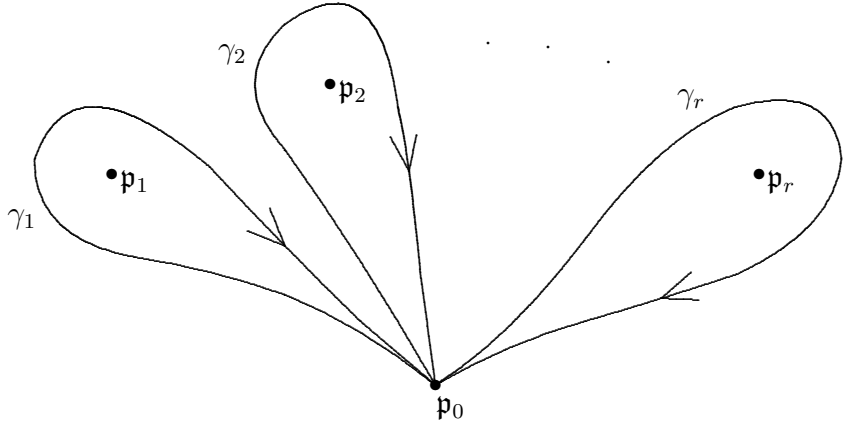
\setlength{\unitlength}{1mm}
\begin{picture}(130,70)

\curve(65,10,55,17,45,22)
\curve(45,22,35,25,25,27)
\curve(25,27,15,31.5,11.7,40)
\curve(11.7,40,20,47,35,39)
\curve(35,39,45,28.5,55,18.5)
\curve(55,18.5,65,10)

\curve(43,33.5,45,28.5)
\curve(45,28.5,40,30.5)

\curve(65,10,55,26.5,45,41.3)
\curve(45,41.3,41,50,45,58)
\curve(45,58,55,59.7,60,47)
\curve(60,47,61.5,38.5,63,25)
\curve(63,25,65,10)

\curve(59,43,61.5,38.5)
\curve(61.5,38.5,62.5,43.3)

\curve(65,10,75,18,85,29.5)
\curve(85,29.5,95,40.7,105,47)
\curve(105,47,115,47,119,40)
\curve(119,40,115,31.5,105,24.7)
\curve(105,24.7,95,21.5,85,18.5)
\curve(85,18.5,75,15,65,10)

\curve(99,25,95,21.5)
\curve(95,21.5,100,21.3)

\put(65,10){\circle*{1.5}}
\put(67,7){\makebox(0,0){$\frakp_0$}}

\put(22,38){\circle*{1.5}}
\put(25,37){\makebox(0,0){$\frakp_1$}}
\put(51,50){\circle*{1.5}}
\put(54,49){\makebox(0,0){$\frakp_2$}}
\put(108,38){\circle*{1.5}}
\put(111,37){\makebox(0,0){$\frakp_r$}}

\put(72,55.5){\circle*{.5}}
\put(80,55){\circle*{.5}}
\put(88,53){\circle*{.5}}

\put(10,32){\makebox(0,0){$\gamma_1$}}
\put(38,54){\makebox(0,0){$\gamma_2$}}
\put(99,48){\makebox(0,0){$\gamma_r$}}
\end{picture}

The $\gamma_i$ generate $\Gamma$ with the single relation
$\gamma_1\gamma_2\dots\gamma_r=1$. Clearly $\Gamma$ acts on the fiber
$\pi^{-1}(\frakp_0)$. The induced action gives the group $G$, and the
images of the $\gamma_i$ are the elements $\sigma_i$ as above.
Furthermore, the cycle lengths of $\sigma_i$ on the fiber
$\pi^{-1}(\frakp_0)$ are the multiplicities of the elements in the
fiber $\pi^{-1}(\frakp_i)$, and these cycle lengths are the same as
for the corresponding action on the conjugates of a primitive element
of $E/k(t)$.

For more details about this connection we refer again to \cite{MM} and
\cite{H:Buch}.

\subsection{Branch Cycle Descriptions in Permutation
  Groups}\label{SS:bcd}

Let $G$ be a transitive permutation group of degree $n$, and
$\cE:=(\s1,\s2,\dots,\s r)$ be a generating system with $\s1\s2\dots\s
r=1$. For $\sigma\in G$ define the index $\ind(\sigma)$ as
above. Let the number $g_\cE$ be given by
\[
2(n-1+g_\cE)=\sum_{i=1}^r\ind(\s i).
\]
The topological interpretation from above of the $\sigma_i$ as coming
from a suitable cover of Riemann surfaces shows that $g_\cE$ is a
non--negative integer, because it is the genus of a Riemann
surface. This topological application in a purely group theoretic
context was first made by Ree, see \cite{Ree}. Later, 
Feit, Lyndon, and
Scott gave an elementary group theoretic
argument of this observation, see \cite{FLS}.

In this chapter we will determine such systems $\cE$ for $g_\cE=0$ in
specific groups $G$. According to the previous section, we will call
such systems genus $0$ systems. If we look for $\s i$ in a fixed
conjugacy class $\cC_i$, then it does not matter in which way we order
the classes, for if $\s i$ and $\sigma_{i+1}$ are two consecutive
elements in $\cE$, then we may replace these elements by
$\sigma_{i+1}$ and $\s i^{\sigma_{i+1}}$, respectively.

The strategy of finding such genus $0$ systems in $G$ (or proving that
there are none) depends very much on the specific situation. For many
small groups, we simply check using a program written in GAP
\cite{GAP}. For bigger groups, especially certain sporadic groups, we
can use the character tables in the atlas \cite{ATLAS}. Here, and at
other places, the following easy observation (see \cite[2.4]{PM:Mon})
is useful.

\begin{Lemma}\label{L:Indchi}
  Let $\sigma\in G$, where $G$ is a permutation group of degree $n$,
  then
\[
\ind(\sigma)=n-\frac{1}{\ord{\sigma}}\sum_{k|\ord{\sigma}}
\chi(\sigma^k)\varphi(\frac{\ord{\sigma}}{k}),
\]
where $\chi(\tau)$ is the number of fixed points of $\tau\in G$, and
$\varphi$ is the Euler $\varphi$--function.
\end{Lemma}

\subsection{A Lemma about Genus $0$ Systems}

\begin{Lemma}\label{L:g0:gcd}
  Let $(\sigma_1,\sigma_2,\dots,\sigma_r)$ be a genus $0$ system of a
  transitive permutation group $G$. Suppose that all cycle lengths of
  $\sigma_1$ and $\sigma_2$ are divisible by $d>1$. Then $G$ admits a
  block system of $d$ blocks, which are permuted cyclically.
\end{Lemma}

\begin{proof}
  Let $n$ be the degree of $G$. Let $\cX\to\PP(\CCC)$ be a connected
  cover of the Riemann sphere, such that
  $(\sigma_1,\sigma_2,\dots,\sigma_r)$ is the associated branch cycle
  description. Without loss of generality let $0$ and $\infty$ be
  branch points corresponding to $\sigma_1$ and $\sigma_2$,
  respectively. As our tuple is a genus $0$ system, $\cX$ has genus
  $0$, thus $\cX=\PP(\CCC)$ and the cover is given by a rational
  function $f(X)$. We may assume (by a linear fractional change) that
  $\infty$ is not mapped to $0$ or $\infty$.  Let $\alpha_i$ be the
  elements in $f^{-1}(0)$, and denote the multiplicity of $\alpha_i$
  my $m_i$.  Similarly, let $\beta_i$ have multiplicity $n_i$ in the
  fiber $f^{-1}(\infty)$.  Thus, up to a constant factor, we have
\[
f(X)=\frac{\prod(X-\alpha_i)^{m_i}}{\prod(X-\beta_i)^{n_i}}.
\]
As the $m_i$ and $n_i$ are the cycle lengths of $\sigma_1$ and
$\sigma_2$, respectively, we get $f(X)=g(X)^d$, where $g(X)\in\CCC(X)$
is a rational function. From that the claim follows.
\end{proof}

\begin{Remark*} The completely elementary nature of the lemma makes it
  desirable to have a proof which does not rely on Riemann's existence
  theorem. We sketch an elementary argument, and leave it to the
  reader to fill in the details: First note that if the claimed
  assertion about the permutation action holds for a group containing
  $G$ (and acting on the same set), then it holds for $G$ as well. For
  $i>2$ write $\sigma_i$ as a minimal product of transpositions, and
  replace the element $\sigma_i$ by the tuple of these transpositions.
  This preserves the genus $0$ condition. Also, the product of a
  $k$--cycle with a disjoint $l$--cycle with a transposition which
  switches a point of the $k$--cycle with one of the $l$--cycle is a
  $(k+l)$--cycle. This way, we can assume that all cycle lengths of
  $\sigma_1$ and $\sigma_2$ are $d$, at the cost of extra
  transpositions, but still preserving the genus $0$ property. Write
  $n=md$. Clearly, there are $m-1$ transpositions in our system, such
  that they, together with $\sigma_1$, generate a transitive group.
  Let $\tau_1,\dots,\tau_{m-1}$ be these transpositions. As we have a
  genus $0$ system, the total number of transpositions is $2(m-1)$.
  Using braiding we get an equation of the form
\[
\sigma_1\tau_1\dots\tau_{m-1}=\sigma_2'\tau_1'\dots\tau_{m-1}'=:\rho,
\]
where $\sigma_2'$ is conjugate to $\sigma_2^{-1}$, and the $\tau_i'$
are transpositions. As $\ind(xy)\le\ind(x)+\ind(y)$ and
$(\sigma_1,\tau_1,\dots,\tau_{m-1},\rho^{-1})$ is a genus $\ge0$
system of a transitive subgroup of $G$, we obtain it must be a genus
$0$ system, and $\ind(\rho)=n-1$. Thus $\rho$ is an $n$--cycle.
Inductively, we see that $\lambda:=\sigma_1\tau_1\dots\tau_{m-2}$ is a
product of an $(n-d)$--cycle and a $d$--cycle, and that these two
cycles are fused by $\tau_{m-1}$. Now, by induction on the degree of
$G$, we get that the group generated by the transitive genus $0$
system $(\sigma_1,\tau_1,\dots,\tau_{m-2},\lambda^{-1})$ with respect
to the support of size $n-d$ admits a block system of $d$ blocks being
permuted cyclically. Now extend each block $\Delta$ by a single point
from the remaining $d$ points as follows: Choose $j$ such that
$\tau_{m-1}$ moves a point $\omega$ from $\Delta^{\sigma_1^j}$. Now
append $\omega^{\tau_{m-1}\sigma_1^{-j}}$ to $\Delta$.  One verifies
that this process is well--defined, and gives a block system for
$(\sigma_1,\tau_1,\dots,\tau_{m-1})$ with $d$ blocks being permuted
cyclically. It remains to show that this block system is preserved
also by $(\sigma_2',\tau_1',\dots,\tau_{m-1}')$. At any rate, by
symmetry we get a block system for this tuple too, with $d$ blocks
being permuted cyclically. The point is that the product of the
elements in this tuple is the same $n$--cycle as the product of the
elements in the former tuple, and an $n$--cycle has a unique block
system with $d$ blocks. Therefore the block systems are the same, so
are respected by $G$.
\end{Remark*}

\subsection{The Siegel-Lang Theorem and Hilbert's Irreducibility
  Theorem}
Let $k$ be field which is finitely generated over $\QQQ$, and $R$ a
finitely generated subring. The Siegel-Lang Theorem about points with
coordinates in $R$ on algebraic curves over $k$ has the following
application to Hilbert's irreducibility theorem, see
\cite[2.1]{PM:hitzt}:

Let $f(t,X)\in k(t)[X]$ be irreducible, and $\Red_f(R)$ the set of
those $\bar t\in R$, such that $f(\bar t,X)$ is defined, and reducible 
over $k$. Then, up to a finitely many elements, $\Red_f(R)$ is the
union of finitely many sets of the form $g(k)\cap R$, where $g(Z)\in
k(Z)$  is a rational function.

In view of this result, it is important to know which rational
functions $g(Z)$ have the property that $g(k)\cap R$ is an infinite
set. By another theorem of Siegel-Lang (see \cite[8.5.1]{Lang:Dio}),
this property implies that there are at most two elements of $\bar
k\cup\{\infty\}$ in the fiber $g^{-1}(\infty)$.

The converse is true if we allow to enlarge $R$. More precisely, we
have the following

\begin{Lemma} Let $k$ be a finitely generated extension of $k$,
  $g(Z)\in k(Z)$ a non-constant rational function such that the fiber
  $g^{-1}(\infty)$ has at most two elements. Then there is a finitely
  generated subring $R$ of $k$ with $\abs{g(k)\cap R}=\infty$.
\end{Lemma}

\begin{proof} A linear fractional
  change of the argument of $g$ allows to assume that $g(Z)$ has the
  following shape: There is $m\ge0$ and a polynomial $A(Z)\in k[Z]$,
  such that either $g(Z)=A(Z)/Z^m$, or $g(Z)=A(Z)/(Z^2-d)^m$, where
  $d\in k$ is not a square. In the first case let $R$ be the ring
  generated by $1/2$ and the coefficients of $A(Z)$, then $g(z)\in R$
  if $z=2^r$ for $r\in\ZZZ$, hence $g(k)\cap R$ is an infinite set.
  
  The second case is a little more subtle: For $\alpha,\beta\in k$,
  define sequences $\alpha_n,\beta_n\in k$ for $n\in\NNN$ by
  $\alpha_n+\beta_n\sqrt{d}=(\alpha+\beta\sqrt{d})^n$. Suppose for the
  moment that $\beta_n\ne0$ for all $n$. Define
  $z_n=\alpha_n/\beta_n$.  Then the $z_n$ are pairwise distinct, for
  if $z_j=z_i$ for $j>i$, then
  $\beta_j/\beta_i=(\alpha+\beta\sqrt{d})^{j-i}=
  \alpha_{j-i}+\beta_{j-i}\sqrt{d}$, so $\beta_{j-i}=0$, a
  contradiction. Let $R$ be the ring generated by the coefficients of
  $A(Z)$ and $1/(\alpha^2-\beta^2d)$. Note that
  $(\alpha^2-\beta^2d)^n=\alpha_n^2-\beta_n^2d$, so $g(z_n)\in R$ for
  all $n$.
  
  It remains to show that we can choose suitable $\alpha,\beta\in k$.
  Write $\gamma=\alpha+\beta\sqrt{d}$. Then $\beta_n=0$ is equivalent
  to $\gamma^n\in k$. Thus we have to find $\gamma$ such
  $\gamma^n\not\in k$ for all $n\in\NNN$. Suppose that $\gamma\not\in
  k$, however $\gamma^n\in k$. Let $\sigma$ be the involutory
  automorphism of $k(\sqrt{d})/k$. The minimal polynomial
  $(X-\gamma)(X-\sigma(\gamma))$ of $\gamma$ over $k$ divides
  $X^n-\gamma^n$, so $\sigma(\gamma)=\zeta\gamma$ for $\zeta$ an $n$th
  root of unity. In particular, $\zeta\in k(\sqrt{d})$. But
  $k(\sqrt{d})$ is finitely generated, so this field contains only
  finitely many roots of unity. For a fixed $\gamma\in
  k(\sqrt{d})\setminus k$ consider the elements $\gamma+i$, for
  $i\in\ZZZ$. For each $i$ there is $n_i\in\NNN$ with
  $(\gamma+i)^{n_i}\in k$. By the above, each element
  $(\sigma(\gamma)+i)/(\gamma+i)$ is a root of unity. One of these
  roots of unity appears for infinitely many $i$, which of course is
  nonsense.
\end{proof}

If $k=\QQQ$, then one is mainly interested in the special case
$R=\ZZZ$. Then $\abs{g(\QQQ)\cap\ZZZ}=\infty$ has another strong
consequence, see \cite{Siegel}: If $\abs{g^{-1}(\infty)}=2$, then the
two elements in $g^{-1}(\infty)$ are real and algebraically conjugate.

Motivated by these results, we make the following

\begin{Definition}
  Let $k$ be a field which is finitely generated over $\QQQ$, and
  $g(Z)\in k(Z)$ a non-constant rational function. We say that $g(Z)$
  is a \emph{Siegel function} over $k$, if there is a finitely
  generated subring $R$ of $k$ with $\abs{g(k)\cap R}=\infty$. If
  $k=\QQQ$, then we require more strongly that
  $\abs{g(\QQQ)\cap\ZZZ}=\infty$.
\end{Definition}

In the analysis of Siegel functions we will make use of the fact about
$g^{-1}(\infty)$. For this it will not be necessary to assume that $k$
is finitely generated. Thus we define a more general property, which
holds for Siegel functions.

\begin{Definition} Let $k$ be a field of characteristic $0$, and
  $g(Z)\in k(Z)$ a non-constant rational function. We say that $g(Z)$
  fulfills the \emph{Siegel property}, if $\abs{g^{-1}(\infty)}\le2$.
  
  If $k=\QQQ$, and $\abs{g^{-1}(\infty)}=2$, then we additionally
  require the two elements in $g^{-1}(\infty)$ to be real and
  algebraically conjugate.
\end{Definition}

\subsection{Siegel Functions and Ramification at Infinity}

To ease the language, we start to define monodromy groups of rational
functions.
\begin{Definition} Let $k$ be a field of characteristic $0$, and
  $g(Z)\in k(Z)$ be a non-constant rational function. Denote by $L$ a
  splitting field of $g(Z)-t$ over $k(t)$. Set $A=\Gal(L/k(t))$,
  considered as a permutation group on the roots of $g(Z)-t$. Denote
  by $\hk$ the algebraic closure of $k$ in $L$, and let
  $G\trianglelefteq A$ be the normal subgroup $\Gal(L/\hk(t))$. Then
  $A$ and $G$ are called the \emph{arithmetic monodromy group} and
  \emph{geometric monodromy group} of $g(Z)$, respectively. Note that
  $A/G$ is naturally isomorphic to $\Gal(\hk/k)$.
\end{Definition}

Our goal is to determine the genus $0$ systems and the monodromy
groups of functionally indecomposable rational functions with the
Siegel property. L\"uroth's Theorem shows that functional
indecomposability of $g(Z)$ (over $k$) implies primitivity of $A$.
The following lemma summarizes the properties we will use.

\begin{Lemma}\label{L:Siegel:A=GD}
  With the notation from above let $D\le A$ and $I\trianglelefteq D$
  be the decomposition and inertia group of a place of $L$ lying above
  the place $t\mapsto\infty$ of $k(t)$, respectively. Suppose that
  $g(Z)$ has the Siegel property. Then the following holds.
\begin{itemize}
\item[(a)] The cyclic group $I$ has at most two orbits, with lengths
  equal the multiplicities of the elements in $g^{-1}(\infty)$.
\item[(b)] $A=GD$ and $I\le G\cap D$. In particular, $A=N_A(I)G$.
\item[(c)] If $A$ is primitive, then $G$ is primitive, too.
\item[(d)] $G$ has a genus $0$ system, with a generator of $I$
  belonging to it.
\end{itemize}
\end{Lemma}

\begin{proof} For (a) and (b) see \cite[Lemma 3.4]{PM:hitzt}, for (c)
  see \cite[Theorem 3.5]{PM:hitzt}, and (d) follows from Section
  \ref{SS:g0}.
\end{proof}

\subsection{Monodromy Groups and Ramification of Siegel Functions}

The main result of this section  is

\begin{Theorem}\label{T:g0} Let $g(Z)$ be a non-constant, functionally
  indecomposable rational function over a field $k$ of characteristic
  $0$. Suppose that $\abs{g^{-1}(\infty)}=2$. Let $A$ and $G$ be the
  arithmetic and geometric monodromy group of $g(Z)$, respectively,
  and $(\sigma_1,\sigma_2,\dots,\sigma_r)$ a branch cycle description.
  Let $T$ be the unordered tuple
  $(\ord{\sigma_1},\ord{\sigma_2},\dots,\ord{\sigma_r})$. Then either
  $\alt_n\le G\le A\le\sym_n$, with many possibilities for $T$, or one
  of the following holds, where $G\le A\le A_{\text{max}}$:
\begin{enumerate}
\item\label{T:g0:I} $A$ acts as an affine group, and one of the
  following holds:
\par
\begin{tabular}{|c|c|c|p{57mm}|}
\hline
$n$ & $G$ & $A_{\text{max}}$ & $ $\hfill  $T$ \hfill $ $\\
\hline
$5$      & $\AGL_1(5)$   & $G$          & $(2,4,4)$         \\
$7$      & $\AGL_1(7)$   & $G$          & $(2,3,6)$         \\
$8$      & $\AgL_1(8)$   & $G$          & $(3,3,6)$, $(3,3,7)$ \\
$8$      & $\AGL_3(2)$   & $G$          & many cases \\
$9$      & $\AgL_1(9)$   & $G$  &  $(2,4,8)$ \\
$9$      & $\AGL_2(3)$   & $G$  &  $(2,3,8)$, $(2,6,8)$, $(2,2,2,8)$ \\
$16$     & $C_2^4\rtimes(C_5\rtimes C_4)$ & $G$ & $(2,4,8)$ \\
$16$     & \text{index} $2$ in $A_{\text{max}}$ &
        $(\sym_4\times\sym_4)\rtimes C_2$ & $(2,4,8)$ \\
$16$     & $C_2^4\rtimes\sym_5$ & $G$ & many cases \\
$16$     & $\AgL_2(4)$    & $G$ & $(2,4,15)$ \\
$16$     & $C_2^4\rtimes\alt_7$  & $G$ & $(2,4,14)$ \\
$16$     & $\AGL_4(2)$  & $G$ & many cases \\
$32$     & $\AGL_5(2)$  & $G$ & probably many cases \\
$64$     & $\AGL_6(2)$  & $G$ & probably many cases \\
\hline
\end{tabular}
\item\label{T:g0:II} (Product action) $n=m^2$,
  $G=A=(\sym_m\times\sym_m)\rtimes C_2$, many possibilities for $T$.
\item\label{T:g0:III} $A$ is almost simple, and one of the following
  holds:
\par
\begin{tabular}{|c|c|c|p{57mm}|}
\hline
$n$ & $G$ & $A_{\text{max}}$ & $ $\hfill  $T$ \hfill $ $\\
\hline
$6$      & $\PSL_2(5)$ & $\PGL_2(5)$ & many cases \\
$6$      & $\PGL_2(5)$ & $G$ & $(2,4,5)$, $(4,4,5)$, $(4,4,3)$ \\
$8$      & $\PSL_2(7)$ & $\PGL_2(7)$ & $(2,3,7)$, $(3,3,7)$, $(3,3,4)$ \\
$8$      & $\PGL_2(7)$ & $G$ & $(2,6,7)$, $(2,6,4)$ \\
$10$     & $\alt_5$ & $\sym_5$ & $(2,3,5)$ \\
$10$     & $\sym_5$ & $G$ & $(2,4,5)$, $(2,6,5)$, $(2,2,2,5)$ \\
$10$     & $\PSL_2(9)$ & $\PgL_2(9)$ & $(2,4,5)$ \\
$10$     & $\PsL_2(9)$ & $\PgL_2(9)$ & $(2,6,5)$, $(2,2,2,5)$ \\
$10$     & $\M10$ & $\PgL_2(9)$ & $(2,4,8)$ \\
$10$     & $\PgL_2(9)$ & $G$ & $(2,8,8)$ \\
$12$     & $\M11$ & $G$ & many cases \\
$12$     & $\M12$ & $G$ & many cases \\
$14$     & $\PSL_2(13)$ & $\PGL_2(13)$ & $(2,3,7)$, $(2,3,13)$ \\
$21$     & $\PsL_3(4)$ & $G$ & $(2,4,14)$ \\
$21$     & $\PgL_3(4)$ & $G$ & $(2,3,14)$, $(2,6,14)$, $(2,2,2,14)$ \\
$22$     & $\M22$ & $\M22\rtimes C_2$ & $(2,4,11)$ \\
$22$     & $\M22\rtimes C_2$ & $G$ & $(2,4,11)$, $(2,6,11)$, $(2,2,2,11)$ \\
$24$     & $\M24$ & $G$ & many cases \\
$40$     & $\PSL_4(3)$ & $\PGL_4(3)$ & $(2,3,20)$ \\
$40$     & $\PGL_4(3)$ & $G$ & $(2,4,20)$ \\
\hline
\end{tabular}
\end{enumerate}
\end{Theorem}

For completeness, we state the analogous result if
$\abs{g^{-1}(\infty)}=1$. The proof follows immediately from Lemma
\ref{L:Siegel:A=GD} and the classification result in \cite{PM:Mon}.

\begin{Theorem}\label{T:g0pol} Let $g(Z)$ be a non-constant,
  functionally indecomposable rational function over a field $k$ of
  characteristic $0$ with $\abs{g^{-1}(\infty)}=1$. Let $G$ be the
  geometric monodromy group of $g(Z)$, and
  $(\sigma_1,\sigma_2,\dots,\sigma_r)$ be a branch cycle description.
  Let $T$ be the unordered tuple
  $(\ord{\sigma_1},\ord{\sigma_2},\dots,\ord{\sigma_r})$. Then one of
  the following holds holds:
\begin{enumerate}
\item Infinite series:
\begin{enumerate}
\item $n=p$, $G=C_p$, $T=(p,p)$, $p$ a prime.
\item $n=p$, $G=D_p$, $T=(2,2,p)$, $p$ an odd prime.
\item $G=\alt_n$ ($n$ odd) or $\sym_n$, many possibilities for $T$.
\end{enumerate}
\item Sporadic cases:
\begin{enumerate}
\item $n=6$, $G=\PGL_2(5)$, $T=(2,4,6)$.
\item $n=7$, $G=\PGL_3(2)$, $T=(2,3,7)$, $(2,4,7)$, or $(2,2,2,7)$.
\item $n=8$, $G=\PGL_2(7)$, $T=(2,3,8)$.
\item $n=9$, $G=\PgL_2(8)$, $T=(2,3,9)$ or $(3,3,9)$.
\item $n=10$, $G=\PgL_2(9)$, $T=(2,4,10)$.
\item $n=11$, $G=\PSL_2(11)$, $T=(2,3,11)$.
\item $n=11$, $G=\M11$, $T=(2,4,11)$.
\item $n=13$, $G=\PGL_3(3)$, $T=(2,3,13)$, $(2,4,13)$,
$(2,6,13)$, or $(2,2,2,13)$.
\item $n=15$, $G=\PGL_4(2)$, $T=(2,4,15)$, $(2,6,15)$, or
$(2,2,2,15)$.
\item $n=21$, $G=\PgL_3(4)$, $T=(2,4,21)$.
\item $n=23$, $G=\M23$, $T=(2,4,23)$.
\item $n=31$, $G=\PGL_5(2)$, $T=(2,4,31)$.
\end{enumerate}
\end{enumerate}
\end{Theorem}

\subsection{Proof of Theorem \ref{T:g0}}

The strategy is as follows. Functional indecomposability of $g(Z)$
implies that $A$ is primitive. By Lemma \ref{L:Siegel:A=GD}(a) we can
apply Theorem \ref{T:kl}. It remains to find normal subgroups $G$ of
$A$ for which (b) and (d) of Lemma \ref{L:Siegel:A=GD} hold. For that
it is useful to know that $G$ is primitive as well by Lemma
\ref{L:Siegel:A=GD}(c).

The proof is split up into three sections, according to whether $A$
acts as an affine group, preserves a product structure, or is almost
simple.

\subsubsection{Affine Action}

The proof is based on work by Guralnick, Neubauer, and Thompson on
genus $0$ systems in primitive permutation groups of affine type.

Suppose that $A$ is affine. The cases that $A$ has degree $\le4$ are
immediate, so assume $n\ge5$.

$G$ is primitive by Lemma \ref{L:Siegel:A=GD}. Let $\sigma_r$ be a
generator of $I$, so $\sigma_r$ has two cycles.

Let $N$ be the minimal normal subgroup of $A$. First suppose that
$G''=1$. As $G'$ is abelian, we have $G'=N$, and primitivity of $G$
forces that $G/N$ acts irreducibly on $N$. But $G/N=G/G'$ is abelian,
so $G/N$ is cyclic by Schur's Lemma. More precisely, we can identify
$G$ as a subgroup of $\AGL_1(q)$, where $q=\abs{N}=p^m$ for a prime
$p$. As $q>4$, we have necessarily that $\sigma$ fixes a point and
moves the remaining ones in a $(q-1)$--cycle. An element in $N$ has
index $q(1-1/p)\ge q/2$, whereas an element in $\AGL_1(q)$ of order
$t|q-1$ has index $(q-1)(1-1/t)\ge(q-1)/2$. The index relation gives
$r=3$ and that neither $\sigma_1$ nor $\sigma_2$ is contained in $N$.
So $2(q-1)=q-2+(q-1)(1-1/t_1+1-1/t_2)\ge q-2+(q-1)(1/2+2/3)$, where
$t_i$ is the order of $\sigma_i$. It follows $q\le7$.

Next suppose that $G''>1$. Write $n=p^m$. We use \cite[Theorems 4.1,
5.1]{GT}. If $p>5$, then $p=7$ or $11$, and $m=2$. Furthermore
$T=(2,4,6)$ for $p=7$, or $T=(2,3,8)$ for $p=11$.  So this does not
occur in view of $\ord{\sigma_r}\ge n/2=p^2/2$. Next suppose $p=5$. We
use \cite[Theorem 1.5]{Neubauer:CA1} (the statement is already in
\cite{GT}, but only parts are proven there). Again compare
$\ord{\sigma_r}\ge n/2$ with the possible genus $0$ systems given for
$p=5$. Only $n=25$ with $G=(C_5\times C_5)\rtimes(\SL_2(5)\rtimes
C_2)$ could arise. However, this group does not have an element with
only two cycles by Theorem \ref{T:kl}.

So we have $p=3$ or $2$. Suppose that $p=3$. Use \cite[Theorem
1.5]{Neubauer:CA1} to see that necessarily $n=9$. Check directly that
only the listed degree $9$ cases occur.

Now suppose $p=2$. By \cite{GurNeu:Affine}, we automatically get
$n\le64$. The cases for $n\le16$ are small enough to be checked with
GAP \cite{GAP}. If $n=32$ then $G=A=\AGL_5(2)$ or $n=64$ and
$G=A=\AGL_6(2)$, $\GL_3(4)\le G\le A\le\gL_3(4)$, or $\GL_2(8)\le G\le
A\le\gL_2(8)$ by Theorem \ref{T:kl}. The last two cases are quickly
ruled out with GAP. We did not check for genus $0$ systems in the
former two cases, but we believe that they exist.

\subsubsection{Product Action}

Let $A$ be a non-affine group which preserves a product structure.
Again let $\sigma_r$ be the element with two cycles.

By Theorem \ref{T:kl}, we have have $A=(U\times U)\rtimes C_2$ in
product action, where either $U=\sym_m$, or $U=\PGL_2(p)$ for a prime
$p\ge5$. By primitivity of $G$ we cannot have $G\le (U\times U)$. On
the other hand, the presence of $\sigma_r$ forces $U\times U\le G$,
see the proof of Theorem \ref{T:kl}, so $G=A$.
  
Let $\Delta$ be the set $U$ is acting on, and let
$\Omega:=\Delta\times\Delta$ be the set $G=A$ acts on.
  
We show the existence of a genus $0$ system of the required form for
$U=\sym_m$. Write $\Delta:=\{1,2,\dots,m\}$. Let $\tau\in G$ be the
element which maps $(i,j)$ to $(j,i)$. Let $1\le a<m$ be prime to $m$.
For $\alpha:=(1,2,\dots,m)\in\sym_m$ and
$\beta:=(a,a-1,\dots,2,1)(m,m-1,\dots,a+2,a+1)\in\sym_m$ set
$\sigma_1:=(\alpha,\beta)\in A$, $\sigma_2:=\tau$,
$\sigma_3:=(\sigma_1\sigma_2)^{-1}$. We show that
$(\sigma_1,\sigma_2,\sigma_3)$ is a genus $0$ system of $G$.
  
First we show that $\sigma_1$ and $\sigma_2$ generate $G$. Note that
$a,m-a$, and $m$ are pairwise prime. Let $r$ and $s$ be integers such
that $rm\equiv1\pmod{a(m-a)}$ and $sa(m-a)\equiv1\pmod{m}$.  Then
clearly $\sigma_1^{rm}=(1,\beta)$ and $\sigma_1^{sa(m-a)}=(\alpha,1)$.
Conjugating with $\tau$ shows that also $(\beta,1),(1,\alpha)\in G$.
We are done once we know that $\alpha,\beta$ generate $\sym_m$. But
this is clear, because it is easy to see that the generated group is
doubly transitive and contains the transposition $\alpha\beta=(a,m)$.
  
We compute the index of $\sigma_i$. The element $\sigma_1$ has a cycle
of length $ma$, and another one of length $m(m-a)$. So
$\ind(\sigma_1)=m^2-2$. Furthermore, $\ind(\sigma_2)=(m^2-m)/2$,
because $\sigma_2=\tau$ has exactly $m$ fixed points, and switches the
remaining points in cycles of length $2$. Next,
$\sigma_3:=\tau(\alpha^{-1},\beta^{-1})$. The element $(i,j)\in\Omega$
is a fixed point of $\sigma_3$ if and only if $j=i^\alpha$ and
$i=j^\beta$, hence $j=i+1$ with $i\ne a,m$. Thus there are exactly
$m-2$ fixed points. Now $\sigma_3^2=((a,m),(a+1,1))$ has order $2$ and
exactly $(m-2)^2$ fixed points. Lemma \ref{L:Indchi} gives
\begin{align*}
\ind(\sigma_3) &=
m^2-\frac{1}{4}(\varphi(4)(m-2)+\varphi(2)(m-2)^2+\varphi(1)m^2)\\
&= (m^2+m)/2,
\end{align*}
so the genus of $(\sigma_1,\sigma_2,\sigma_3)$ is $0$.

We now show that $U=\PGL_2(p)$ does not occur. Again, let $\tau$ be
the element which flips the entries of $\Omega$. At least two of the
elements in $\sigma_1,\dots,\sigma_{r-1}$ must be of the form
$\sigma=(\alpha,\beta)\tau$, with $\alpha,\beta\in\PGL_2(p)$. This
$\sigma$ is conjugate in $G$ to $(1,\alpha\beta)\tau$.  If
$\alpha\beta=1$, then $\ind(\sigma)=((p+1)^2-(p+1))/2$. Otherwise,
$\ind(\sigma)\ge2((p+1)^2-4)/3$, because
$\sigma^2\sim(\alpha\beta,\alpha\beta)$ has at most $4$ fixed points.

If $\sigma$ has the form $(\alpha,\beta)$, then $\sigma$ has at most
$4(p+1)$ fixed points, so $\ind(\sigma)\ge((p+1)^2-4(p+1))/2$.

As $\sum_{i=1}^{r-1}\ind(\sigma_i)=(p+1)^2$, it follows from these
index bounds that $r=3$, so $\sigma_1$ and $\sigma_2$ have the
$\tau$--part. Because not both $\sigma_1$ and $\sigma_2$ can be
involutions (for $G$ is not dihedral), we obtain
$(p+1)^2\ge((p+1)^2-(p+1))/2+2((p+1)^2-4)/3$, so $p<5$, a
contradiction.

\subsubsection{Almost Simple Action}

Let $\Si$ be the simple non-anabelian group with $\Si\le G\le
A\le\Aut(\Si)$, and $\sigma_r$ again the element with two cycles. We
have to check the groups in Theorem \ref{T:kl}(III) for the existence
of genus $0$ systems of the required form.
  
If $\Si=\alt_n$ ($n$ even) in natural action, then it is easy to check
that there are many such genus $0$ systems, and it is obviously not
possible to give a reasonable classification of them.
  
Next, the cases except the infinite series \ref{T:kl:III:PSL2(p)} and
\ref{T:kl:III:Singer/2} of Theorem \ref{T:kl} are easily dealt with,
using the atlas \cite{ATLAS} and some ad hoc arguments, or more
conveniently using \cite{GAP}.
  
Now assume $\PSL_2(q)\le G\le\PgL_2(q)$ in the natural action, with
$q\ge5$ a prime power. Note that $q$ is odd. As $n=q+1$ and
$\ind(\sigma_r)=n-2$, the index relation gives
\[
q+1=\sum_{k=1}^{r-1}\ind(\sigma_k).
\]
We distinguish two cases.

First assume $G\le\PGL_2(q)$. For $\sigma\in\PGL_2(q)$ we easily
obtain (see e.~g.\ \cite{PM:Mon}) that
$\ind(\sigma)\ge(q-1)(1-1/\ord{\sigma})$. So the index relation gives
\[
\sum_{k=1}^{r-1}(1-1/\ord{\sigma_k})\le\frac{q+1}{q-1}.
\]
As $G$ is not dihedral, either $r\ge4$, or $r=3$ and $\sigma_1$ and
$\sigma_2$ are not both involutions. In the first case, we obtain
$q=5$, and in the second case,
$\sum_{k=1}^2(1-1/\ord{\sigma_k})\ge(1-1/2)+(1-1/3)$ gives $q\le13$.
Check these cases directly.

Next suppose that $G\not\le\PGL_2(q)$, but $G\le\PgL_2(q)$. Check the
case $q=9$, $\sigma_r\not\in\PGL_2(9)$ directly and exclude it in the
following. Thus $\sigma_r\in\PGL_2(q)$ by Lemma \ref{L:PgL:kl}.
Denote by $\bar\sigma_k$ the image of $\sigma_k$ in the abelian group
$\PgL_2(q)/\PGL_2(q)$.  Then the elements $\bar\sigma_k$ for
$k=1,\dots,r-1$ are not all trivial and have product $1$. Thus the
order of two of the elements $\sigma_1,\sigma_2,\dots,\sigma_{r-1}$
have a common divisor $\ge2$. Furthermore, for $\sigma\in\PgL_2(q)$,
we have the index bound
$\ind(\sigma)\ge(1-1/\ord{\sigma})(q-\sqrt{q})$, see \cite{PM:Mon}.
This information, combined with the index relation, gives
\begin{align*}
\frac{5}{4} &=    (1-\frac{1}{2})+(1-\frac{1}{4})\\
            &\le  \sum_{k=1}^{r-1}(1-1/\ord{\sigma_k})\\
            &\le  \frac{1}{q-\sqrt{q}}\sum_{k=1}^{r-1}\ind(\sigma_k)\\
            &=    \frac{q+1}{q-\sqrt{q}}.
\end{align*}
Hence $q\le5\sqrt{q}+4$, so $q=9$, $25$, or $27$. If $q=27=3^3$, then
the above argument shows that the common divisor can be chosen to be
$3$, so the analogous calculation gives $4/3\le(27+1)/(27-\sqrt{27})$,
which does not hold. Similarly, refine the argument (using
\cite{PM:Mon}) or simply check with GAP \cite{GAP} that $q=25$ does
not occur.

The main case which is left to investigate is case III(d) of Theorem
\ref{T:kl}, namely that $\PSL_m(q)\le G\le\PgL_m(q)$ acts naturally on
the projective space, $q$ is an odd prime power, $m\ge2$ is even, and
$\sigma_r$ is the square of a Singer cycle. The case $m=2$ has been
done above. The case $m\ge4$, which is somewhat involved, will be
handled in the remaining part of this section. In order to finish the
almost simple case, we need to show that $m=4$, $q=3$, giving the
degree $n=40$ cases in Theorem \ref{T:g0}.

For this we need the following index bounds.

\begin{Lemma}\label{L:Indbound}
  Let $q$ be a prime power, and $1\ne\sigma\in\PgL_m(q)$, where
  $m\ge4$. Then the following holds:
\begin{enumerate}
\item\label{L:Indbound_a}
  $\ind(\sigma)\ge(1-1/\ord{\sigma})(q^{m-1}-1)$.
\item\label{L:Indbound_b} If $\ord{\sigma}$ is a prime not dividing
  $q(q-1)$, and $\sigma\in\PGL_m(q)$, then
  $\ind(\sigma)\ge(1-1/\ord{\sigma})q^{m-2}(q+1)$.
\item\label{L:Indbound_c} If $\ord{\sigma}$ is a prime dividing $q$,
  and $\sigma\in\PGL_m(q)$, then
  $\ind(\sigma)=(1-1/\ord{\sigma})(q^m-q^j)/(q-1)$ for some $1\le j\le
  m-1$.
\end{enumerate}
\end{Lemma}

\begin{proof}
  For \ref{L:Indbound_a} see \cite{PM:Mon}.

Set $N:=(q^m-1)/(q-1)$, and let $s$ be the order of $\sigma$.

Now assume the hypothesis in \ref{L:Indbound_b}. Let $\chi(\sigma)$ be
the number of fixed points of $\sigma$. Then clearly
$\ind(\sigma)=(N-\chi(\sigma))(1-1/s)$. Let $\hat\sigma\in\GL_m(q)$ be
a preimage of $\sigma$ of order $s$. For $\alpha\in\FFF_q$, let
$d(\alpha)$ be the dimension of the eigenspace of $\hat\sigma$ with
eigenvalue $\alpha$. Clearly
\[
\chi(\sigma)=\sum_{\alpha\in\FFF_q}\frac{q^{d(\alpha)}-1}{q-1}.
\]
So $\chi(\sigma)\le(q^d-1)/(q-1)$, where $d=\sum_\alpha d(\alpha)$. On
the other hand, as $s$ does not divide $q-1$, $\hat\sigma$ must have
eigenvalues not in $\FFF_q$. So $d\le m-2$, and the claim follows.

To prove \ref{L:Indbound_c}, note that a preimage of order $s$ of
$\sigma$ in $\GL_m(q)$ admits Jordan normal form over $\FFF_q$.
\end{proof}

Recall that $N=(q^m-1)/(q-1)$. Note that $\ind(\sigma_r)=N-2$, so the
index relation gives
\begin{equation}
\label{qm_ind}
\sum_{k=1}^{r-1}\ind(\sigma_k)=N.
\end{equation}

\begin{Claim}
  $r=3$.
\end{Claim}

\begin{proof}
  Suppose that $r\ge4$. From \ref{L:Indbound_a} in Lemma
  \ref{L:Indbound} we have
  $\ind(\sigma_k)\ge(1-1/\ord{\sigma_k})(q^{m-1}-1)$, hence

\begin{align}
  \sum_{k=1}^{r-1}(1-1/\ord{\sigma_k}) &\le
  \frac{N}{q^{m-1}-1}\notag\\
  &=
  1+\frac{1}{q-1}+\frac{1}{q^{m-1}-1}\label{1/gk_2}\\
  &\le 1+\frac{1}{q-1}+\frac{1}{q^3-1}\notag\\
  &< 1+\frac{2}{q-1}\notag.
\end{align}

First note that if $r\ge4$, then $3/2<1+\frac{2}{q-1}$, so $q<5$ and
hence $q=3$. We get more precisely $\sum_{k=1}^{r-1}(1-1/\ord{\sigma_k})\le
1+1/(3-1)+1/(27-1)=20/13$. However, $2(1-1/2)+(1-1/3)=5/3>20/13$, so
besides $q=3$ we obtain $r=4$, and $\sigma_1$, $\sigma_2$, $\sigma_3$ are
involutions. Note that $\sigma_4$ has cycles of even length, as $4|N$. So
these involutions do have fixed points by Lemma \ref{L:g0:gcd}. Let
$\hat\sigma$ be a preimage in $\GL_m(3)$ of an involution in $\PGL_m(3)$
with fixed points. Thus ${\hat\sigma}^2$ has eigenvalue $1$ on the one
hand, but is also scalar. So $\hat\sigma$ has only the eigenvalues $1$ and
$-1$, and both eigenvalues occur. This shows $\chi(\sigma)\equiv2\pmod{3}$,
hence $\ind(\sigma)\equiv (N-2)/2\equiv1\pmod{3}$. So
\[
1\equiv N=\sum_{k=1}^3\ind(\sigma_k)\equiv0\pmod{3},
\]
a contradiction.
\end{proof}

\begin{Claim}
  $q\le7$.
\end{Claim}

\begin{proof}
  From \eqref{1/gk_2} and $r=3$ we obtain
\begin{equation}
\label{g1g2}
\frac{1}{\ord{\sigma_1}}+\frac{1}{\ord{\sigma_2}}\ge
1-\frac{1}{q-1}-\frac{1}{q^{m-1}-1} \ge
1-\frac{1}{q-1}-\frac{1}{q^3-1}.
\end{equation}
$\sigma_1$ and $\sigma_2$ are not both involutions (because $G$ is not
dihedral). This gives $1/2+1/3\ge1-1/(q-1)-1/(q^3-1)$, so $q<8$.
\end{proof}

In the following we assume $\ord{\sigma_1}\le\ord{\sigma_2}$.

\begin{Claim} $q\ne7$.\end{Claim}

\begin{proof}
  Suppose $q=7$. From \eqref{1/gk_2} we obtain
  $1/\ord{\sigma_1}+1/\ord{\sigma_2}\ge1-1/6-1/(7^3-1)>3/4$, hence
  $\ord{\sigma_1}=2$, $\ord{\sigma_2}=3$. Again, as
  $2|(N/2)=\ord{\sigma_3}$, we get that $\sigma_1$ has fixed points,
  and so $\chi(\sigma_1)\equiv2\pmod{7}$, hence
  $\ind(\sigma_1)\equiv3\pmod{7}$. From
  $3+2(N-\chi(\sigma_2))/3=3+\ind(\sigma_2)\equiv N\equiv1\pmod{7}$ it
  follows that $\chi(\sigma_2)\equiv4\pmod{7}$. So a preimage
  $\hat\sigma_2\in\GL_m(7)$ of $\sigma_2$ has exactly $4$ different
  eigenvalues $\lambda$ in $\FFF_7$. Let $\hat\sigma_2^3$ be the
  scalar $\rho$. The equation $X^3-\rho$ has at most $3$ roots in
  $\FFF_7$, a contradiction.
\end{proof}

\begin{Claim} $q\ne5$. \end{Claim}

\begin{proof}
  Suppose $q=5$. The proof is similar to the argument in the previous
  claim, so we only describe the steps which differ from there. We
  obtain $\ord{\sigma_1}=2$ and $\ord{\sigma_2}=3$ or $4$.
  
  First assume that $\ord{\sigma_2}=3$. As $3|N$, we obtain that
  $\sigma_2$ has fixed points by Lemma \ref{L:g0:gcd}, so a preimage
  $\hat{\sigma_2}\in\GL_m(5)$ has eigenvalues in $\FFF_5$. Suppose
  (without loss, as $\gcd(q-1,3)=1$) that $1$ is one of the
  eigenvalues. As $(X^3-1)/(X-1)$ is irreducible in $\FFF_5$, this is
  the only $\FFF_5$--eigenvalue of $\hat{\sigma_2}$. So
  $\chi(\sigma_2)\equiv1\pmod{5}$, hence
  $\ind(\sigma_2)\equiv0\pmod{5}$. This gives
  $\chi(\sigma_1)\equiv{4}\pmod{5}$, which is clearly not possible.
  
  Now assume that $\ord{\sigma_2}=4$. The index relation together with
  Lemma \ref{L:Indchi} gives
\begin{equation}
\label{2^2}
2\chi(\sigma_1)+2\chi(\sigma_2)+\chi(\sigma_2^2)=N.
\end{equation}
Clearly, $\chi(\sigma_2^2)\ge\chi(\sigma_2)$. If $\chi(\sigma_2^2)=0$,
then $\chi(\sigma_1)\equiv3\pmod{5}$, which is not possible.  Thus
$\sigma_2^2$ has fixed points.

First assume that $\sigma_2$ has no fixed points. Then $\sigma_1$ has
fixed points by Lemma \ref{L:g0:gcd}, so
$\chi(\sigma_1)\equiv2\pmod{5}$. From that we obtain
\[
2((5^a-1)+(5^{m-a}-1))+((5^b-1)+(5^{m-b}-1))=5^m-1
\]
for suitable $1\le a,b\le m-1$. However, $5^a+5^{m-a}\le 5+5^{m-1}$,
and similarly for $b$, so $3(5+5^{m-1})\ge 5(5^{m-1}+1)$. This gives
$5^{m-1}\le5$, a contradiction.

So $\sigma_2$ has fixed points as well, therefore all eigenvalues of a
preimage $\hat{\sigma_2}\in\GL_m(5)$ are in $\FFF_5$. Without loss
assume that $1$ is an eigenvalue of $\hat{\sigma_2}$, and denote by
$a,b,c,d$ the multiplicity of the the eigenvalue $1,2,3,4\in\FFF_5$,
respectively. Clearly $b+c>0$, as $\hat{\sigma_2}$ has order $4$.
Also, $a>0$ by our choice. We obtain that
$\chi(\sigma_2^2)=(5^{a+d}-1)/4+(5^{b+c}-1)/4$, hence
$\chi(\sigma_2^2)\equiv2\pmod{5}$. Relation \eqref{2^2} gives
$\chi(\sigma_1)+\chi(\sigma_2)\equiv2\pmod{5}$. If $\sigma_1$ has
fixed points, then $\chi(\sigma_1)\equiv2\pmod{5}$, hence
$\chi(\sigma_2)\equiv0\pmod{5}$, which is not the case. Thus
$\chi(\sigma_1)=0$ and $\chi(\sigma_2)\equiv2\pmod{5}$, so $d=0$ and
either $b=0$ or $c=0$. Suppose without loss $c=0$. Hence
$\chi(\sigma_2)=\chi(\sigma_2^2)$, and we obtain
\[
N=\frac{5^m-1}{4}=2\chi(\sigma_2)+\chi(\sigma_2^2)=3\chi(\sigma_2)=
  3\left(\frac{5^a-1}{4}+\frac{5^{m-a}-1}{4}\right),
\]
so
\[
5^m+5=3(5^a+5^{m-a})\le3(5+5^{m-1}),
\]
a contradiction as previously.
\end{proof}

\begin{Claim}
  If $q=3$, then $m=4$ and $(\ord{\sigma_1},\ord{\sigma_2})=(2,3)$ or
  $(2,4)$.
\end{Claim}

\begin{proof}
  As $\ind(\sigma_k)\ge (3^{m-1}-1)/2$, and $\ind(\sigma_k)\ge
  2(3^{m-2}-1)$ unless $\sigma_k$ is an involution in $\PGL_m(3)$ of
  minimal possible index, we obtain from the index relation
  \eqref{qm_ind} that
\begin{equation}
\label{q=3}
\ind(\sigma_k)\le
\begin{cases}
3^{m-1}\text{ in any case},\\
\frac{5\cdot3^{m-2}+3}{2}\text{ for $k=2$ if $\sigma_1$ has not minimal
possible index.}
\end{cases}
\end{equation}

We first note that no prime $s\ge5$ does divide $\ord{\sigma_k}$, for
\eqref{q=3} and Lemma \ref{L:Indbound}\ref{L:Indbound_b} would give
\[
(1-\frac{1}{5})3^{m-2}4\le\ind(\sigma_k)\le 3^{m-1},
\]
which is nonsense.

Similarly, we see that $9$ does not divide $\ord{\sigma_k}$. Let
$\sigma\in\PGL_m(3)$ have order $9$, and let $\hat\sigma\in\GL_m(3)$
be a preimage of order $9$. So $\hat\sigma$ admits Jordan normal form
over $\FFF_3$, and there must be at least one Jordan block of size
$\ge4$ by Lemma \ref{L:unip:ord}. Thus $\chi(\sigma)\le(3^{m-3}-1)/2$,
and also $\chi(\sigma^3)\le(3^{m-1}-1)/2$. Now
\[
\ind(\sigma)=(1-\frac{1}{9})N-\frac{2}{3}\chi(\sigma)-
\frac{2}{9}\chi(\sigma^3)
\]
by Lemma \ref{L:Indchi}. Use the above estimation to obtain after some
calculation that $\ind(\sigma)\ge 32\cdot3^{m-4}>3^{m-1}$, contrary to
\eqref{q=3}.

Now suppose that $4$ divides the order of $\sigma_k$. Let $\sigma$ be
a power of $\sigma_k$ of order $4$. As $\sigma_k$ must have a cycle of
odd length by Lemma \ref{L:g0:gcd}, $\sigma$ must have a fixed point.
Thus there is a preimage $\hat\sigma\in\GL_m(3)$ of $\sigma$ with
$\hat\sigma^4=1$. Let $a$ and $b$ be the number of Jordan blocks of
size $1$ with eigenvalue $1$ and $-1$, respectively, and let $j$ be
the number of square blocks of size $2$.  The square of such a block
matrix is a scalar with eigenvalue $-1$. We have $a+b+2j=m$, and $2\le
a+b\le m-2$. Also, $\chi(\sigma)=(3^a-1+3^b-1)/2$ and
$\chi(\sigma^2)=(3^{a+b}-1+3^{2j}-1)/2$. From that we obtain
\begin{align*}
  \ind(\sigma) &= \frac{3}{4}N-\frac{1}{2}\chi(\sigma)-
  \frac{1}{4}\chi(\sigma^2)\\
  &=
  \frac{3}{4}N-\frac{3^a+3^b-2}{4}-\frac{3^{a+b}+3^{m-a-b}-2}{8}\\
  &\ge\frac{3}{4}N-\frac{3^{m-2}-1}{4}-\frac{3^{m-2}+7}{8}\\
  &= 3^{m-1}-1.
\end{align*}
Note that $\ind(\sigma_k)\ge\ind(\sigma)$. From that we see that
$k=2$, and by \eqref{q=3} it follows that $\sigma_1$ is an involution
with minimal possible index. Thus $\ind(\sigma_2)=3^{m-1}$ again by
\eqref{q=3}. This shows that $\ord{\sigma_2}$ is not divisible by $3$,
because then a cycle of $\sigma_2$ of length divisible by $3$ would
break up into at least $3$ cycles of $\sigma$, so
$\ind(\sigma_2)\ge2+\ind(\sigma)\ge1+3^{m-1}$, a contradiction to
\eqref{q=3}.

Similarly, we see that $8$ does not divide $\ord{\sigma_2}$. Suppose
otherwise. Then we get the same contradiction unless
$\ord{\sigma_2}=8$ and $\sigma_2$ has exactly $1$ cycle of length $8$.
But then $\sigma_2^4$ has $N-8$ fixed points, however
$\chi(\sigma_2^4)\le(3^{m-1}+1)/2$, so $(3^m-1)/2-8\le (3^{m-1}+1)/2$,
so $3^{m-1}\le9$, a contradiction.

So $\ord{\sigma_2}=4$, and $\ind(\sigma_2)=3^{m-1}$ by what we have
seen so far. Express $\ind(\sigma_2)$ in terms of $a$ and $b$ as
above. As $\hat\sigma_1$ fixes a hyperplane pointwise, and
$\gen{\hat\sigma_1,\hat\sigma_2}$ is irreducible, we infer that
$a,b\le1$.  Also, $a+b>0$, so $a=b=1$ because $a+b$ is even.
Substitute $a=b=1$ in the relation $\ind(\sigma_2)=3^{m-1}$ to get
$3^{m-1}=27$, so $m=4$.  This case indeed occurs.

Next we look at elements of order $6$. Let $\sigma\in\PGL_m(3)$ have
order $6$, and $\hat\sigma\in\GL_m(3)$ be a preimage. We have
\[
\ind(\sigma)=\frac{5}{6}N-\frac{1}{3}\chi(\sigma)-
\frac{1}{3}\chi(\sigma^2)-\frac{1}{6}\chi(\sigma^3).
\]
Clearly
\[
\chi(\sigma^2)\le\frac{3^{m-1}-1}{2}
\]
and
\[
\chi(\sigma^3)\le\frac{3^{m-1}+1}{2}.
\]
If $\sigma$ has no fixed points, then $\hat\sigma^6=-\bn1$, and
therefore $\sigma^3$ has no fixed points as well. In this case, we
thus obtain $\ind(\sigma)\ge5N/6-\chi(\sigma^2)/3\ge
(13\cdot3^{m-1}-3)/12>3^{m-1}$. This, in conjunction with \eqref{q=3},
shows that if $\ord{\sigma_k}=6$, then $\sigma_k$ has a fixed point.
Suppose that $\sigma=\sigma_k$ has order $6$ and a fixed point. Then
$\hat\sigma$ admits Jordan normal form over $\FFF_3$, and one realizes
easily that
\[
\chi(\sigma)\le\frac{3^{m-2}-1+3^1-1}{2}=\frac{3^{m-2}+1}{2}.
\]
Using this, one obtains after some calculation
\[
\ind(\sigma)\ge\frac{17\cdot3^{m-1}-9}{18}.
\]
However, $(17\cdot3^{m-1}-9)/18>(5\cdot3^{m-2}+3)/2$, so we get from
\eqref{q=3} that $k=2$ and $\sigma_1$ is an involution with minimal
index.  So $\ind(\sigma_2)=3^{m-1}$ by \eqref{qm_ind}, and $\sigma_1$
leaves a hyperplane invariant. The irreducibility of
$\gen{\sigma_1,\sigma_2}$ forces that $\hat\sigma_2$ has eigenspaces
of dimension at most $1$. On the other hand, the Jordan blocks of
$\hat\sigma_2$ have size at most $3$. As $m\ge4$, there is thus
exactly one Jordan block with eigenvalue $1$, and exactly one with
eigenvalue $-1$. Let $u$ and $m-u$ be the size of these blocks,
respectively.  Clearly $\chi(\sigma_2)=2$, $\chi(\sigma_2^2)=4$, and
$\chi(\sigma_2^3)=(3^u+3^{m-u}-2)/2$. From that one computes
\[
\ind(\sigma_2)=\frac{5\cdot3^m-3^u-3^{m-u}-27}{12}.
\]
Now $\ind(\sigma_2)=3^{m-1}$ yields the equation $3^m=3^u+3^{m-u}+27$,
which gives $3^{m-u}=(3^u+27)/(3^u-1)$. Check that the right hand side
is never a power of $3$ for $u=1,2,3$.

It remains to look at $\ord{\sigma_2}=3$. Then $\ord{\sigma_1}=2$ or
$3$. Note that $\ind(\sigma_2)=3^{m-1}-3^{j_2-1}$ by Lemma
\ref{L:Indbound}\ref{L:Indbound_c}, where $j_2$ is the number of
Jordan blocks. Suppose that $\ord{\sigma_1}=3$, and let $j_1$ be the
number of Jordan blocks. The index relation yields
$3^{m-1}+1=2(3^{j_1-1}+3^{j_2-1})$. Looking modulo $3$ shows that
$j_1=j_2=1$. But this gives $m=2$, a contradiction.

Finally, suppose $\ord{\sigma_1}=2$. As the cycles of $\sigma_3$ are
divisible by $2$, Lemma \ref{L:g0:gcd} shows that $\sigma_1$ has fixed
points. Then $\ind(\sigma_1)=(3^m-3^i-3^{m-i}+1)/4$, where $1\le i\le
m-1$ is the multiplicity of the eigenvalue $1$ of an involutory
preimage of $\sigma_1$ in $\GL_m(3)$. The index relation yields
\[
3^{m-1}=3^i+3^{m-i}+4\cdot3^{j_2-1}-3.
\]
If $i=1$ or $m-1$, then the right hand side is bigger than the left
hand side. Thus $2\le i\le m-2$. Looking modulo $9$ then shows that
$j_2=2$, so we get $3^{m-1}=3^i+3^{m-i}+9$. Looking modulo $27$ reveals
that $3^i=3^{m-i}=9$, thus $m=4$. This occurs indeed.
\end{proof}

\section{Siegel Functions over the Rationals}\label{S:Rat}

\subsection{Monodromy Groups and Ramification}

The main arithmetic constraint on monodromy groups is given in the
following lemma, see \cite[Lemma 3.4]{PM:hitzt}:

\begin{Lemma}\label{L:Siegel:rat}
  Let $g(Z)\in\QQQ(Z)$ be a rational function of degree $n=2m\ge2$,
  such that $g^{-1}(\infty)$ consists of two real elements, which are
  algebraically conjugate in $\QQQ(\sqrt{d})$, for $d>1$ a square-free
  integer. Let $t$ be a transcendental over $\QQQ$, and $L$ a
  splitting field of $g(Z)-t$ over $\QQQ(t)$.
  
  Let $D\le A$ and $I\trianglelefteq D$ be the decomposition and
  inertia group of a place of $L$ lying above the place
  $t\mapsto\infty$ of $\QQQ(t)$, respectively.
  
  Then $I=\gen{\sigma}$ for some $\sigma\in G$, and the following
  holds.
\begin{enumerate}
\item $\sigma$ is a product of two $m$--cycles.
\item $\sigma^k$ is conjugate in $D$ to $\sigma$ for all $k$ prime to
  $m$.
\item $D$ contains an element which switches the two orbits of $I$.
\item $D$ contains an involution $\tau$, such that
  $\sigma^\tau=\sigma^{-1}$, and $\tau$ fixes the orbits of $I$
  setwise.
\item If $\sqrt{d}\not\in\QQQ(\zeta_m)$ (with $\zeta_m$ a primitive
  $m$--th root of unity), then the centralizer $C_D(I)$ contains an
  element which interchanges the two orbits of $I$.
\end{enumerate}
\end{Lemma}

The main result about the monodromy groups of Siegel functions over
$\QQQ$ is.

\begin{Theorem}\label{T:g0_Q}
  Let $g(Z)\in\QQQ(Z)$ be a functionally indecomposable rational
  function of degree $n\ge2$ with $\abs{g^{-1}(\infty)}=2$.  Let $A$
  and $G$ be the arithmetic and geometric monodromy group of $g$,
  respectively.  Let $T$ be the ramification type of $g$. Then one of
  the following holds:
\begin{enumerate}
\item $n$ is even, $\alt_n\le G\le A\le\sym_n$, many possibilities for
  $T$; or
\item $n=6$, $G=\PSL_2(5)$, $A=\PGL_2(5)$, $T=(2,5,3)$ and
  $(2,2,2,3)$; or
\item $n=6$, $G=\PGL_2(5)=A$, $T=(4,4,3)$; or
\item $n=8$, $G=\AGL_3(2)=A$, $T=(2,2,3,4)$, $(2,2,4,4)$, and
  $(2,2,2,2,4)$; or
\item $n=10$, $\Si\le G\le A\le\Aut(\Si)$, where $\Si=\alt_5$ or
  $\alt_6$, with many possibilities for $T$; or
\item $n=16$, $G=(S_4\times S_4)\rtimes C_2=A$, $T=(2,6,8)$,
  $(2,2,2,8)$; or
\item $n=16$, $G=C_2^4\rtimes S_5=A$, $T=(2,5,8)$, $(2,6,8)$, and
  $(2,2,2,8)$.
\end{enumerate}
\end{Theorem}

The analogue of the previous theorem for Siegel functions with
$\abs{g^{-1}(\infty)}=1$ follows from the classification of the
monodromy groups of polynomials. For completeness, we give the result
from \cite{PM:Mon}.

\begin{Theorem}\label{T:g0pol_Q}
  Let $g(Z)\in\QQQ(Z)$ be a functionally indecomposable rational
  function with $\abs{g^{-1}(\infty)}=1$ and of degree $n\ge2$.  Let
  $A$ and $G$ be the arithmetic and geometric monodromy group of $g$,
  respectively.  Let $T$ be the ramification type of $g$. Then one of
  the following holds:
\begin{enumerate}
\item $n$ is a prime, $C_n=G\le A=\AGL_1(n)$, $T=(n,n)$.
\item $n\ge3$ is a prime, $D_n=G\le A=\AGL_1(n)$, $T=(2,2,n)$.
\item $n\ge4$, $\alt_n\le G\le A\le\sym_n$, many possibilities for
  $T$.
\item $n=6$, $G=\PGL_2(5)=A$, $T=(2,4,6)$.
\item $n=9$, $G=\PgL_2(8)=A$, $T=(3,3,9)$.
\item $n=10$, $G=\PgL_2(9)=A$, $T=(2,4,10)$.
\end{enumerate}
\end{Theorem}

\subsection{Proof of Theorem \ref{T:g0_Q}}

  Let $\cE=(\sigma_1,\sigma_2,\dots,\sigma_r)$ be a genus $0$ system
  of $G$, and $T$ its type, such that $\sigma_r$ is the element
  $\sigma$ from Lemma \ref{L:Siegel:rat}. So $n=2m$, where $\sigma_r$
  has two cycles, both of length $m$.
  
  We denote by $L$ a splitting field of $g(Z)-t$ over $\QQQ(t)$, and
  if $U$ is a subgroup of $A=\Gal(L/\QQQ(t))$, then $L_U$ is the fixed
  field of $U$ in $L$.
  
  First suppose that $A$ is an affine permutation group (different
  from $\alt_4$ and $\sym_4$). Theorems \ref{T:kl} and \ref{T:g0}
  gives the candidates for $G$ and $A$ and genus $0$ systems. The only
  possible degrees are $8$ and $16$.
  
  Suppose $n=8$. The only possible candidate with a genus $0$ system
  is $G=\AGL_3(2)=A$. The rational genus $0$ systems in $G$ have type
  $(3,4,4)$, $(4,4,4)$, $(2,2,4,4)$, $(2,2,3,4)$, or $(2,2,2,2,4)$.
  
  The $(3,4,4)$--tuple must have all branch points rational. By
  \cite{Malle:FoD}, the minimal field of definition of such a cover
  has degree $2$ over $\QQQ$, so this case is out.
  
  In the $(4,4,4)$ case, a minimal field of definition has degree $4$
  over $\QQQ$ if all branch points are rational. There could possibly
  be two of the branch points conjugate, which would lower the degree
  of the minimal field of definition by at most a factor $2$, so this
  does not occur as well.
  
  The cases with $4$ and $5$ branch points all occur, see Section
  \ref{SS:Comp}.

Now suppose $n=16$. The only cases where $G$ has a genus $0$ system of
the required form, and $\sigma_r$ fulfills the necessary properties in
Lemma \ref{L:Siegel:rat}, are the following:
\begin{itemize}
\item[(a)] $G$ has index $2$ in $(\sym_4\times\sym_4)\rtimes C_2$,
  $T=(2,4,8)$.
\item[(b)] $G=A=(\sym_4\times\sym_4)\rtimes C_2$, $T=(2,6,8)$ or
  $T=(2,2,2,12)$. (This is case $m=4$ in \ref{T:g0:II} of Theorem
  \ref{T:g0}.)
\item[(c)] $G=A=C_2^4\rtimes\sym_5$, $T=(2,5,8)$, $(2,6,8)$, and
  $(2,2,2,8)$.
\end{itemize}

We start excluding case (a), where $G_1=(C_3\times C_3)\rtimes
C_4$, and $\cE$ has type $(2,4,8)$. Here $[A:G]\le2$. The group $G$
has, up to conjugacy, a unique subgroup $U$ of index $8$. Set $\tilde
U:=N_A(U)$. Then $A=\tilde UG$, so the fixed field $L_{\tilde U}$ is a
regular extension of $\QQQ(t)$. Look at the action of $A$ on $A/\tilde
U$. With respect to this action, the elements in $\cE$ have cycle
types $2-2$, $2-2-4$, $8$. From that we get that $L_{\tilde U}$ has
genus $0$, and because of the totally ramified place at infinity, we
have $L_{\tilde U}=\QQQ(x)$ where $t=f(x)$ with $f\in\QQQ[X]$. Now
$A$, in this degree $8$ action, preserves a block system of blocks of
size $4$, and the last element in $\cE$ leaves the two blocks
invariant. Suppose without loss that $\sigma_2$ corresponds to $0$.
Then this yields (after linear fractional changes) $f(X)=h(X)^2$ with
$h\in\QQQ[X]$, where $h(X)=X^2(X^2+pX+p)$, where the ramification
information tells us that $h$ has, besides $0$, two further branch
points which are additive inverses to each other. This gives the
condition $27p^2-144p+128=0$, so $p\in\QQQ(\sqrt{3})\setminus\QQQ$, a
contradiction.

Cases (b) and (c) however have the required arithmetic realizations.
As the proof involves a considerable amount of computations, we
postpone the analysis to Section \ref{SS:Comp}.

None of the product action cases in \ref{T:g0:II} with $m\ge5$ can
occur, because by \ref{T:kl} there is no element with two cycles of
equal lengths.

Now assume that $A$ is an almost simple group. Suppose that $A$ is
neither the alternating nor the symmetric group in natural action.
Theorem \ref{T:g0} lists those cases where a transitive normal
subgroup $G$ has a genus $0$ system. In our case, the permutation
degree $n=2m$ is even, and one member $\sigma_r$ of the genus $0$
system is a product of two $m$--cycles. The condition (b) in Lemma
\ref{L:Siegel:rat}, namely that $\sigma_r$ is rational in $A$, already
excludes most examples. The two biggest degrees which survive that
condition are $n=22$ with $G=\M22$, $A=\M22\rtimes C_2$ and $n=24$
with $G=A=\M24$. However, $\sigma_r$ violates condition (d) of Lemma
\ref{L:Siegel:rat} in both cases.

Excluding the case $n=12$, $G=\M12$ for a moment, the next smallest
cases with rational $\sigma_r$ have degree $n\le10$. We go through the
possibilities which fulfill the necessary properties from Lemma
\ref{L:Siegel:rat}, starting with the small degrees.

Let $n=6$. Then $A=\PGL_2(5)$, and $G=\PSL_2(5)$ or $G=\PGL_2(5)$. If
$G=A$, then $T=(4,4,3)$, and an example is given by
\[
g(Z)=\frac{Z^4(13Z^2-108Z+225)}{(Z^2-15)^3}.
\]
Next suppose $G=\PSL_2(5)$. There is the possibility $T=(2,5,3)$, with
an example
\[
g(Z)=\frac{Z^5(Z-2)}{(Z^2-5)^3},
\]
or $T=(2,2,2,3)$, where
\[
g(Z)=\frac{(Z^2-2Z+2)(Z^2-16Z+14)^2}{(Z^2-2)^3}
\]
is an example.

Let $n=8$. Then $A=\PGL_2(7)$, and $[A:G]\le2$. First suppose
$G=\PSL_2(7)$. Then $T=(3,3,4)$. Suppose the required $g(Z)$ exists.
Without loss assume that $\infty$ is the branch point corresponding to
$\sigma_3$. The two finite branch points could be algebraically
conjugate. But there is a Galois extension $K/\QQQ$ of degree dividing
$4$, such that the branch points are in $K$, and
$g^{-1}(\infty)\subset K$. So, by linear fractional twists over $K$,
we can pass from $g$ to
\[
\tilde g(Z)=\frac{(Z^2+a_1Z+a_0)(Z^2+p_1Z+p_0)^3}{Z^4}.
\]
If $a_1\ne0$, then we may assume that $a_1=1$. If however $a_1=0$,
then $p_1=0$ cannot hold, because $\tilde g$ were functionally
decomposable. Thus if $a_1=0$, we may assume that $p_1=1$. Thus we
have two cases to consider. Together with the obvious requirement
$a_0p_0\ne0$, and the ramification information in the other finite
branch point, this gives a $0$--dimensional quasi affine variety. See
\cite[Sect.~I.9]{MM} where this kind of computation is explained in
detail. By computing a Gr\"obner bases with respect to the
lexicographical order we can solve the system. We obtain an empty set
in the second case, and a degree $4$ equation over $\QQQ$ for $p_1$ in
the first case. However, this degree $4$ polynomial turns out to be
irreducible over $\QQQ$ with Galois group $D_4$, hence $p_1\not\in K$,
a contradiction.

Now assume $G=A$. Then $T=(2,6,4)$. The corresponding triple is
rationally rigid and $\sigma_2$ has a single cycle of length $6$, so
there exists a rational function $g(Z)$ with the required ramification
data. Still, we need to decide about the fiber $g^{-1}(\infty)$. We do
this by explicitly computing $g$, getting
$g(Z)=\frac{Z^6(9Z^2-6Z+49)}{(Z^2+7)^4}$. So the fiber
$g^{-1}(\infty)$ is not real, contrary to our requirement.

Let $n=10$. Then $\Si\le A\le\Aut(\Si)$ with $\Si=\alt_5$ or
$\Si=\alt_6$. In view of the results we want to achieve, there is
little interest in investigating these cases more closely.

Finally, we have to rule out the case $n=12$, $G=A=\M12$. We have the
following possibilities for $T$: $(2,5,6)$, $(3,4,6)$, $(3,3,6)$,
$(4,4,6)$, $(2,6,6)$, $(2,8,6)$, and $(2,2,2,6)$.

In the cases with three branch points, explicit computations are
feasible, and it turns out that only the two cases $(3,3,6)$ and
$(4,4,6)$ give Galois realizations over $\QQQ(t)$. However, in both
cases the subfields of degree $12$ over $\QQQ(t)$ are not rational.
Indeed, in the first case, we get the function field of the quadratic
$X^2+Y^2+1=0$, and in the second case, the function field of the
quadratic $X^2+3Y^2+5=0$. In Section \ref{SS:Comp} we give explicit
polynomials over $\QQQ(t)$ of degree $12$ with Galois group $\M12$ and
ramification type $(3,3,6)$ or $(4,4,6)$, respectively. However, a
variation of the argument below could be used as an alternative.

So we need to worry about the ramification type $T=(2,2,2,6)$. The
criterion in Lemma \ref{L:Siegel:rat} is too coarse in order to rule
out that case.  However, we still get rid of this case by considering
the action of complex conjugation, and what it does to a genus $0$
system. Let $\cE$ be a genus $0$ system of type $T$, and suppose that
a function $g(Z)$ exists as required. By passing to a real field $k$
containing $g^{-1}(\infty)$, we may assume that $g(Z)=h(Z)/Z^6$, where
$h[Z]\in k[Z]$ is a monic polynomial of degree $12$ and $h(0)\ne0$. If
$h(0)<0$, then $h(Z)-t_0Z^6$ has exactly $2$ real roots for $t_0\ll0$
(by a straightforward exercise in calculus). However, $\M12$ does not
have an involution with only $2$ fixed points, so this case cannot
occur.

Thus $h(0)>0$. Then, for $t_0\gg0$, $h(Z)-t_0Z^6$ has precisely $4$
real roots. Choose such a $t_0\in k$ with $\Gal(h(Z)-t_0Z^6/k)=\M12$.
By a linear fractional change over $k$, we can arrange the following:
$t_0$ is mapped to $\tilde t_0$, the branch points of the
corresponding rational function $\tilde g$ are all finite, and the
real branch points of $\tilde g$ are smaller than $\tilde t_0$. Let
$\tilde t_0$ be the base point of a branch cycle description
$\mathbf{\sigma}=(\sigma_1,\sigma_2,\sigma_3,\sigma_4)$ coming from
the ``standard configuration'' as in \cite[Sect.~I.1.1]{MM} or
\cite[\S2]{FD:real}. Note that the order of the conjugacy classes here
must \emph{not} be chosen arbitrarily. So the element of order $6$ is
one of the $\sigma_i$. As $k\subset\RRR$, complex conjugation $\rho$
leaves the set of branch points invariant, but reflects the paths at
the real axis, inducing a new branch cycle description
$\mathbf{\sigma}^\rho$. For instance, if all branch points are real,
we get
\[
\mathbf{\sigma}^\rho=(\sigma_1^{-1},(\sigma_2^{-1})^{\sigma_1^{-1}},
(\sigma_3^{-1})^{\sigma_2^{-1}\sigma_1^{-1}},
(\sigma_4^{-1})^{\sigma_3^{-1}\sigma_2^{-1}\sigma_1^{-1}}),
\]
and a similar transformation formula holds if there is a pair of
complex conjugate branch points. For this old result by Hurwitz
, see \cite[Theorem I.1.2]{MM}, \cite{FD:real}.
  
Identify the Galois group $\Gal(\tilde g(Z)-\tilde t/k(\tilde t))$
with $\Gal(\tilde g(Z)-\tilde t_0/k)$, so that they are permutation
equivalent on the roots of $\tilde g(Z)-\tilde t$ and $\tilde
g(Z)-\tilde t_0$, respectively. Let $\psi$ be the complex conjugation
on the splitting field of $\tilde g(Z)-\tilde t_0$.  Then, under this
identification, $\mathbf{\sigma}^\psi=\mathbf{\sigma}^\rho$.  (Here
$\mathbf{\sigma}^\psi$ means simultaneously conjugating the components
with $\psi$.)  This is a result by D\`ebes and Fried, extending a more
special result by Serre \cite[8.4.3]{Serre:Galois} (which does not
apply here), see \cite{FD:real} and \cite[Theorem I.10.3]{MM}.
  
Now, for instance using GAP, one checks that in all possible
configurations for $\mathbf{\sigma}$ and possibilities of real and
complex branch points, an element $\psi$ as above either does not
exist, or is a fixed point free involution. However, as we have chosen
$\tilde t_0$ such that $\tilde g(Z)-\tilde t_0$ has precisely $4$ real
roots, the case that $\psi$ has precisely $4$ fixed points should also
occur.  As this is not the case, we have ruled out the existence of
$\M12$ with this specific arithmetic data.

\subsection{Computations}\label{SS:Comp}

This section completes the proof of Theorem \ref{T:g0_Q} in those
cases which require or deserve some explicit computations besides
theoretical arguments. We continue to use the notation from there.

\subsubsection{$n=8$, $G=\AGL_3(2)$.}
Here $n=8$, and $G=A=\AGL_3(2)$. We have already seen that the only
possible ramification types could be $T=(2,2,2,2,4)$, $(2,2,4,4)$, and
$(2,2,3,4)$. We will establish examples for all three cases. While
deriving possible forms of $g(Z)$ we do not give complete
justification for each step, because the required properties of $g(Z)$
can be verified directly from the explicit form. Thus the description
of the computation is only meant as a hint to the reader how we got the
examples.

In the construction of examples we employ a $2$--parametric family of
polynomials of degree $7$ over $\QQQ(t)$ with a $(2,2,2,2,4)$
ramification type and Galois group $\PSL_2(7)$. This family is due to
Malle, see \cite{Malle:PSL2_7}. Define
\begin{multline*}
f_{\alpha,\beta}(X) := \frac{(X^3+2(\beta-1)X^2+
(\alpha+\beta^2-4\beta)X-2\alpha)}{X^2(X-2)}\cdot\\
 (X^4-2(\beta+2)X^2+4\beta X-\alpha).
\end{multline*}
One verifies that for all $(\alpha,\beta)\in\QQQ^2$ in a non--trivial
Zariski open subset of $\QQQ^2$, the following holds:
$f_{\alpha,\beta}$ has arithmetic and geometric monodromy group
$\PSL_2(7)$ with ramification type $(2,2,2,2,4)$. The elements of
order $2$ are double transpositions, while the element of order $4$
has type $1-2-4$. We take the composition $f_{\alpha,\beta}(r(X))$,
where $r\in\QQQ(X)$ has degree $2$, and is ramified in $0$ and $1$.
Multiplying $r$ with a suitable constant (depending on $\alpha$ and
$\beta$), one can arrange that the discriminant of the numerator of
$f_{\alpha,\beta}(r(X))-t$ is a square. This can be used to show that
the arithmetic and geometric monodromy group of
$f_{\alpha,\beta}(r(X))$ is $\AGL_3(2)$ in the degree $14$ action. One
can now pass to the fixed field $E$ of $\GL_3(2)<\AGL_3(2)$ in a
splitting field $L$ of $f_{\alpha,\beta}(r(X))-t$ over $\QQQ(t)$. A
minimal polynomial $F_{\alpha,\beta}(X,t)$ for a primitive element of
$E/\QQQ(t)$ can be computed, we do not print it here because it is very
long. For that we used a program written by Cuntz
based on KASH\cite{kash} which computes subfields in algebraic
function fields.

It turns out that the degree in $t$ of $F_{\alpha,\beta}(X,t)$ is $2$.
So we can easily derive a condition for the genus $0$ field $E$ to be
rational. In this case, we get that $E$ is rational if and only if
$-\alpha$ is a sum of two squares in $\QQQ$. For instance, the choice
$\alpha:=-1/2$, $\beta=1$ yields
\[
g(Z)=\frac{(13Z^4+60Z^3+100Z^2+72Z+20)(11Z^4+8Z^3-12Z^2-16Z+12)}{(Z^2-2)^4}.
\]

Next we want to see how to get the cases with $4$ branch points. Let
$\Delta_{\alpha,\beta}(t)$ be the discriminant of a numerator of
$f_{\alpha,\beta}-t$ with respect to $X$. A necessary condition for
having only $4$ branch points is that the discriminant of
$\Delta_{\alpha,\beta}(t)$ with respect to $t$ vanishes. This gives a
condition on $\alpha$ and $\beta$, and if one performs the
computation, it follows that this condition is given by the union of
two genus $0$ curves which are birationally isomorphic to $\PP(\QQQ)$
over $\QQQ$. For the computation of such a birational map, we made use
of the Maple package \texttt{algcurves} by Mark van
Hoeij (available at
http://klein.math.fsu.edu/{$\tilde{\;}$}hoeij, also implemented in
Maple V Release 5).

An example for the ramification type $(2,2,4,4)$ is
\[
g(Z)=\frac{(3Z^2-15Z+20)Z^2}{(Z^2-5)^4},
\]
whereas
\[
g(Z)=\frac{(11Z^2+30Z+18)(3Z^2+30Z-46)^3}{(Z^2-2)^4}
\]
is an example of ramification type $(2,2,3,4)$.

\subsubsection{$n=16$, $G=(S_4\times S_4)\rtimes C_2$}
Here $n=16$, and $G=A=(S_4\rtimes S_4)\rtimes C_2$ in product action
of the wreath product $S_4\wr C_2$. First suppose that $\cE$ has type
$(2,6,8)$. There are two such possibilities, both being rationally
rigid. The first has fine type $(2-2-2-2,3-6-6,8-8)$, and the second
one has fine type $(2-2-2-2-2-2,2-3-3-6,8-8)$. From this we can
already read off that there is a rational function $g(Z)\in\QQQ(Z)$ of
degree $16$ and the ramification data and monodromy groups given as
above. Let $\sigma_3$ correspond to the place at infinity. One
verifies that the centralizer $C_A(\sigma_3)$ is intransitive, so
$g^{-1}(\infty)\subset K\cup\{\infty\}$, where $K$ is a quadratic
subfield of $\QQQ(\zeta_8)$, so $K=\QQQ(\sqrt{-1})$,
$K=\QQQ(\sqrt{-2})$, or $K=\QQQ(\sqrt{2})$. The first two
possibilities cannot hold, because complex conjugation would yield an
involution in $A$, which inverts $\sigma_3$, and interchanges the two
cycles of $\sigma_3$. One verifies that such an element does not
exist. Let $\tilde D$ be the normalizer in $A$ of $I:=\gen{\sigma_3}$.
Then $\tilde D$ contains a decomposition group $D$ of a place of $L$
lying above the infinite place of $\QQQ(t)$. Also, $[D:I]\ge4$ by
rationality of $\sigma_3$. On the other hand, $[\tilde D:I]=4$. Thus
$D=\tilde D$. But $\tilde D$ interchanges the two cycles of
$\sigma_3$, so the elements in $g^{-1}(\infty)$ cannot be rational.
This establishes the existence of $g$ of the required type.

In this situation, we were lucky that theoretical arguments gave a
positive existence result. However, it is also quite amusing to take
advantage of the specific form of $A$ and compute an explicit example
from the data given here.

Recall that $G=A=S_4\wr C_2$ is in product action. To this wreath
product there belongs a subgroup $U$ of index $8$, which is a point
stabilizer corresponding to the natural imprimitive action of $A$. The
fine types of the two $(2,6,8)$--tuples with respect to this degree
$8$ action are $(2,2-6,8)$ and $(2-2-2-2,2-3,8)$, respectively. One
verifies immediately that $L_U$ is a rational field, indeed
$L_U=\QQQ(x)$, where $t=h(x)^2$ with $h\in\QQQ[X]$. The idea is to
compute this field, and then extract from that the degree $16$
extension we are looking for.

In the first case, we may assume $h$ of the form $h(X)=X^3(X-1)$,
whereas $h(X)=X^3(X-8)+216$ (note that $h(X)+216=(X-6)^2(X^4+4X+12)$)
in the second case.

We have $h(x)^2=t$. Set $y:=h(x)$, and let $x'$ be a root of
$h(X)=-y$. Then also $h(x')^2=t$. However, $x+x'$ is fixed under a
suitable point stabilizer of $A$ with respect to the degree $16$
action of the wreath product $S_4\wr C_2$ in power action.

Take the first possibility for $h$. Using resultants, one immediately
computes a minimal polynomial $H(W,t)$ of $w:=x+x'$ over $\QQQ(t)$:
\begin{eqnarray*}
  H(W,t) &=&
  W^{16}-8W^{15}+27W^{14}-50W^{13}+55W^{12}-36W^{11}+13W^{10}\\
  & &
  -2W^9+136tW^8-544tW^7+892tW^6-744tW^5+315tW^4\\
  & & -54tW^3+16t^2.
\end{eqnarray*}
Here, however, $t$ appears quadratic, so this does not immediately
yield the function $g$ we are looking for. However, it is easy to
write down a parametrization for the curve $H(W,t)=0$:
\begin{align*}
W &= \frac{Z(2Z+3)}{Z^2-2}\\
t &= -\frac{1}{16}\frac{Z^6(Z+2)^6(2Z+3)^3}{(Z^2-2)^8}=:g(Z).
\end{align*}
The function $g(Z)$ parameterizing $t$ is the function we are looking
for.

Similarly, the second possibility of $h$ gives a function
\[
g(Z)=\frac{(Z^2+4Z+6)(Z-2)^2(3Z^2-4Z+2)^3}{(Z^2-2)^8}.
\]

By Theorem \ref{T:g0}, there is, for this setup, also the possibility
of a $(2,2,2,8)$ system. This is no longer rigid. But even if we could
show, for instance using a braid rigidity criterion as in
\cite[Chapt.~III]{MM}, the existence of a regular Galois extension
$L/\QQQ(t)$ with the correct data, we would not be able to decide
about rationality of the degree $16$ subfield we are after. However,
the following computations will display and solve the problem.

With $s\in\QQQ$ arbitrary set $h(X):=X^4+2sX^2+(8s+32)X+s^2-4s-24$.
One verifies that the splitting field of $h(X)^2-t$ over $\QQQ(t)$ is
regular with Galois group $A$, and that we have the ramification given
by the $(2,2,2,8)$ system, provided that $s\not\in{-4,-3,-12}$. (The
cases $s=-3$ and $s=-12$ give the first and second possibilities from
above, whereas for $s=-4$ the monodromy group of $h$ is $D_4$ rather
than $S_4$.) Again, let $x$ be with $h(x)^2=t$, and $x'$ be with
$h(x')=-h(x)$. As above, derive a minimal polynomial $H(W,t)$ for
$x+x'$ over $\QQQ(t)$. One calculates that the curve $H(W,t)=0$ is
birationally isomorphic to the quadratic $U^2-2V^2=4s+16$. Of course,
it depends on $s$ whether this quadratic has a rational point, which
in turn is equivalent that $L_U$ ($U$ from above) is a rational
function field. But if one chooses $s$ such that $4s+16=u_0^2-2v_0^2$
for $u_0,v_0\in\QQQ$, then $L_U$ is rational, and from the explicit
choice of a rational point on the quadratic we get $g(Z)$,
parameterized by $(u_0,v_0)$, where two such pairs give the same
function if $u_0^2-2v_0^2={u_0'}^2-2{v_0'}^2$. Up to the details which
are routine, this shows that the ramification type $(2,2,2,8)$ appears
as well.

\subsubsection{$n=16$, $G=C_2^4\rtimes S_5$}
Now $G=A=C_2^4\rtimes S_5$, where the action of $S_5$ is on the
$S_5$--invariant hyperplane of the natural permutation module for
$S_5$ over $\FFF_2$. We verify that the genus $0$ systems of type
$(2,5,8)$ and $(2,6,8)$ are rationally rigid, also, it follows from
the ramification type, that the degree $16$ field we are looking for
is rational. As in the previous case, we can recognize the
decomposition group (belonging to the inertia group
$I:=\gen{\sigma_3}$) as the normalizer of $I$ in $A$, and from the
properties of $N_A(I)$ we can read off, exactly as in the previous
case, that $g(Z)=h(Z)/(Z^2-2)^8$ exists as required.

Explicit computation is different from the previous case. Suppose we
have the ramification type $(2,5,8)$. As an abstract group,
$A=V\rtimes S_5$, where $V<\FFF_2^5$ is the hyperplane of vectors with
coordinate sum $0$, and $S_5$ permutes the coordinates naturally. This
interpretation of $A$ as a subgroup of the wreath product $C_2\wr S_5$
gives an imprimitive faithful degree $10$ action of $A$. Let $U$ be
the corresponding subgroup of index $10$. One verifies that $L_U$ is
the root field of $h(X^2)-t$, where $h(Y)=(Y-1)^5/Y$. Let $y_i$ be the
roots of $h(Y)-t$, $i=1,\dots,5$, and for each $i$, let $x_i$ be a
square root of $y_i$. Set $w=x_1+x_2+\dots+x_5$. We compute a minimal
polynomial $H(W,t)$ for $w$. Namely consider
$H(W,t):=\prod(X+\epsilon_1x_1+\epsilon_2x_2+\dots+\epsilon_5x_5)$,
where the product is over $\epsilon_i\in\{-1,1\}$, such that the sum
of the entries for each occuring tuple is $0$. Obviously, $H(w,t)=0$,
and $H(W,t)\in\QQQ[W,t]$. As to the practical computation, we computed
the solutions of $h(X)-t$ in Laurent series in $1/t^{1/5}$ around the
place with inertia group order $5$. Eventually, after calculations
similar as above, we get
\[
g(Z)=\frac{(Z-1)(Z^2+Z-1)^5}{(Z^2-2)^8}.
\]

If the ramification type is $(2,6,8)$, then $L$ is the splitting field
of $h(X^2)-t$, with $h(Y)=(2Y^2-27)^2(Y^2-1)^3/Y^2$, and after similar
computations we get
\[
g(Z)=\frac{(5Z^2+4Z-10)(Z+2)^2(5Z^2-12Z+6)^3}{(Z^2-2)^8}.
\]

Also, the case $(2,2,2,8)$ is not hard to establish by the procedure
described above. An example (as part of a $1$--parameter family) is
\begin{multline*}
  g(Z)=\frac{(15Z^4-74Z^3+140Z^2-124Z+44)^2}
  {(Z^2-2)^8}\cdot\\
  (47Z^8-472Z^7+1912Z^6-4272Z^5+4840Z^4-1824Z^3-288Z^2-64Z-16).
\end{multline*}

\subsubsection{$n=12$, $G=\M12$}
In order to rule out the ramification types $T=(3,3,6)$ and $(4,4,6)$,
we computed explicitly polynomials $F(X,t)$ of degree $12$ over
$\QQQ(t)$, such that the splitting field $L$ has Galois group $\M12$
over $\QQQ(t)$, and the ramification type $T$. From the explicit form
of $F(X,t)$ we can read off that a degree $12$ extension $E$ in $L$ of
$\QQQ(t)$ cannot be a rational field. Nowadays such computations are
routine, so we just give the polynomials.

For $T=(3,3,6)$ we obtain
\begin{eqnarray*}
F(X,t) &=& X^{12}+396X^{10}+27192X^9+933174X^8+20101752X^7+\\
       & & (-2t+169737744)X^6+16330240872X^5+\\
       & & (8820t+538400028969)X^4+(92616t+8234002812376)X^3+\\
       & & (-3895314t+195276967064388)X^2+\\
       & & (-48378792t+3991355037576144)X+\\
       & & t^2+62267644t+30911476378259268,
\end{eqnarray*}
and for $T=(4,4,6)$ we get
\begin{eqnarray*}
F(X,t) &=& X^{12}+44088X^{10}+950400X^9+721955520X^8+\\
       & & 31696106112X^7+(2t+5460734649920)X^6+\\
       & & 393700011065856X^5+\\
       & & (-120180t+20231483772508800)X^4+\\
       & & (-2587680t+911284967252689920)X^3+\\
       & & (137561760t+21295725373309787136)X^2+\\
       & & (4418468352t+183784500436675461120)X+\\
       & & t^2+31440107840t+3033666001201482093568.
\end{eqnarray*}

As $t$ is quadratic in both cases, it is easy to compute a quadratic
$Q$ such that $E$ is the field of rational functions on $Q$. Then $E$
is rational if and only if $Q$ has a rational point. However, in both
cases there is not even a real point on $Q$. This actually indicates
that the argument we used to exclude $T=(2,2,2,6)$ might be applicable 
here as well. One can verify that this is indeed the case.

\section{Applications to Hilbert's Irreducibility
  Theorem}\label{S:App}

Immediate consequences of Theorems \ref{T:g0}, \ref{T:g0pol},
\ref{T:g0_Q}, and \ref{T:g0pol_Q} are

\begin{Theorem}\label{T:comp} Let $k$ be a field of characteristic $0$,
  and $g(Z)\in k(Z)$ a rational function with the Siegel property.
  Then each non-abelian composition factor of $\Gal(g(Z)-t/k(t))$ is
  isomorphic to one of the following groups: $\alt_j$ ($j\ge5$),
  $\PSL_2(7)$, $\PSL_2(8)$, $\PSL_2(11)$, $\PSL_2(13)$, $\PSL_3(3)$,
  $\PSL_3(4)$, $\PSL_4(3)$, $\PSL_5(2)$, $\PSL_6(2)$, $\M11$, $\M12$,
  $\M22$, $\M23$, $\M24$.
\end{Theorem}

\begin{Theorem}\label{T:comp_Q} Let $g(Z)\in\QQQ(Z)$ be a Siegel
  function over $\QQQ$. Then each non-abelian composition factor of
  $\Gal(g(Z)-t/\QQQ(t))$ is isomorphic to one of the following groups:
  $\alt_j$ ($j\ge5$), $\PSL_2(7)$, $\PSL_2(8)$.
\end{Theorem}

\begin{Theorem} Let $g(Z)\in\QQQ(Z)$ be a Siegel function over
  $\QQQ$. Assume that $A=\Gal(g(Z)-t/\QQQ(t))$ is a simple group. Then
  $A$ is isomorphic to an alternating group or $C_2$.
\end{Theorem}

In \cite{PM:hitzt} we showed that this latter theorem implies

\begin{Corollary}
  Let $f(t,X)\in\QQQ(t)[X]$ be irreducible with Galois group $G$,
  where $G$ is a simple group not isomorphic to an alternating group
  or $C_2$. Then $\Gal(f(\bar t,X)/\QQQ)=G$ for all but finitely many
  specializations $\bar t\in\ZZZ$.
\end{Corollary}

Similarly, Theorems \ref{T:comp} and \ref{T:comp_Q} have the following 
application to Hilbert's irreducibility theorem. See \cite{PM:hitzt},
where we also have results of this kind which do not rely on
group-theoretic classification results.

\begin{Corollary}
  Let $f(t,X)\in\QQQ(t)[X]$ be irreducible, and assume that the Galois
  group of $f(t,X)$ over $\QQQ(t)$ acts primitively on the roots of
  $f(t,X)$ and has a non-abelian composition factor which is not
  alternating and not isomorphic to $\PSL_2(7)$ or $\PSL_2(8)$. Then
  $f(\bar t,X)$ remains irreducible for all but finitely many $\bar
  t\in\ZZZ$.
\end{Corollary}

\begin{Corollary}
  Let $k$ be a finitely generated field of characteristic $0$, and $R$
  a finitely generated subring of $k$. Let $f(t,X)\in k(t)[X]$ be
  irreducible, and assume that the Galois group of $f(t,X)$ over
  $k(t)$ acts primitively on the roots of $f(t,X)$ and has a
  non-abelian composition factor which is not alternating and is not
  isomorphic to one of the following groups: $\PSL_2(7)$, $\PSL_2(8)$,
  $\PSL_2(11)$, $\PSL_2(13)$, $\PSL_3(3)$, $\PSL_3(4)$, $\PSL_4(3)$,
  $\PSL_5(2)$, $\PSL_6(2)$, $\M11$, $\M12$, $\M22$, $\M23$, $\M24$.
  Then $f(\bar t,X)$ remains irreducible for all but finitely many
  $\bar t\in R$.
\end{Corollary}

\noindent{\sc IWR, Universit\"at Heidelberg, Im Neuenheimer Feld 368,\\
69120 Heidelberg, Germany}\par
\noindent{\sl E-mail: }{\tt Peter.Mueller@iwr.uni-heidelberg.de}
\end{document}